\documentclass[]{article}

\RequirePackage{xcolor,graphicx,subfig}

\RequirePackage{amsmath,amsfonts,amssymb,amsthm,mathtools}
\RequirePackage{bm}

\RequirePackage[a4paper]{geometry}
\RequirePackage{siunitx}

\DeclareMathOperator*{\tr}{\operatorname{tr}}

\DeclareMathOperator*{\sign}{\operatorname{sign}}
\newcommand{\defeq}{\vcentcolon=}

\newtheorem{scheme}{Scheme}
\newtheorem{remark}{Remark}

\usepackage[normalem]{ulem}
\usepackage{enumitem}

\newcommand{\sib}[1]{[\si{#1}]}

\title{A mathematical model for thermal single-phase flow and reactive transport in fractured porous media}

\author{Alessio Fumagalli \and Anna Scotti}

\date{MOX - Dipartimento di Matematica ``F. Brioschi'', Politecnico di
Milano, via Bonardi 9, 20133 Milan, Italy.}

\begin{document}

\maketitle

\begin{abstract}
    In this paper we present a mathematical model and a numerical workflow for
    the simulation of a thermal single-phase flow with reactive transport in
    porous media, in the presence of fractures. The latter are thin regions which might behave as high or
    low permeability channels depending on their physical parameters, and are thus of
    paramount importance in underground flow problems.
    Chemical reactions may alter the local properties of the porous media as well as the fracture walls,
    changing the flow path and possibly occluding some portions of the fractures
    or zones in the porous media. To solve numerically the coupled problem we
    propose a temporal splitting scheme so that the equations describing each physical
    process are solved sequentially. Numerical tests shows the accuracy of the
    proposed model and the ability to capture complex phenomena, where one or
    multiple fractures are present.
\end{abstract}






\section{Introduction}

The presence of fractures has an impact on subsurface flows at all scales: flow
tends to focus along highly permeable fractures, which can create shortcuts in
the domain, or, in the case of cemented or low permeable fractures, they might
create barriers in the domain.
In the context of reactive transport fractures can be
responsible for fast transport of fluid with different chemical composition with
respect to the surrounding matrix: this occurs for instance in geothermal
reservoirs where water with different salinity, solutes and temperature is
injected in the subsurface. These differences in composition and temperature can
trigger transformations such as mineral precipitation, dissolution or
replacement, with an impact on porosity and fracture aperture. The effective
exploitation of the geothermal system can be jeopardized by such phenomena.

Because of their thickness or aperture, fractures are usually represented as
lower dimensional objects and new equations along with interface conditions with
the surrounding porous media are derived. This procedure is usually referred to as
model reduction and the resulting model is named mixed-dimensional or
hybrid-dimensional problem. Seminal works dealing with single-phase flow are for
example \cite{Alboin2000,Alboin2002,Faille2002,Angot2003,Martin2005}. During the
years new models have been developed based on this idea, in particular for
multi-phase flow \cite{Fumagalli2012d,Ahmed2017}, transport
\cite{Stefansson2018,Chave2019}, and faults flow \cite{Tunc2012,Faille2014a,Faille2014}.
The geometrical complexity of the fracture networks requires to handle in an
accurate way also the intersection between them, indeed the intersection may
have different physical parameters than the incident fractures. In this case new
models have been derived where the intersections is part of the problem, see for
example \cite{Formaggia2012,Boon2018,Schwenck2015}.
In the special case of high speed circulation of the liquid in the fractures,
the Darcy model may be not appropriate. Thus several authors proposed a new
model based on Forchheimer or even more advanced flow model. Refer to
\cite{Frih2008,Morales2010,Morales2012,Knabner2014,Ahmed2018}.
Finally, we refer to \cite{Berre2019b} for a more detailed review on different strategies
to handle the complex problem of fractured porous medium.

The numerical solution of these problems is challenging due to several aspects,
in fact the fracture networks may pose severe constraints in the grid
generation resulting in poor quality and too many elements. Since this work
is more focused on the modelling side, we refer to the main works that dealt with
different classes of numerical schemes: classical mixed finite elements \cite{Martin2005},
hybrid high-order \cite{Chave2018},
discontinuous Galerkin \cite{Antonietti2019}, mimetic finite differences
\cite{Antonietti}, extended finite elements
\cite{Formaggia2012,Fumagalli2012g,Schwenck2015,Flemisch2016},
virtual element method \cite{Fumagalli2016a,Fumagalli2017a}, and  references therein.
Important benchmark studies to validate the effectiveness of the numerical
schemes are \cite{Dreuzy2013,Flemisch2016a,Scialo2017,Berre2020}.
Finally, a unified approach for numerical frameworks to solve such problems is
presented in \cite{Nordbotten2018}.

The aim of our work is to propose a model to account explicitly for the presence
of fractures and their impact on the flow, temperature, transport and reactions.
 The equations describing flow and transport are thus a coupled system of mixed-dimensional
PDEs which will be approximated by means of lowest order mixed finite elements
or mixed virtual elements, depending on the geometrical complexity of the
computational grid. We will consider a simple model for mineral precipitation
and dissolution following the model presented, among others, in \cite{Agosti2016}. To avoid
the occurrence of negative concentrations and oscillations when the amount of
precipitate approaches zero we adopt an event detection/location strategy to
detect the discontinuity in the ODE describing the reaction part, which is, for
this reason, split from advection and diffusion by means of a first-order
operator splitting. Several numerical examples will show the validity of our
approach for increasing level of geometrical difficulty of the fracture network.

The paper is organized as follow. In Section \ref{sec:model_porous_media} we
introduce the mathematical model to describe fluid flow,
heat transport, and solute transport with chemical reactions in porous media. The latter are
particularized in Section \ref{sec:chemical_model}. The mixed-dimensional
problem to describe the physical processes in the fractures is discussed in
Section \ref{sec:model_fracture}. Section \ref{sec:discretization} presents the
discretization considered to approximate the models, in particular a
splitting scheme is detailed that allows for a sequential resolution of each physical
process involved in the simulation. In Section
\ref{sec:numerical_examples} we run different examples to show the validity and
accuracy of the proposed approach. Finally, Section \ref{sec:conclusion} is
devoted to the conclusions.

\section{Model in the porous media} \label{sec:model_porous_media}

In this section we describe the mathematical model for our problem. The
physical processes are described separately but coupled together by suitable
constitutive relations. We first focus on the {model in the porous matrix},
leaving for a subsequent section the introduction of the reduced model to include
the fracture effects. First, we introduce the Darcy flow in Subsection \ref{subsec:darcy_flow}
followed by the heat equation in Subsection \ref{subsec:heat_equation}. The
models for the solute and precipitate are presented in Subsection
\ref{subsec:solute_precipitate}, and the section concludes with the constitutive
relations in Subsection \ref{subsec:constitutive_relations}.

\begin{figure}[tb]
    \centering
    \resizebox{0.33\textwidth}{!}{\fontsize{1cm}{2cm}\selectfont
\begingroup%
  \makeatletter%
  \providecommand\color[2][]{%
    \errmessage{(Inkscape) Color is used for the text in Inkscape, but the package 'color.sty' is not loaded}%
    \renewcommand\color[2][]{}%
  }%
  \providecommand\transparent[1]{%
    \errmessage{(Inkscape) Transparency is used (non-zero) for the text in Inkscape, but the package 'transparent.sty' is not loaded}%
    \renewcommand\transparent[1]{}%
  }%
  \providecommand\rotatebox[2]{#2}%
  \newcommand*\fsize{\dimexpr\f@size pt\relax}%
  \newcommand*\lineheight[1]{\fontsize{\fsize}{#1\fsize}\selectfont}%
  \ifx\svgwidth\undefined%
    \setlength{\unitlength}{368.96932983bp}%
    \ifx\svgscale\undefined%
      \relax%
    \else%
      \setlength{\unitlength}{\unitlength * \real{\svgscale}}%
    \fi%
  \else%
    \setlength{\unitlength}{\svgwidth}%
  \fi%
  \global\let\svgwidth\undefined%
  \global\let\svgscale\undefined%
  \makeatother%
  \begin{picture}(1,0.83213788)%
    \lineheight{1}%
    \setlength\tabcolsep{0pt}%
    \put(0,0){\includegraphics[width=\unitlength,page=1]{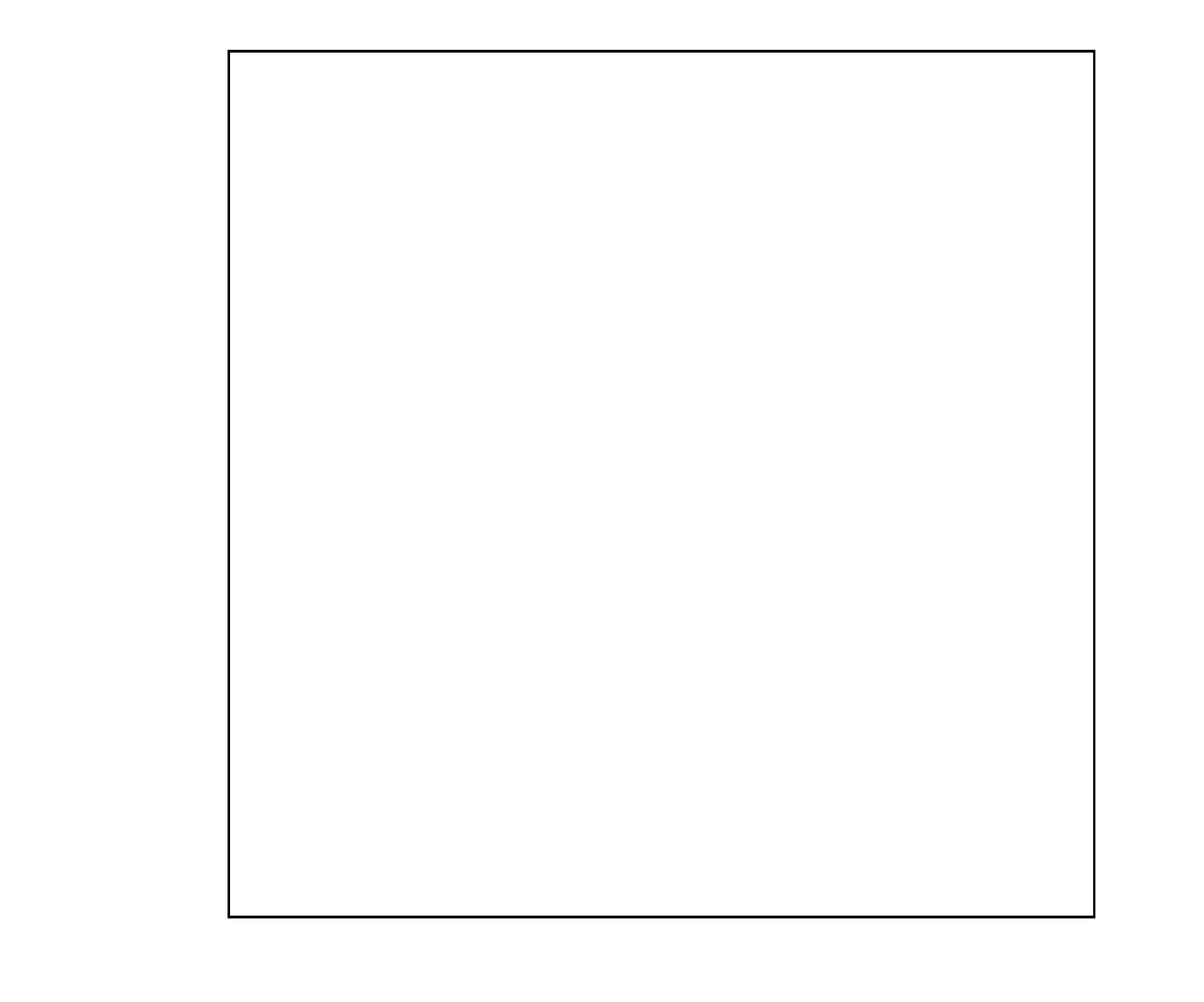}}%
    \put(0.21825759,0.71124567){\color[rgb]{0,0,0}\makebox(0,0)[lt]{\lineheight{1.25}\smash{\begin{tabular}[t]{l}$\Omega$\end{tabular}}}}%
    \put(0.50371328,0.01647731){\color[rgb]{0,0,0}\makebox(0,0)[lt]{\lineheight{1.25}\smash{\begin{tabular}[t]{l}$\partial_n \Omega$\end{tabular}}}}%
    \put(0.50371328,0.80969217){\color[rgb]{0,0,0}\makebox(0,0)[lt]{\lineheight{1.25}\smash{\begin{tabular}[t]{l}$\partial_n \Omega$\end{tabular}}}}%
    \put(0.06617176,0.41767952){\color[rgb]{0,0,0}\makebox(0,0)[lt]{\lineheight{1.25}\smash{\begin{tabular}[t]{l}$\partial_e \Omega$\end{tabular}}}}%
    \put(0.92585059,0.41767952){\color[rgb]{0,0,0}\makebox(0,0)[lt]{\lineheight{1.25}\smash{\begin{tabular}[t]{l}$\partial_e \Omega$\end{tabular}}}}%
    \put(0,0){\includegraphics[width=\unitlength,page=2]{domain.pdf}}%
    \put(0.0166007,0.68339009){\color[rgb]{0,0,0}\makebox(0,0)[lt]{\lineheight{1.25}\smash{\begin{tabular}[t]{l}$\bm{n}_{\partial \Omega}$\end{tabular}}}}%
  \end{picture}%
\endgroup%
}
    \caption{Example of porous media with some nomenclature considered.}%
    \label{fig:domain}
\end{figure}

The porous media occupies the domain $\Omega \subset \mathbb{R}^n$, with
$n=2$ or $3$, with Lipschitz continuous external boundary $\partial \Omega$. The
latter has been divided into two disjoint, possibly empty, subsets $\partial_e \Omega$ and
$\partial_n \Omega$, such that $\overline{\partial \Omega} =
\overline{\partial_e \Omega} \cup \overline{\partial_n \Omega}$ and
$\mathring{\partial_e \Omega} \cap \mathring{\partial_u \Omega} = \emptyset$.
For simplicity we assume that $\partial_n \Omega \neq \emptyset$.
The outward unit normal of $\partial \Omega$ is indicated as $\bm{n}_{\partial
\Omega}$, see Figure \ref{fig:domain}.
Finally, the final time is indicated as $T>0$.

It is important to note that, since the flow equation is the main driving force of the process,
the boundary conditions of the other problems {are conforming with} the ones imposed on the flow problem.

In the subsequent parts we will explicitly indicate the units of measure of each
variable and coefficient. We will make use of the notation 1\sib{\cubic\metre\of{\phi}} to be the
{unit cubic meter of pore space}, defined as \si{\cubic\metre\of{\phi}}=\si{\phi \cubic\metre} with $\phi$ the porosity.

General references for the following models are, for example, \cite{Bear1972,Bear1990,Helmig1997,Nordbotten2011}.

\subsection{Darcy flow model}\label{subsec:darcy_flow}

We consider a Darcy model to describe the flow of the water in the underground.
We are interested in the evolution of the Darcy velocity
$\bm{q}$ in \sib{\metre\per\second}
and pressure $p$ in \sib{\pascal} described by the system
\begin{subequations}\label{eq:problem_darcy}
\begin{align}\label{eq:darcy}
    &\begin{aligned}
        &\mu \bm{q} + k(\phi) \nabla p =\bm{0}\\
        &\partial_t \phi + \nabla \cdot \bm{q} + f =0
    \end{aligned}
    &\text{in } \Omega \times (0, T).
\end{align}
The scalar source or sink term is
denoted by $f$ in \sib{\per\second}.
To keep the model
simpler, we assume that the water
viscosity $\mu$ in \sib{\pascal\second} can be considered constant.  The permeability
$k$ in \sib{\square\metre} is a symmetric, isotropic and positive definite tensor which
depends on the porosity $\phi$ in \sib{\cubic\metre\of{\phi}\per\cubic\metre} which in turn, as we
will see in a subsequent model in Subsection
\ref{subsec:constitutive_relations}, depends on the precipitate concentration.
Thus also the porosity is a variable of the previous system. We assume
that $\phi \in [0, 1]$.

Boundary conditions are coupled to \eqref{eq:darcy} to close the system. In
particular, we have
\begin{gather}\label{eq:darcy_bc}
    \begin{aligned}
        &\tr \bm{q} \cdot \bm{n}_{\partial \Omega} = q_{\partial \Omega}&& \text{on } \partial_e
        \Omega\\
        &\tr p = p_{\partial \Omega} && \text{on }\partial_n \Omega
    \end{aligned},
\end{gather}
\end{subequations}
where $\tr$ denotes an abstract trace operator, $p_{\partial\Omega}$ in \sib{\pascal} and
$q_{\partial\Omega}$ in \sib{\meter\per\second} are the
pressure and normal flux given data. System \eqref{eq:problem_darcy} forms the
Darcy flow problem.

\subsection{Heat model} \label{subsec:heat_equation}

The heat equation models thermal conduction (Fourier's law) and convection of
 heat in the porous media.
A complete model should consider two sets of equations: one for the rock matrix and one for the
water, coupled with a suitable transfer function.
However we assume local thermal equilibrium, meaning that
the rock matrix and water are in thermal equilibrium so we can use only
one common set of primary variables to describe the process.
The temperature field is indicated as $\theta$ in \sib{\kelvin} and its evolution is described by
\begin{subequations}\label{eq:problem_heat}
\begin{align}\label{eq:heat}
    &\begin{aligned}
        &\bm{\tau} - \rho_w c_w \bm{q} \theta + \lambda(\phi) \nabla \theta = \bm{0}\\
        &\partial_t [c(\phi) \theta] + \nabla
        \cdot \bm{\tau} + j = 0
    \end{aligned}
    &\text{in } \Omega \times (0, T),
\end{align}
where $\bm{\tau}$ in \sib{\joule\per\square\metre\per\second} is the total heat flux,
$c$ in \sib{\joule\per\cubic\metre\per\kelvin}
is the effective thermal capacity
which is defined as the porosity weighted
average between the water $c_w$ and solid $c_s$ specific thermal capacity, both
in
\sib{\joule\per\kilogram\per\kelvin}. We have thus the expression of $c$ given
by
\begin{gather}\label{eq:effective_thermal_capacity}
    c(\phi) = \phi \rho_w c_w + (1-\phi) \rho_s c_s.
\end{gather}
$\rho_w$ and $\rho_s$ are the densities, both in
\sib{\kilogram\per\cubic\metre}, of the water and solid
phase respectively. Finally, $\lambda$ is the effective thermal conductivity
measured in \sib{\watt\per\metre\per\kelvin}. For simplicity, we assume that the
densities and the specific thermal capacities are given and constant.
Following \cite{Cote2005}, being the porous media saturated with water,
we model the effective thermal conductivity as
\begin{gather}\label{eq:effective_thermal_conductivity}
    \lambda(\phi) = \lambda_w^\phi \lambda_s^{1-\phi},
\end{gather}
where $\lambda_s$ and $\lambda_s$ are water and solid thermal conductivity, both
in \sib{\watt\per\metre\per\kelvin}, and for simplicity assumed to be constant.
Finally, $j$ in \sib{\joule\per\cubic\metre\per\second} models a source or sink of heat in the system.

In addition to system \eqref{eq:heat}, we consider suitable boundary and initial
conditions,
\begin{gather}\label{eq:heat_bc}
    \begin{aligned}
        &\tr \bm{\tau} \cdot \bm{n}_{\partial \Omega} = \tau_{\partial \Omega}
        && \text{on }\partial_n \Omega \times (0, T) \\
        &\tr \theta = \theta_{\partial \Omega} && \text{on } \partial_e \Omega
        \times (0, T)\\
        &\theta(t=0) = \theta_0 && \text{in } \Omega \times \{0\}
    \end{aligned},
\end{gather}
\end{subequations}
where $\tau_{\partial\Omega}$ in \sib{\joule\per\square\metre\per\second}
and $\theta_{\partial\Omega}$ in \sib{\kelvin} are suitable boundary data for the
heat flux and temperature, respectively. Finally, $\theta_0$ in \sib{\kelvin} is the initial
condition for the temperature. System \eqref{eq:problem_heat} forms the heat
problem.

\subsection{Solute and precipitate model} \label{subsec:solute_precipitate}

We consider the passive scalar model to describe the evolution of the solute $u$
in \sib{\mole\per\cubic\metre\of{\phi}} in the porous medium. Note that solute concentration
is expressed in terms of moles per unit pore volume, so we should have $u \geq0$. The system is written as
\begin{subequations}\label{eq:problem_solute_precipitate}
\begin{align}\label{eq:solute}
    &\begin{aligned}
        &\bm{\chi} - \bm{q} u + \phi d \nabla u = \bm{0}\\
        &\partial_t (\phi u) + \nabla \cdot \bm{\chi} + \phi r_w(u, w; \theta)
        = 0
    \end{aligned}
    &\text{in } \Omega\times (0, T),
\end{align}
with $\bm{\chi}$ in \sib{\mole\per\square\metre\per\second} the total flux given by a
combination of the advective field and Fick's law. The symmetric, isotropic and positive definite tensor $d$ in \sib{\square\metre\per\second}
represents the molecular diffusivity of $u$ in
the water and $r_w$ in \sib{\mole\per\cubic\metre\of{\phi}\per\second}
is a reaction term which involves also the precipitate $w$, as well as the
temperature $\theta$ {and will be detailed later on}.

To model the concentration of the precipitate $w\geq0$, expressed in
\sib{\mole\per\cubic\metre\of{\phi}}, we consider the
following ordinary differential equation which models its evolution in time
\begin{align}\label{eq:precipitate}
    &\partial_t (\phi w) - \phi r_w(u, w; \theta) = 0 && \text{in } \Omega\times (0, T).
\end{align}
The actual expression of $r_w$ depends on several aspects, see Section \ref{sec:chemical_model} for a more
detailed discussion.

Boundary and initial conditions are supplied to equations
\eqref{eq:solute} and \eqref{eq:precipitate} as
\begin{gather}\label{eq:solute_precipitate_bc}
    \begin{aligned}
        &\tr \bm{\chi} \cdot \bm{n}_{\partial \Omega} = \chi_{\partial\Omega}
        &&\text{on } \partial_e \Omega \times (0, T)\\
        &\tr u = u_{\partial\Omega} &&\text{on } \partial_n \Omega \times (0, T)\\
        &u(t=0) = u_0 && \text{in } \Omega \times \{0\}\\
        &w(t=0) = w_0 &&\text{in } \Omega \times \{0\}
    \end{aligned},
\end{gather}
\end{subequations}
where the values of $\chi_{\partial\Omega}$ in \sib{\mole\per\square\metre\per\second} and
$u_{\partial\Omega}$ in \sib{\mole\per\cubic\metre\of{\phi}} are the normal component of
total flux and solute data imposed at the boundary $\partial \Omega$. Finally,
$u_0$ in \sib{\mole\per\cubic\metre\of{\phi}} and $w_0$ in
\sib{\mole\per\cubic\metre\of{\phi}} specify the initial conditions for both the solute and
precipitate. System \eqref{eq:problem_solute_precipitate} forms the solute and
precipitate problem.

\subsection{Permeability and porosity model}\label{subsec:constitutive_relations}

\begin{figure}[tbp]
    \centering
    \includegraphics[width=0.33\textwidth]{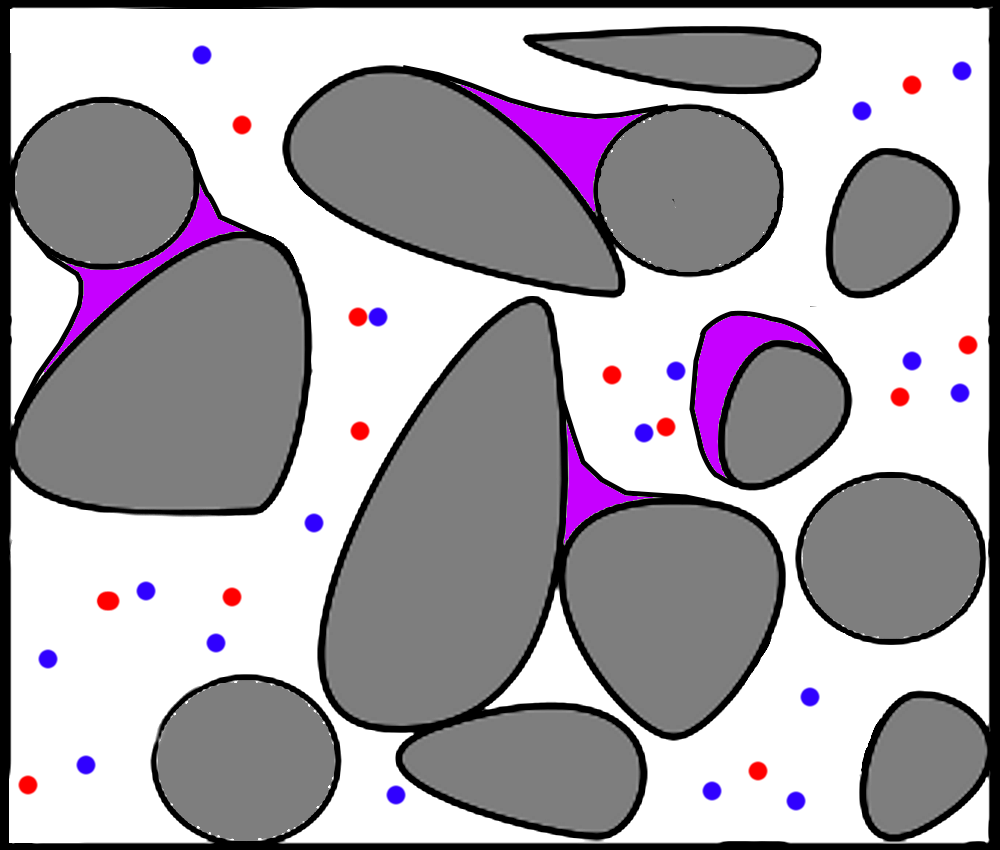}
    \caption{Graphical representation that shows the two solutes (blue and red
    dots) reacting and forming the precipitate (violet) around the
    grains.}
    \label{fig:deposition_grains}
\end{figure}

We consider a Kozeny-type relationship \cite{Bear1972} to link the permeability
with the porosity, namely
\begin{gather}\label{eq:kozeny}
    k(\phi) = k_0 \dfrac{\phi^2}{\phi_0^2}
\end{gather}
where ${k}_0$ in \sib{\square\metre} and $\phi_0$ in
\sib{\cubic\metre\of{\phi}\per\cubic\metre} are the given reference permeability
and porosity, respectively. Other and more
general relationships
between $k$ and $\phi$ are possible.
Finally, the model that links the porosity $\phi$ to the
precipitate $w$ is given by  the following ODE,
\begin{gather}\label{eq:porosity}
    \begin{aligned}
        &\partial_t \phi + \eta_\Omega \phi \partial_t w = 0 && \text{in } \Omega\times
        (0, T)\\
        &\phi(t=0) = \phi_0 && \text{in } \Omega \times \{0\}
    \end{aligned},
\end{gather}
where $\eta_\Omega$ in \sib{\cubic\metre\of{\phi}\per\mole} is a
proportionality parameter associated with molar volume of the mineral, see \cite{Noorden2009}, that determines the
rate of deposition of the solute around the grains as shown in Figure
\ref{fig:deposition_grains}.

\subsection{Complete model}\label{subsec:complete_model_bulk}

The model to describe flow, heat conduction and convection, precipitation and
dissolution of the chemical species as well as permeability and porosity
alteration is given by the system of equations \eqref{eq:problem_darcy},
\eqref{eq:problem_heat}, \eqref{eq:problem_solute_precipitate},
\eqref{eq:kozeny} and \eqref{eq:porosity}. The system is fully coupled,
non-linear and possibly non-smooth due to the discontinuous reaction rate.

\subsection{Non-dimensional equations}

We now derive the non-dimensional version of the solute equation to identify some
non-dimensional numbers that could characterize the solution behaviour.  Let us
define some reference quantities: $L$ for length, $Q$ for velocity,
$\phi_{\Omega, 0}$ for porosity, $u_{e}$ for molar concentration. Let us denote
with $\cdot^\prime$ nondimensional quantities, so that $x^\prime={x}/{L}$,
$t^\prime={tU}/{L}$, so that $\partial_t=\partial t^\prime {U}/{L}$ and $\partial
x={\partial x'}{1}/{L}$.
In primary form, the solute equation \eqref{eq:solute} can be rewritten as
\begin{gather*}
{\partial_{t^\prime} (u' \phi)} +  \nabla' \cdot \left( u'  \bm{q}^\prime\right) -
\dfrac{D }{LQ}  \nabla' \cdot \left( \phi  \nabla' u'\right) = -\dfrac{\phi
L}{Q} \lambda r'_w.
\end{gather*}
The above equation simplifies if we consider constant porosity, $\phi=\phi_0$.
\begin{gather*}
    \partial_{t^\prime} u' +  \dfrac{1}{\phi_0}\nabla' \cdot \left( u'
    \bm{q}^\prime\right) - \dfrac{D }{LQ}  \nabla' \cdot \left( \nabla'
    u'\right) = -\dfrac{L \lambda}{Q} r'_w.
\end{gather*}
Note that we can define a ``Reynolds number'' as
$\mathbb{R}e={LQ\phi_0}/{D}$ and a Damk\"ohler number $\mathbb{D}a=
{L \lambda \phi_0}/{Q}$.  A large $\mathbb{D}a$ corresponds to fast
reactions with respect to advection, while a small $\mathbb{D}a$ corresponds to
fast advection with respect to reaction speed.

\section{Chemical model}\label{sec:chemical_model}

Our idealized model for chemistry considers two reactions: precipitation and
dissolution, which can be written as
\begin{gather}
    \alpha U + \beta^+ V \rightarrow W + \beta^- V\label{eq:chemical_first}\\
    W + \beta^- V \rightarrow \alpha U +\beta^+ V\label{eq:chemical_second}
\end{gather}
where $U$ and $V$ are two solutes (ions) that can precipitate to form a solid
(salt) $W$, and $\alpha$ and $\beta^{\pm}$ are integer stoichiometric
coefficients.
According to the mass action law \cite{Guldberg1864}, the precipitation rate depends
on the rate $\lambda^+$, which is the reaction constant for
\eqref{eq:chemical_first}, and on the
concentrations of the two ions raised to the power indicated by the
stoichiometric coefficients $\alpha$ and $\beta^+$, conversely, the rate of
dissolution depends only on the reactant $v$ if $\beta^-$ is greater than zero.
The latter being the reaction constant for \eqref{eq:chemical_second}. In our
model, these two
coefficients might depend on the temperature $\theta$ of the system
$\lambda^\pm = \lambda^\pm(\theta)$.
The net rate of precipitation is thus a function of $u$ and $v$, and it is given by
\begin{gather*}
    r_w (u, v; \theta) =  \lambda^+(\theta) u^\alpha v^{\beta^+} -
    \lambda^-(\theta) v^{\beta^-}
\end{gather*}
where $u$, $v$ are the molar concentrations in \sib{\mole\per\cubic\metre} of
$U$ and $V$, respectively. If we consider this simplified set of reactions
\begin{gather*}
    \alpha U + \beta V \rightarrow W \\
    W \rightarrow \alpha U +\beta V
\end{gather*}
the net precipitation rate becomes $r_w(u, v; \theta) = \lambda^+(\theta) u^\alpha v^{\beta} -
\lambda^-(\theta)$,
i.e. we have a dissolution rate $\lambda^-$ that is independent from $u$ and $v$ since it is established in literature that
the activity of a pure crystalline solid is a constant \cite{Knabner1995}.
At equilibrium the net rate of precipitation is zero, yielding
\begin{gather*}
    \lambda^+(\theta) u_{e}^\alpha v_e^{\beta} = \lambda^-(\theta),
\end{gather*}
where $u_e$ and $v_e$ are the molar concentrations of $u$ and $v$ at
equilibrium, respectively.
The net precipitation rate can be rewritten as
\begin{gather*}
    r_w(u, v; \theta)=\lambda^-(\theta) \left( \dfrac{u^\alpha v^\beta}{u_e^\alpha v_e^\beta} -1\right).
\end{gather*}
We further assume that there is
electrical equilibrium in the system, i.e. the number of cations equals the number of anions,
$u=v$ and we can consider only one of the two variables to describe the reaction
rate. For example, in the special case of $\alpha=\beta=1$ the reaction becomes
$U+V\leftrightarrow W$ and the reaction rate can be written as
\begin{gather*}
    r_w(u; \theta) =\lambda^-(\theta) \left[ \left(\dfrac{u}{u_e}\right)^2 -1\right].
\end{gather*}
In our work we consider problems
with a reaction rate is given by a function $r$ which depends on one of the two
solute
\begin{gather*}
    r_w(u; \theta) =\lambda^-(\theta) [r(u) -1].
\end{gather*}
Finally, note that the rate of dissolution does not depend on the concentration
of the solid salt $W$, in other words it does not vanish when $w=0$. This
condition must be enforced in the model, so that dissolution stops {only}
when $w\leq0$ {and} the net precipitation rate is negative, i.e. when the
precipitate is no longer present but the solute concentration is such that we
should have dissolution. This can be summarized in the following more general expression
\begin{gather}\label{eq:chimica}
    r_w(u, w; \theta) =\lambda^-(\theta)
    \left\{\max[r(u)-1,0]+H(w)\min[r(u)-1,0]\right\},
\end{gather}
where $H$ is the Heaviside function. The latter expression of $r_w$ is the one
used in this work.

\section{Model in the fracture}\label{sec:model_fracture}

\begin{figure}[tb]
    \centering
    \resizebox{0.33\textwidth}{!}{\fontsize{1cm}{2cm}\selectfont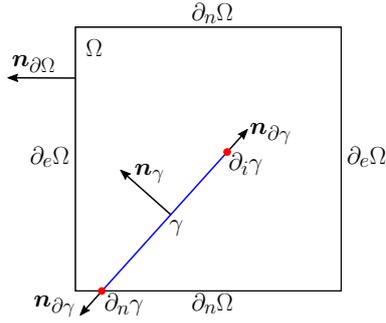}
    \caption{Example of fractured porous media
    with some nomenclature considered.}%
    \label{fig:fracture_domain}
\end{figure}

We here introduce a mixed-dimensional model to approximate the problem described in the previous sections in the
presence of fractures. We start with the simplified assumption of a single
fracture $\gamma$ cutting the domain as shown in Figure \ref{fig:fracture_domain}. The case of multiple intersecting fractures will be described later.
Since a fracture is an equi-dimensional region, i.e. 2D or 3D region with a small thickness, we
adopt the model reduction strategy to reduce the geometrical complexity and we
approximate it as a one co-dimensional
object. The fracture aperture $\epsilon_\gamma$ in \sib{\metre} will become a model
parameter and not a geometrical constraint,
which, in our case, is even more beneficial since it might change during the
simulation. For more references on this approach see
\cite{Martin2005,Morales2010,Frih2011,Jaffre2011,Tunc2012,Morales2012,Sandve2012,Fumagalli2012g,Dreuzy2013,Schwenck2015,Fumagalli2017a,Nordbotten2018,Stefansson2018,Chave2018,Ahmed2018,Berre2019b,Fumagalli2020a}
to name a few.

We assume that the fracture is open without any presence of infilling porous material,
getting unitary porosity in the fracture itself. So we get that
\si{\cubic\metre\of{\phi}}= \si{\cubic\metre} in the fracture, i.e. the pore volume and total volume are equivalent.
Now the role played by $\phi$
in the surrounding porous medium is given by the aperture $\epsilon_\gamma$.

The fracture $\gamma$ is a non self-intersecting piece-wise $C^2$ curve (if
$n=2$) or surface (if $n=3$). We indicate with $\partial \gamma$ the boundary of
$\gamma$, which can be divided into three disjoint, possible empty, parts:
$\partial_e \gamma$, $\partial_n\gamma$, and $\partial_i \gamma$. The latter is
the part of $\partial\gamma$ internal to the porous media. We clearly have
$\overline{\partial \gamma} = \overline{\partial_e \gamma} \cup
\overline{\partial_n \gamma} \cup \overline{\partial_i \gamma}$ and also that
$\mathring{\partial_* \gamma} \cap \mathring{\partial_{*^\prime} \gamma} =
\emptyset$ for any combination of non-equal elements in $*, *^\prime \in \{e, n,
i\}$. We can define an outward unit normal to $\partial \gamma$, which is tangent
to $\gamma$ itself and orthogonal to $\partial \gamma$ named $\bm{n}_{\partial
\gamma}$.

{At the interface with the surrounding medium,} the fracture has two different sides
$\gamma_+$ and $\gamma_-$ with associated normal vectors $\bm{n}_+$ and
$\bm{n}_-$.
We assign a {unique} normal $\bm{n}_\gamma$ to the fracture which is
associated with both fracture sides (i.e., $\bm{n}_\gamma =
\bm{n}_{\gamma,+}$ on $\gamma_+$ and {$\bm{n}_\gamma = -\bm{n}_{\gamma,-}$} on
$\gamma_-$). See Figure \ref{fig:fracture_domain} as an example.

\subsection{Reduced variables}

Variables and data associated with the fracture will be denoted with
a subscript $\gamma$, while we use a subscript $\Omega$ to indicate variables and
data in the surrounding porous media.
We introduce the fracture vector variables which come from the
integration over each normal section of the fracture of the corresponding equi-dimensional
variables, for $\bm{x} \in \gamma$ the curvilinear abscissa, as
\begin{gather*}
    \bm{q}_\gamma(\bm{x}) \defeq \int_{\epsilon_\gamma(\bm{x})} T(\bm{x}) \bm{q}(\bm{x}, s) d s \quad
    \bm{\tau}_\gamma(\bm{x}) \defeq \int_{\epsilon_\gamma(\bm{x})} T(\bm{x}) \bm{\tau}(\bm{x}, s) d s \quad
    \bm{\chi}_\gamma(\bm{x}) \defeq \int_{\epsilon_\gamma(\bm{x})} T(\bm{x})
    \bm{\chi}(\bm{x}, s) d s,
\end{gather*}
with $T\defeq I - N$ and $N \defeq \bm{n}_\gamma \otimes \bm{n}_\gamma$ the
tangential and normal projection
matrices, respectively.
The unit of measure of the previous variables are:
$\bm{q}_\gamma$ in \sib{\square\metre\per\second}, $\bm{\tau}_\gamma$ in
\sib{\joule\per\metre\per\second}, and
$\bm{\chi}_\gamma$ in \sib{\mole\per\metre\per\second}.
Moreover, the reduced scalar fields are defined
from their respective equi-dimensional variable as the average for each section,
\begin{gather*}
    p_\gamma(\bm{x}) \defeq \dfrac{1}{\epsilon_\gamma(\bm{x})}\int_{\epsilon_\gamma(\bm{x})}
    p(\bm{x}, s) ds \quad
    \theta_\gamma(\bm{x}) \defeq \dfrac{1}{\epsilon_\gamma(\bm{x})}
    \int_{\epsilon_\gamma(\bm{x})} \theta(\bm{x}, s) ds\\
    u_\gamma(\bm{x}) \defeq \dfrac{1}{\epsilon_\gamma(\bm{x})}\int_{\epsilon_\gamma(\bm{x})}
    u(\bm{x}, s) ds \quad
    w_\gamma(\bm{x}) \defeq \dfrac{1}{\epsilon_\gamma(\bm{x})}\int_{\epsilon_\gamma(\bm{x})}
    w(\bm{x}, s) ds.
\end{gather*}
In this case, the units of measures are the same as the original variables: $p_\gamma$ in \sib{\pascal}, $\theta_\gamma$ in \sib{\kelvin},
$u_\gamma$ in \sib{\mole\per\cubic\meter}, and $w_\gamma$ in \sib{\mole\per\cubic\meter}.
With an abuse in notation, we introduce the flux-based variable compounds as
\begin{gather*}
    \bm{q} \defeq (\bm{q}_\Omega, \bm{q}_\gamma)\quad
    \bm{\tau} \defeq (\bm{\tau}_\Omega, \bm{\tau}_\gamma)\quad
    \bm{\chi} \defeq (\bm{\chi}_\Omega, \bm{\chi}_\gamma)
\end{gather*}
even if that the units of measure of each compound are heterogeneous.

When the nabla operator is applied to a fracture variable, we implicitly assume that
it is defined on the tangential
space of the fracture itself, i.e. $\nabla\cdot\bm{\nu}_\gamma \defeq T :
\tilde{\nabla} \bm{\nu}_\gamma$,  $\tilde{\nabla}$ being the standard
gradient and $\bm{\nu}_\gamma$ a regular enough vector function defined on $\gamma$.
Analogously for the gradient of a fracture variable, which is defined as $\nabla
\nu_\gamma \defeq T \tilde{\nabla} \nu_\gamma$ with $\nu_\gamma$ a regular enough scalar function defined
on $\gamma$.   The conservation operators or mixed-dimensional divergences
are defined on the compounds and are given by
\begin{gather*}
    \nabla_\Omega \cdot \bm{\nu} \defeq \nabla \cdot \bm{\nu}_\Omega
    \quad \text{and} \quad
    \nabla_\gamma \cdot \bm{\nu} \defeq \nabla \cdot \bm{\nu}_\gamma - \tr
    \bm{\nu}_\Omega \cdot \bm{n}_\gamma
\end{gather*}
with $\bm{\nu} = (\bm{\nu}_\Omega, \bm{\nu}_\gamma)$ a generic compound of vector variables.
The $\nabla_\gamma \cdot$ considers also the contribution from the surrounding
porous media into the
fracture as flux exchange, and $\tr$ indicates the trace operator from $\Omega$ to each side of the
fracture, $\gamma_+$ and $\gamma_-$, viewed from $\Omega$.

Following the idea in
\cite{Martin2005}, we require that the coefficients associated to the diffusion
coefficients $k$ and $d$ can be decomposed as
\begin{gather}\label{eq:reduced_perm}
    k = \kappa N + k_\gamma T
    \quad\text{and}\quad
    d = \delta N + d_\gamma T,
\end{gather}
where the first relation
implies that the permeability $k$ can be decomposed in its normal
$\kappa$ (in \sib{\square\metre}) and tangential $k_\gamma$ (in \sib{\square\metre})
parts with respect to the fracture geometry. The same applies for  molecular
diffusivity, while the effective thermal conductivity \eqref{eq:effective_thermal_conductivity} is a scalar
value, so it does not need this decomposition.
Note that for some specific cases the value of the normal and tangential
component might be equal, however to keep a more general setting we  use different
symbols.

\subsection{Reduced Darcy flow model}

The reduced model for the Darcy flow describes the evolution of the reduced
Darcy velocity $\bm{q}_\gamma$ and pressure $p_\gamma$ in the fracture,
and reads
\begin{subequations}\label{eq:reduced_problem_darcy}
\begin{align}\label{eq:reduced_darcy}
    &\begin{aligned}
        &\mu \bm{q}_\gamma + \epsilon_\gamma k_\gamma(\epsilon_\gamma) \nabla p_\gamma =\bm{0}\\
        &\partial_t \epsilon_\gamma + \nabla_\gamma \cdot \bm{q} +
        \epsilon_\gamma f_\gamma =0
    \end{aligned}
    &\text{in } \gamma \times (0, T),
\end{align}
where the reduced source or sink term $f_\gamma$ \sib{\per\second} is computed as
$f_\gamma(\bm{x}) \defeq \epsilon_\gamma^{-1}(\bm{x})\int_{\epsilon_\gamma(\bm{x})}
f(\bm{x}, s) ds$. Following lubrication theory, the fracture tangential
permeability $k_\gamma$ is expressed as a function of the aperture, as described
in more detail in Subsection \ref{eq:reduced_permeability_aperture}.
Due to the introduction of the mixed-dimensional divergences
the form of \eqref{eq:reduced_darcy} is similar to \eqref{eq:darcy}.
At the fracture-matrix interface we consider a discrete version of Darcy's law in the normal direction, see for
example \cite{Martin2005}, which is given by
\begin{gather}\label{eq:reduced_darcy_cc}
    \mu \epsilon_\gamma \tr \bm{q}_\Omega \cdot \bm{n}_\gamma + \kappa_\gamma(\epsilon_\gamma) (p_\gamma - \tr p_\Omega) =
    0 \quad \text{on } \gamma \times (0, T).
\end{gather}
The latter relation models the fact that the flux exchange between the fracture
and the surrounding porous media is related to the pressure jump via
$\kappa_\gamma$.
This parameter, defined by \eqref{eq:permeability_aperture}, is
related to the aperture with a power law of exponent greater than one. Thus, if the
aperture goes to zero the flux exchange vanishes and the fracture and porous
media become decoupled.

Finally, we need to supply boundary conditions also to the fracture tips. In
particular, we distinguish between immersed tips and tips touching the domain boundary. In the first case, the so-called tip conditions are considered
while in the latter case we inherit the boundary conditions from the
equi-dimensional problem. In formula
\begin{gather}\label{eq:bc_reduced_darcy_model}
    \begin{aligned}
        &\tr \bm{q}_\gamma \cdot \bm{n}_{\partial \gamma} =
        q_{\partial\gamma}&&\text{on } \partial_e \gamma\times(0, T)\\
        &\tr p_\gamma = p_{\partial \gamma} && \text{on } \partial_n \gamma\times(0, T)\\
        &\tr \bm{q}_\gamma \cdot \bm{n}_{\partial \gamma} = 0
        &&\text{on } \partial_i \gamma\times(0, T)
    \end{aligned},
\end{gather}
\end{subequations}
with $q_{\partial \gamma}$ in \sib{\square\metre\per\second} and $p_{\partial
\gamma}$ in \sib{\pascal} are the given flux and
pressure at the fracture boundary. The last condition is the tip condition which
imposes no flow, see \cite{Angot2003}. System \eqref{eq:reduced_problem_darcy}
is the reduced Darcy flow
problem.

\subsection{Reduced heat model}

As we did in Subsection \ref{subsec:heat_equation},
we assume also in the fracture local thermal equilibrium. However, in a fracture we can reach high
speed circulation of water which might invalidate  this assumption. To keep
the presentation simple, we leave this case for future investigations.
Following the idea in \cite{Fumagalli2012a}, the heat equation which models the
thermal flux $\bm{\tau}_\gamma$ and temperature $\theta_\gamma$ in $\gamma$ is written as
\begin{subequations}\label{eq:reduced_problem_heat}
\begin{align}\label{eq:reduced_problem_heat_model}
    &\begin{aligned}
        &\bm{\tau}_\gamma - \rho_w c_w \bm{q}_\gamma \theta_\gamma + \epsilon_\gamma
        \lambda_w \nabla
        \theta_\gamma = \bm{0}\\
        &\rho_w c_w\partial_t (\epsilon_\gamma  \theta_\gamma) + \nabla_\gamma \cdot
        \bm{\tau} + j_\gamma = 0
    \end{aligned}
    &\text{in } \gamma \times (0, T).
\end{align}
Being the fracture open, the effective thermal capacity and conductivity are simplified
to $\lambda = \lambda_w$ as well as $c = \rho_w c_w$. The source term
$j_\gamma$ in \sib{\joule\per\cubic\metre\per\second} is computed as
$j_\gamma(\bm{x}) \defeq \epsilon_\gamma^{-1}(\bm{x})\int_{\epsilon_\gamma(\bm{x})}
j(\bm{x}, s) ds$. At the interface between
the fracture and surrounding porous media the coupling conditions are
\begin{align}\label{eq:reduced_problem_heat_cc}
    \epsilon_\gamma (\tr \bm{\tau}_\Omega \cdot \bm{n}_\gamma - \rho_w c_w \tr
    \bm{q}_\Omega \cdot \bm{n}_\gamma \tr \theta_\Omega) + \lambda_w
    (\theta_\gamma - \tr \theta_\Omega) = 0\quad \text{on } \gamma \times (0,
    T).
\end{align}
The boundary and initial conditions for the reduced heat equation are inherit
from the equi-dimensional problem and they are given by
\begin{gather}\label{eq:reduced_problem_heat_bc}
    \begin{aligned}
        &\tr \bm{\tau}_\gamma \cdot \bm{n}_{\partial \gamma} = \tau_{\partial
        \gamma}
        &&\text{on } \partial_e \gamma \times (0, T)\\
        &\tr \theta_\gamma = \theta_{\partial \gamma} && \text{on } \partial_n
        \gamma\times(0, T)\\
        &\tr \bm{\tau}_\gamma \cdot \bm{n}_{\partial \gamma} =0
        &&\text{on } \partial_i \gamma\times(0, T)\\
        &\theta_\gamma(t=0) = \theta_{\gamma, 0} && \text{in } \gamma \times
        \{0\}
    \end{aligned},
\end{gather}
\end{subequations}
where $\tau_{\partial\gamma}$ in \sib{\joule\per\metre\per\second}
and $\theta_{\partial\gamma}$ in \sib{\kelvin} are the thermal flux and
temperature boundary data, respectively, and the third condition is the internal
tip condition. The data $\theta_{\gamma, 0}$ in \sib{\kelvin} is the initial
temperature distribution in the fracture.
System \eqref{eq:reduced_problem_heat} is the reduced system for
temperature.

\subsection{Reduced solute and precipitate model}

The reduced model that describes the evolution of the solute $u_\gamma$ and its
 flux $\bm{\chi}_\gamma$ can be written in the following way
\begin{subequations}\label{eq:reduced_solute_precipitate}
\begin{align}\label{eq:reduced_solute}
    &\begin{aligned}
        &\bm{\chi}_\gamma - \bm{q}_\gamma u_\gamma + \epsilon_\gamma d_\gamma \nabla
        u_\gamma = \bm{0}\\
        &\partial_t (\epsilon_\gamma u_\gamma) + \nabla_\gamma \cdot
        \bm{\chi}
        + \epsilon_\gamma r_w(u_\gamma, w_\gamma; \theta_\gamma)= 0
    \end{aligned}
    &&\text{in } \gamma \times (0, T).
\end{align}
Now the reaction term $r_w$ has units of measure equal to
\sib{\mole\per\cubic\metre\per\second}. To couple the solute in the fracture
with the one in the surrounding porous media, we consider the following
interface condition
\begin{gather}\label{eq:reduced_solute_cc}
    \begin{aligned}
        &\epsilon_\gamma (\tr \bm{\chi}_\Omega \cdot \bm{n}_\gamma - \tr
        \bm{q}_\Omega \cdot \bm{n}_\gamma \tr u_\Omega) + \delta_\gamma
        (u_\gamma - \tr u_\Omega) = 0&& \text{on } \gamma \times (0,
        T).
    \end{aligned}
\end{gather}
For the precipitate in the fracture $w_\gamma$, being the original model an
ordinary differential equation valid for each point in $\Omega$ the reduced
model  becomes simply
\begin{align}\label{eq:reduced_reaction}
    &\partial_t (\epsilon_\gamma w_\gamma) - \epsilon_\gamma r_w(u_\gamma,
    w_\gamma; \theta_\gamma) = 0 &&
    \text{in } \gamma\times (0, T).
\end{align}
The boundary and initial conditions for the reduced solute and precipitate are
finally given by
\begin{gather}\label{eq:reduced_solute_bc}
    \begin{aligned}
        &\tr \bm{\chi}_\gamma \cdot \bm{n}_{\partial \gamma} = \chi_{\partial
        \gamma} && \text{on } \partial_e \gamma\times(0, T)\\
        &\tr u_\gamma = u_{\partial \gamma} && \text{on } \partial_n
        \gamma\times(0, T)\\
        &\tr \bm{\chi}_\gamma \cdot \bm{n}_{\partial \gamma} =0&&
        \text{on } \partial_i \gamma \times (0, T)\\
        &u_\gamma(t=0) = u_{\gamma, 0} && \text{in } \gamma \times \{0\}\\
        &w_\gamma(t=0) = w_{\gamma, 0} && \text{in } \gamma \times \{0\}
    \end{aligned},
\end{gather}
\end{subequations}
where $\chi_{\partial \gamma}$ in \sib{\mole\per\metre\per\second} and
$u_{\partial \gamma}$ in \sib{\mole\per\cubic\metre} are the given boundary conditions
for the flux and solute, respectively. The third conditions is valid for the
internal tips, while $u_{\gamma, 0}$ in \sib{\mole\per\cubic\metre} and $w_{\gamma,
0}$ in \sib{\mole\per\cubic\metre} are the initial conditions in the fracture for
the solute and precipitate, respectively. System
\eqref{eq:reduced_solute_precipitate} describes the evolution of the solute and precipitate
in the fracture.

\subsection{Permeability and aperture model}\label{eq:reduced_permeability_aperture}

We assume that both components of $k$ follow a cubic law which relates them to the aperture, more precisely
\begin{gather}\label{eq:permeability_aperture}
    k_\gamma(\epsilon_\gamma) = k_{\gamma,0}
    \dfrac{\epsilon_\gamma^2}{\epsilon_{\gamma, 0}^2}
    \quad \text{and} \quad
    \kappa_\gamma(\epsilon_\gamma) =
    \kappa_{\gamma, 0}\dfrac{\epsilon_\gamma^2}{\epsilon_{\gamma, 0}^2},
\end{gather}
where $k_{\gamma,0}$ \sib{\square\metre} and $\kappa_{\gamma, 0}$ \sib{\square\metre}
are reference coefficients along and across the
fracture, respectively, and $\epsilon_{\gamma, 0} > 0$ in \sib{\metre} is the initial aperture.
Finally, we consider a similar law of \eqref{eq:porosity} to describe the
evolution of the fracture aperture $\epsilon_\gamma$. We have
\begin{gather}\label{eq:reduced_aperture}
    \begin{aligned}
        &\partial_t \epsilon_\gamma  + \eta_\gamma \epsilon_\gamma \partial_t
        w_\gamma = 0&&
        \text{in } \gamma \times (0, T)\\
        &\epsilon_\gamma(t=0) = \epsilon_{\gamma, 0} && \text{in } \gamma \times \{0\}
    \end{aligned},
\end{gather}
here $\eta_\gamma$ in \sib{\cubic\metre\per\mole} represents the rate of deposition of the solute at the fracture
walls.

\subsection{Complete reduced model}\label{subsec:reduced_complete_model}

The complete setting for a fractured porous media considers the model presented
in Subsection \ref{subsec:complete_model_bulk} for the surrounding porous media
$\Omega$. To describe the evolution in $\gamma$ of the flow, heat, precipitate and solute as
well as the permeability and fracture aperture variations the model is coupled
with \eqref{eq:reduced_problem_darcy}, \eqref{eq:reduced_problem_heat},
\eqref{eq:reduced_solute_precipitate}, \eqref{eq:permeability_aperture}, and
\eqref{eq:reduced_aperture}. Also in this case the system is fully coupled,
non-linear and possibly non-smooth.

\subsection{Multiple fracture intersections}

The previously described models extend straightforwardly to the case of multiple fractures that do
not intersect to each other. In this case equations \eqref{eq:reduced_problem_darcy},
\eqref{eq:reduced_problem_heat},
\eqref{eq:reduced_solute_precipitate}, \eqref{eq:permeability_aperture}, and
\eqref{eq:reduced_aperture} are written for each fracture which, separately, are
coupled with the porous media. A different case is when two or more fractures
intersect each other and a flow exchange between them can take place. Moreover, due to specific physical
properties in the vicinity of the intersection it is a common approach to
allow for different data at the intersection, see
\cite{Formaggia2012,Fumagalli2012g,Fumagalli2017a,Boon2018,Flemisch2016a} to
name a few.

\begin{figure}[tb]
    \centering
    \resizebox{0.33\textwidth}{!}{\fontsize{0.75cm}{2cm}\selectfont
\begingroup%
  \makeatletter%
  \providecommand\color[2][]{%
    \errmessage{(Inkscape) Color is used for the text in Inkscape, but the package 'color.sty' is not loaded}%
    \renewcommand\color[2][]{}%
  }%
  \providecommand\transparent[1]{%
    \errmessage{(Inkscape) Transparency is used (non-zero) for the text in Inkscape, but the package 'transparent.sty' is not loaded}%
    \renewcommand\transparent[1]{}%
  }%
  \providecommand\rotatebox[2]{#2}%
  \newcommand*\fsize{\dimexpr\f@size pt\relax}%
  \newcommand*\lineheight[1]{\fontsize{\fsize}{#1\fsize}\selectfont}%
  \ifx\svgwidth\undefined%
    \setlength{\unitlength}{311.96932983bp}%
    \ifx\svgscale\undefined%
      \relax%
    \else%
      \setlength{\unitlength}{\unitlength * \real{\svgscale}}%
    \fi%
  \else%
    \setlength{\unitlength}{\svgwidth}%
  \fi%
  \global\let\svgwidth\undefined%
  \global\let\svgscale\undefined%
  \makeatother%
  \begin{picture}(1,0.98417801)%
    \lineheight{1}%
    \setlength\tabcolsep{0pt}%
    \put(0,0){\includegraphics[width=\unitlength,page=1]{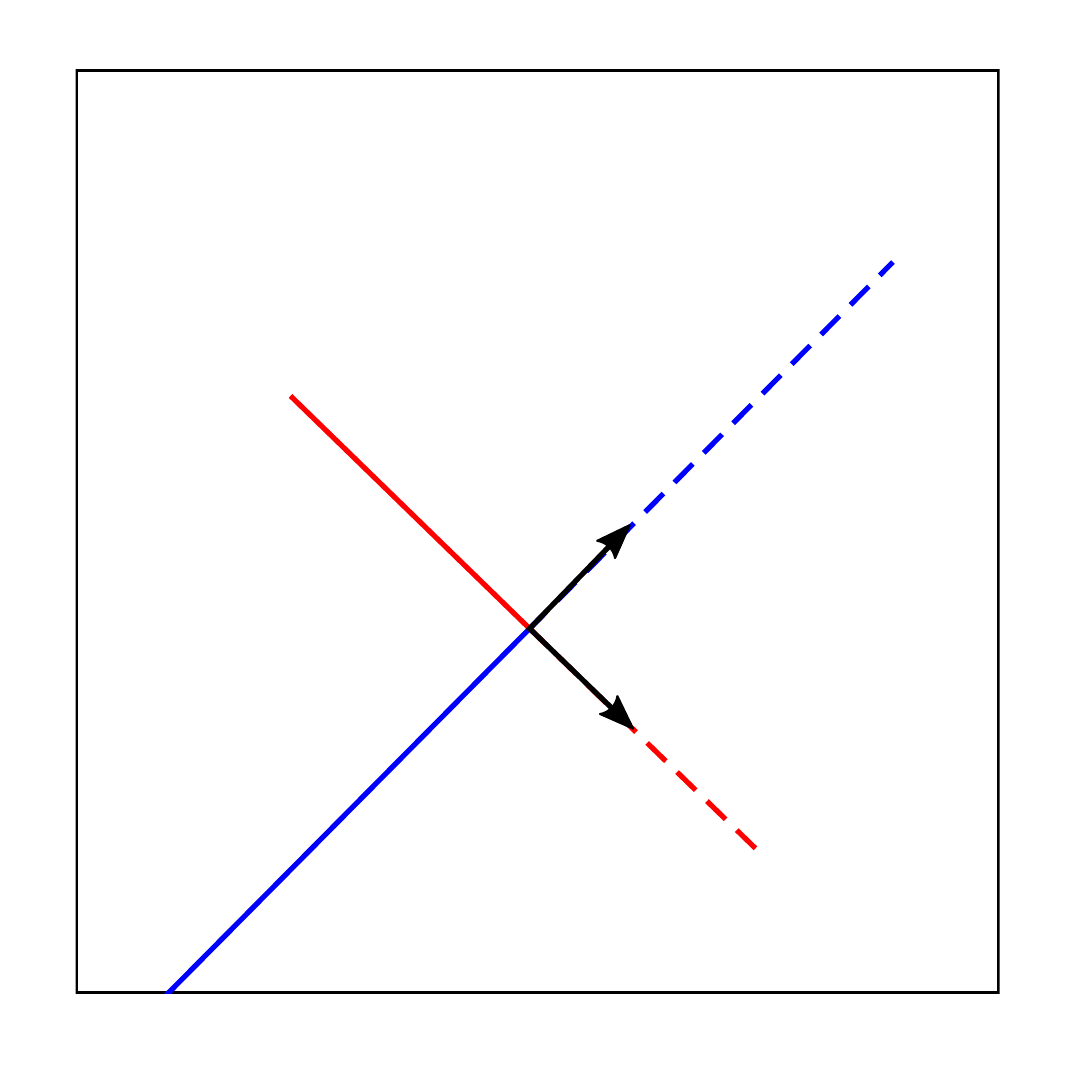}}%
    \put(0.24602719,0.11439133){\color[rgb]{0,0,0}\makebox(0,0)[lt]{\lineheight{1.25}\smash{\begin{tabular}[t]{l}$\gamma_1$\end{tabular}}}}%
    \put(0.26938949,0.64229609){\color[rgb]{0,0,0}\makebox(0,0)[lt]{\lineheight{1.25}\smash{\begin{tabular}[t]{l}$\gamma_2$\end{tabular}}}}%
    \put(0.43304714,0.38483897){\color[rgb]{0,0,0}\makebox(0,0)[lt]{\lineheight{1.25}\smash{\begin{tabular}[t]{l}$\iota$\end{tabular}}}}%
    \put(0.5871651,0.29656798){\color[rgb]{0,0,0}\makebox(0,0)[lt]{\lineheight{1.25}\smash{\begin{tabular}[t]{l}$\bm{n}_{\partial \gamma_2}$\end{tabular}}}}%
    \put(0.58182635,0.49585967){\color[rgb]{0,0,0}\makebox(0,0)[lt]{\lineheight{1.25}\smash{\begin{tabular}[t]{l}$\bm{n}_{\partial \gamma_1}$\end{tabular}}}}%
    \put(0,0){\includegraphics[width=\unitlength,page=2]{domain_intersection.pdf}}%
    \put(0.09978823,0.83831918){\color[rgb]{0,0,0}\makebox(0,0)[lt]{\lineheight{1.25}\smash{\begin{tabular}[t]{l}$\Omega$\end{tabular}}}}%
  \end{picture}%
\endgroup%
}
    \caption{Example of two intersecting fractures with some nomenclature considered.}%
    \label{fig:domain_intersection}
\end{figure}

In this work we consider this approach applied to each equation of the model.
We denote an intersection with $\iota$, which can
be a line if $n=3$ or a point if $n=2$.
Refer to Figure \ref{fig:domain_intersection}.
If the intersection is mono-dimensional then it can be seen as a channel where fluid flow can occur, it is thus approximated with a reduced model similar to \eqref{eq:reduced_problem_darcy},
\eqref{eq:reduced_problem_heat},
\eqref{eq:reduced_solute_precipitate}, \eqref{eq:permeability_aperture}, and
\eqref{eq:reduced_aperture} where the aperture now represents the
cross sectional area of the intersection. If the intersection is
zero-dimensional, resulting also from multiple-intersections of two-dimensional
fractures, the treatment is the following.

As did before, we indicate with $\bm{n}_{\partial \gamma_i}$ the normal
pointing outward from $\gamma_i$ towards the intersection and such that it is also
tangent to the fracture $\gamma_i$.
At the intersection we impose the following conditions for the flow
problem
\begin{gather*}\label{eq:reduced_darcy_intersection}
    \begin{aligned}
        &\mu \epsilon_\iota \tr \bm{q}_{\gamma_i} \cdot \bm{n}_{\partial
        \gamma_i}
        + \kappa_{\gamma_i}(\epsilon_\iota) (p_\iota - \tr p_{\gamma_i}) = 0\\
        &\partial_t \epsilon_\iota + \nabla_\iota \cdot \bm{q} = 0
    \end{aligned}
    \quad \text{on } \iota \times (0, T).
\end{gather*}
where $\nabla_\iota \cdot $ is the zero dimensional conservation operator,
defined as
\begin{gather*}
    \nabla_\iota \cdot \bm{\nu} = -\sum_{i\in N_\iota} \tr \bm{\nu}_{\gamma_i}
    \cdot \bm{n}_{\partial \gamma_i}
\end{gather*}
with $N_\iota$ the set of fractures
meeting in $\iota$ and $\bm{\nu} = (\bm{\nu}_\Omega, \bm{\nu}_{\gamma_1},
\ldots, \bm{\nu}_{\gamma_N})$ the extended compound with $N$ fractures in the
problem. In the previous equation
\eqref{eq:reduced_darcy_intersection} $p_\iota$ is
the pressure at the intersection and $ \kappa_{\gamma_i}$ a fracture specific
permeability associated to $\iota$ which can depend on $\epsilon_\iota$.
Where the latter is the measure of
the intersection in the original equi-dimensional framework.

For the heat equation the coupling conditions are given by
\begin{gather*}
    \begin{aligned}
        &\epsilon_\iota (\tr \bm{\tau}_{\gamma_i} \cdot \bm{n}_{\partial
        \gamma_i} - \rho_w c_w \tr
        \bm{q}_{\gamma_i} \cdot \bm{n}_{\partial \gamma_i} \tr \theta_{\gamma_i})
        + \lambda_w
        (\theta_\iota - \tr \theta_{\gamma_i}) = 0\\
        &\rho_w c_w \partial_t (\epsilon_\iota \theta_\iota) + \nabla_\iota \cdot \bm{\tau} = 0
    \end{aligned}
    \quad \text{on } \iota \times (0, T),
\end{gather*}
where in this case $\theta_\iota$ is the temperature associated to the
intersection $\iota$. For the solute equations the coupling conditions in $\iota$ are given
by
\begin{gather*}
    \begin{aligned}
        &\epsilon_\iota (\tr \bm{\chi}_{\gamma_i} \cdot \bm{n}_{\partial
        \gamma_i} - \tr
        \bm{q}_{\gamma_i} \cdot \bm{n}_{\partial \gamma_i} \tr
        u_{\gamma_i}) + \delta_\iota
        (u_\iota - \tr u_{\gamma_i}) = 0\\
        &\partial_t (\epsilon_\iota u_\iota) + \nabla_\iota \cdot
        \bm{\chi}
        + \epsilon_\iota r_w(u_\iota, w_\iota; \theta_\iota)= 0
    \end{aligned}
    \quad \text{on } \iota \times (0, T),
\end{gather*}
with $u_\iota$ and $w_\iota$ the solute and precipitate in the intersection, while for the precipitate we have
\begin{gather*}
    \partial_t (\epsilon_\iota w_\iota)
    - \epsilon_\iota r_w(u_\iota, w_\iota; \theta_\iota)= 0
    \quad \text{on } \iota \times (0, T).
\end{gather*}
Finally, we can extend also the model for the parameter $\epsilon_\iota$ as well
as for the permeability associated to the intersections. We have the following
\begin{gather*}
    \kappa_{\gamma_i}(\epsilon_\iota) =
    \kappa_{\gamma_i, 0}\dfrac{\epsilon_\iota^2}{\epsilon_{\iota, 0}^2}
    \qquad \text{and} \qquad
    \begin{aligned}
        &\partial_t \epsilon_\iota + \eta_\iota \epsilon_\iota \partial_t
        w_\iota = 0&&
        \text{in } \iota \times (0, T)\\
        &\epsilon_\iota(t=0) = \epsilon_{\iota, 0} && \text{in } \iota \times \{0\}
    \end{aligned},
\end{gather*}
where $\kappa_{\gamma_i, 0}$ and $\epsilon_{\iota, 0}$ are initial values of
$\kappa_{\gamma_i}$ and $\epsilon_\iota $, respectively.

\section{Discretization}\label{sec:discretization}

As discussed in Subsection \ref{subsec:reduced_complete_model},
the problem is fully coupled. We adopt here a (first order in time) splitting scheme
such that legacy codes can be used for its numerical solution. Due to the
relation \eqref{eq:chimica}, it is common to solve the reaction step with an
explicit scheme and split the diffusion and advection parts of
\eqref{eq:solute} and \eqref{eq:reduced_solute}.

The rock domain $\Omega$ is approximated by a grid $\Omega_h$ of non-overlapping elements, whose
regularity is related to the chosen spatial numerical scheme, that completely
cover $\Omega$ itself. We consider here a conforming \cite{Martin2005,Nordbotten2018} approximation of
the fracture grids with respect to the surrounding porous media grid, meaning
that fracture cells are geometrically identical to faces (or edges in 2D) of the porous media
grid. We indicate a generic fracture grid as $\gamma_h$.
The time interval $(0, T)$ is divided, for simplicity, in equally spaced time
steps $\Delta t$ such that $N \Delta t = T$, with $N$ the number of time steps.
Finally, we indicate with $t^n = n \Delta t$ and with a super-script $n$, or $n+1$, the value of
variables or data computed at time $t^n$, or $t^{n+1}$.

\subsection{The temporal splitting scheme}\label{sec:splitting}

In this part we introduce the splitting scheme to solve the global problem
introduced in Subsection \eqref{subsec:reduced_complete_model}. Variables and
data are considered semi-discretized in time but not yet in space.
\begin{scheme}[Temporal splitting scheme]\label{scheme:splitting}
We set the initial condition as
\begin{gather*}
    \phi^0_\Omega = \phi_{\Omega, 0} \quad \epsilon^0_\gamma=\epsilon_{\gamma,
    0} \quad
    \theta^0_\Omega=\theta_{\Omega, 0} \quad \theta^0_\gamma=\theta_{\gamma,
    0}\quad
    u^0_\Omega=u_{\Omega, 0} \\ u^0_\gamma = u_{\gamma, 0} \quad
    w^{-1}_\Omega = w^0_\Omega=w_{\Omega, 0} \quad w^{-1}_\gamma =
    w^0_\gamma=w_{\gamma, 0}.
\end{gather*}
In each time step from $t^n$ to $t^{n+1}$ we perform the following steps:
\begin{enumerate}
    \item Extrapolate the precipitate concentration to get a better estimate of porosity,
        see \cite{Giovanardi2015,Agosti2015}, obtaining
        \begin{gather*}
            w^*_\Omega=2w^n_\Omega-w^{n-1}_\Omega \quad \text{and} \quad
            w^*_\gamma=2w^n_\gamma-w^{n-1}_\gamma.
        \end{gather*}
    \item \label{step:porosity} Compute the corresponding porosity and aperture with an implicit
    discretization of \eqref{eq:porosity} and \eqref{eq:reduced_aperture}, respectively to get
        \begin{gather*}
            \phi_\Omega^*=\frac{\phi_\Omega^n}{1+ \eta_\Omega(w^*_\Omega-w^{n}_\Omega)}
            \quad \text{and} \quad
            \epsilon^*_\gamma=\frac{\epsilon^n_\gamma}{1+ \eta_\gamma(w^*_\gamma-w^{n}_\gamma)}.
        \end{gather*}
        {Note that, with this approximation, these extrapolated values of porosity and aperture cannot become negative when the precipitate increases.}
    \item Update porous media permeability $k_\Omega(\phi_\Omega^*)$ and fracture normal and
    tangential permeabilities $k_\gamma(\epsilon^*_\gamma)$ and
    $\kappa_\gamma(\epsilon^*_\gamma)$ according to \eqref{eq:kozeny} and
    \eqref{eq:permeability_aperture}, respectively.
    \item \label{step:darcy} With $(\phi_\Omega^*, \epsilon^*_\gamma)$ and the computed permeabilities, solve
        the Darcy problem, \eqref{eq:problem_darcy} and \eqref{eq:reduced_problem_darcy},
        to get $(\bm{q}^{n+1}_\Omega, p^{n+1}_\Omega)$ and $(\bm{q}^{n+1}_\gamma,
        p^{n+1}_\gamma)$, respectively.
    \item \label{step:heat} With the advective fields $(\bm{q}^{n+1}_\Omega, \bm{q}^{n+1}_\gamma)$
        computed in the previous point,
        solve the heat equations \eqref{eq:problem_heat} and \eqref{eq:reduced_problem_heat}
        to obtain temperature distribution $\theta^{n+1}_\Omega$ and
        $\theta^{n+1}_\gamma$, respectively.
    \item \label{step:solute_ad} With the advective fields $(\bm{q}^{n+1}_\Omega, \bm{q}^{n+1}_\gamma)$,
        solve the advection-diffusion parts for the solute \eqref{eq:problem_solute_precipitate} and
        \eqref{eq:reduced_solute_precipitate} to get the
        intermediate $u^{n+\frac{1}{2}}_\Omega$ and $u^{n+\frac{1}{2}}_\gamma$, respectively.
    \item Since in the previous advection-diffusion step we have accounted for porosity changes
        using $\phi_\Omega^*$, the new concentration $u_\Omega^{n+\frac{1}{2}}$ and
        $u_\gamma^{n+\frac{1}{2}}$ accounts for the change in pore
        volume, thus, the precipitate needs to be updated accordingly computing
        \begin{gather*}
            w^{n+\frac{1}{2}}_\Omega=w^n_\Omega \frac{\phi_\Omega^n}{\phi_\Omega^*}
            \quad \text{and} \quad
            w^{n+\frac{1}{2}}_\gamma=w^n_\gamma
            \frac{\epsilon^n_\gamma}{\epsilon^*_\gamma}.
        \end{gather*}
    \item \label{step:solute_r} Solve the reaction step starting from
    $(w^{n+\frac{1}{2}}_\Omega,
    w^{n+\frac{1}{2}}_\gamma)$ and $(u^{n+\frac{1}{2}}_\Omega,
    u^{n+\frac{1}{2}}_\gamma)$
    to get $(w^{**}_\Omega, w^{**}_\gamma)$ and $(u^{**}_\Omega, u^{**}_\gamma)$.
    \item At this point we can update the porosity and aperture with the true concentration of
        precipitate at time $n+1$, that is
        \begin{gather*}
            \phi_\Omega^{n+1}=\frac{\phi_\Omega^n}{1+ \eta_\Omega(w^{**}_\Omega-w^{n}_\Omega)}
            \quad \text{and} \quad
            \epsilon^{n+1}_\gamma=\frac{\epsilon^n_\gamma}{1+
            \eta_\gamma(w^{**}_\gamma-w^{n}_\gamma)}.
        \end{gather*}
    \item Finally, we correct the concentrations to account for the difference between the extrapolated and
        ``true'' new porosity and aperture, namely $\phi_\Omega^*$ and $\phi_\Omega^{n+1}$ and
        $\epsilon^*_\gamma$ and $\epsilon^{n+1}_\gamma$, respectively. To ensure mass
        conservation, we thus compute
        \begin{gather*}
            w^{n+1}_\Omega=w^{**}_\Omega \frac{\phi_\Omega^*}{\phi_\Omega^{n+1}}
            \quad
            u^{n+1}_\Omega=u^{**}_\Omega
            \frac{\phi_\Omega^*}{\phi_\Omega^{n+1}}
            \quad
            w^{n+1}_\gamma=w^{**}_\gamma
            \frac{\epsilon^*_\gamma}{\epsilon^{n+1}_\gamma}
            \quad
            u^{n+1}_\gamma=u^{**}_\gamma
            \frac{\epsilon^*_\gamma}{\epsilon^{n+1}_\gamma}.
        \end{gather*}
\end{enumerate}
\end{scheme}

In the following parts we detail the temporal solution of some of the previous
points. In section \ref{subsec:temporal_discretization_darcy}
we discuss Step \ref{step:darcy}, in section
\ref{subsec:temporal_discretization_heat}  the discretized heat model
from Step \ref{step:heat},
in section \ref{subsec:temporal_discretization_solute} the Step
\ref{step:solute_ad}, finally in Subsection \ref{sec:reaction_discr} the
reaction part in Step \ref{step:solute_r}.

\begin{remark}
    Since the order of convergence in time is bounded by the chosen splitting scheme,
    thus of order one, we generally consider low order schemes when high accuracy is not
    needed.
\end{remark}

\begin{remark}
    Steps \ref{step:heat} and \ref{step:solute_ad} can be solve in parallel increasing the
    performances of the code.
    Moreover, the computation of the porosity, aperture, and reaction parts are done
    cell by cell and they are thus embarrassingly parallel.
\end{remark}

\subsubsection{Temporal discretization of Darcy model}\label{subsec:temporal_discretization_darcy}

The Implicit Euler (IE) scheme is used to discretize the temporal derivative of the
porosity, with $\phi_\Omega^*$ from Step \ref{step:porosity} and $\phi_\Omega^n$ from the
previous time step. In \eqref{eq:darcy}, the bulk part of Step \ref{step:darcy}, i.e. the Darcy problem in $\Omega$ is thus solved as
\begin{subequations}\label{eq:problem_darcy_time_discr}
\begin{align}
    &\begin{aligned}
        &\mu \bm{q}^{n+1}_\Omega + k(\phi_\Omega^*) \nabla p^{n+1}_\Omega =\bm{0}\\
        &\phi_\Omega^* - \phi_\Omega^{n} + \Delta t \nabla_\Omega \cdot \bm{q}^{n+1} + \Delta t
        f^{n+1}_\Omega =0
        \end{aligned}
    &\text{in } \Omega\times(t^n, t^{n+1})
\end{align}
The term involving the porosity is now a given data and can be considered as
additional source term. Boundary conditions \eqref{eq:darcy_bc} are adapted accordingly.

For the solution in fracture of Step \ref{step:darcy} we proceed similarly, the
aperture time derivative is dicretized with IE by considering $\epsilon^*_\gamma$ and
$\epsilon^n_\gamma$ from the previous time step. From \eqref{eq:reduced_darcy} the
scheme becomes
\begin{align}
    &\begin{aligned}
        &\mu \bm{q}_\gamma^{n+1} + \epsilon^*_\gamma k_\gamma(\epsilon^*_\gamma) \nabla p_\gamma^{n+1} =\bm{0}\\
        &\epsilon^*_\gamma - \epsilon^{n}_\gamma + \Delta t\nabla_\gamma \cdot \bm{q}^{n+1} + \Delta t
        \epsilon^*_\gamma f_\gamma^{n+1} = 0
    \end{aligned}
    &\text{in } \gamma\times(t^n, t^{n+1})
\end{align}
Also in this case the boundary conditions \eqref{eq:bc_reduced_darcy_model} are discretized likewise.

The discretization of the coupling term \eqref{eq:reduced_darcy_cc}, between the
fracture and the bulk for the flow model, {uses the extrapolated value}  of the aperture
\begin{gather}
    \mu \epsilon^*_\gamma \tr \bm{q}_\Omega^{n+1} \cdot \bm{n}_\gamma +
    \kappa_\gamma(\epsilon^*_\gamma) (p_\gamma^{n+1} - \tr p_\Omega^{n+1}) =
    0 \quad \text{on } \gamma\times(t^n, t^{n+1}).
\end{gather}
\end{subequations}

Because problems in the fractures and surrounding porous media are coupled,
problem \eqref{eq:problem_darcy_time_discr} is solved to obtain the final value
of $(\bm{q}^{n+1}_\Omega, p^{n+1}_\Omega)$ and $(\bm{q}_\gamma^{n+1},
p_\gamma^{n+1})$.

\subsubsection{Temporal discretization of heat model}\label{subsec:temporal_discretization_heat}

By considering the temperature at previous time step $\theta^n_\Omega$ and
$\theta^{n}_\gamma$, also the heat equation is discretized in Step \ref{step:heat}
with the IE scheme. We obtain the following expression for the heat
problem \eqref{eq:heat} in $\Omega$
\begin{subequations}\label{eq:problem_heat_time_discr}
\begin{align}
    &
    \begin{aligned}
        &\bm{\tau}^{n+1}_\Omega - \rho_w c_w \bm{q}^{n+1}_\Omega
        \theta_\Omega^{n+1} + \lambda(\phi_\Omega^*) \nabla \theta^{n+1}_\Omega = \bm{0}\\
        &c(\phi_\Omega^*) \theta^{n+1}_\Omega
        - c(\phi_\Omega^n) \theta^{n}_\Omega
        + \Delta t\nabla_\Omega
        \cdot \bm{\tau}^{n+1} + \Delta t j_\Omega^{n+1} = 0
    \end{aligned}
    &\text{in } \Omega \times (t^n, t^{n+1}),
\end{align}
where $c(\phi_\Omega^*)$ and $\lambda(\phi_\Omega^*)$ are the extrapolated values of the effective
thermal capacity and conductivity at time $n+1$ computed using the value of $\phi_\Omega^*$ in
\eqref{eq:effective_thermal_capacity} and
\eqref{eq:effective_thermal_conductivity} as
\begin{gather*}
    c(\phi_\Omega^*) = \phi_\Omega^* \rho_w c_w + (1-\phi_\Omega^*) \rho_s c_s
    \quad \text{and} \quad
    \lambda(\phi_\Omega^*) = \lambda_w^{\phi_\Omega^*} \lambda_s^{1-\phi_\Omega^*}.
\end{gather*}

For the fracture part of Step \ref{step:heat}, the time derivative discretized
with IE considers the extrapolated value of aperture $\epsilon^*_\gamma$. We obtain the
discretized version of \eqref{eq:reduced_problem_heat_model} given by
\begin{align}
    &\begin{aligned}
        &\bm{\tau}_\gamma^{n+1} - \rho_w c_w \bm{q}_\gamma^{n+1}
        \theta_\gamma^{n+1} + \epsilon^*_\gamma
        \lambda_w \nabla
        \theta_\gamma^{n+1} = \bm{0}\\
        &\rho_w c_w (\epsilon^*_\gamma  \theta_\gamma^{n+1} -
        \epsilon^n_\gamma  \theta_\gamma^{n})
        + \Delta t \nabla_\gamma \cdot
        \bm{\tau}^{n+1} + \Delta t j_\gamma^{n+1} = 0
    \end{aligned}
    &\text{in } \gamma \times (t^n, t^{n+1}).
\end{align}
The discretization of boundary conditions \eqref{eq:reduced_problem_heat_bc}
follows immediately.

The temporal discretization of the interface condition \eqref{eq:reduced_problem_heat_cc}
between the fracture and surrounding porous media in Step \ref{step:heat} is
the following
\begin{align}
    \epsilon^*_\gamma (\tr \bm{\tau}_\Omega^{n+1} \cdot \bm{n}_\gamma - \rho_w c_w \tr
    \bm{q}_\Omega^{n+1} \cdot \bm{n}_\gamma \tr \theta_\Omega^{n+1}) + \lambda_w
    (\theta_\gamma^{n+1} - \tr \theta_\Omega^{n+1}) = 0\quad \text{on } \gamma
    \times (t^n, t^{n+1}).
\end{align}
\end{subequations}

The coupled problem \eqref{eq:problem_heat_time_discr} is solved to obtain the final value
of the primary variables $\theta^{n+1}_\Omega$ and $\theta_\gamma^{n+1}$.

\subsubsection{Temporal discretization of advection-diffusion solute model}\label{subsec:temporal_discretization_solute}

The temporal discretization of the advection and diffusion parts of the solute in
Step \ref{step:solute_ad} is obtained by setting
the corresponding reaction term to zero. We consider the IE scheme for the
temporal discretization of \eqref{eq:solute}, which, in the porous matrix, reads
\begin{subequations}\label{eq:problem_solute_time_discr}
\begin{align}
    &\begin{aligned}
        &\bm{\chi}^{n+1}_\Omega - \bm{q}^{n+1}_\Omega u^{n+1} + \phi_\Omega^* d \nabla
        u^{n+1}_\Omega =
        \bm{0}\\
        &\phi_\Omega^* u^{n+1}_\Omega - \phi_\Omega^n u^n_\Omega + \Delta t \nabla_\Omega \cdot \bm{\chi}^{n+1}
        = 0
    \end{aligned}
    &\text{in } \Omega\times(t^n, t^{n+1}),
\end{align}
with boundary conditions of the solute in $\Omega$
\eqref{eq:solute_precipitate_bc} easily implemented.

The fracture part of Step \ref{step:solute_ad} consists in solving equation
\eqref{eq:reduced_solute} with null reaction term. Also in this case we take
extrapolated aperture $\epsilon^*_\gamma$ in the time discretization. The equations
become
\begin{align}
    &\begin{aligned}
        &\bm{\chi}_\gamma^{n+1} - \bm{q}_\gamma^{n+1} u_\gamma^{n+1} +
        \epsilon^*_\gamma d_\gamma \nabla u^{n+1}_\gamma =
        \bm{0}\\
        &\epsilon^*_\gamma u_\gamma^{n+1} - \epsilon^n_\gamma u_\gamma^n + \Delta t
        \nabla_\gamma \cdot \bm{\chi}^{n+1}
        = 0
    \end{aligned}
    &\text{in } \gamma\times(t^n, t^{n+1}).
\end{align}

Finally, the coupling conditions \eqref{eq:reduced_solute_cc} between the
fracture and surrounding porous media are discretized as
\begin{gather}
    \begin{aligned}
        &\epsilon^*_\gamma (\tr \bm{\chi}_\Omega^{n+1} \cdot \bm{n}_\gamma - \tr
        \bm{q}_\Omega^{n+1} \cdot \bm{n}_\gamma \tr u_\Omega^{n+1}) + \delta
        (u_\gamma^{n+1} - \tr u_\Omega^{n+1}) = 0&& \text{on } \gamma \times
        (t^n, t^{n+1}).
    \end{aligned}
\end{gather}
\end{subequations}

Coupled problem \eqref{eq:problem_solute_time_discr} is solved to obtain the final value
of the primary variables $u^{n+1}_\Omega$ and $u_\gamma^{n+1}$.

\subsection{Integration in time of the discontinuous reaction
problem}\label{sec:reaction_discr}

\begin{figure}[tb]
    \centerline{\resizebox{0.5\textwidth}{!}{\fontsize{0.8cm}{2cm}\selectfont
\begingroup%
  \makeatletter%
  \providecommand\color[2][]{%
    \errmessage{(Inkscape) Color is used for the text in Inkscape, but the package 'color.sty' is not loaded}%
    \renewcommand\color[2][]{}%
  }%
  \providecommand\transparent[1]{%
    \errmessage{(Inkscape) Transparency is used (non-zero) for the text in Inkscape, but the package 'transparent.sty' is not loaded}%
    \renewcommand\transparent[1]{}%
  }%
  \providecommand\rotatebox[2]{#2}%
  \newcommand*\fsize{\dimexpr\f@size pt\relax}%
  \newcommand*\lineheight[1]{\fontsize{\fsize}{#1\fsize}\selectfont}%
  \ifx\svgwidth\undefined%
    \setlength{\unitlength}{504.67212106bp}%
    \ifx\svgscale\undefined%
      \relax%
    \else%
      \setlength{\unitlength}{\unitlength * \real{\svgscale}}%
    \fi%
  \else%
    \setlength{\unitlength}{\svgwidth}%
  \fi%
  \global\let\svgwidth\undefined%
  \global\let\svgscale\undefined%
  \makeatother%
  \begin{picture}(1,0.77131937)%
    \lineheight{1}%
    \setlength\tabcolsep{0pt}%
    \put(0,0){\includegraphics[width=\unitlength,page=1]{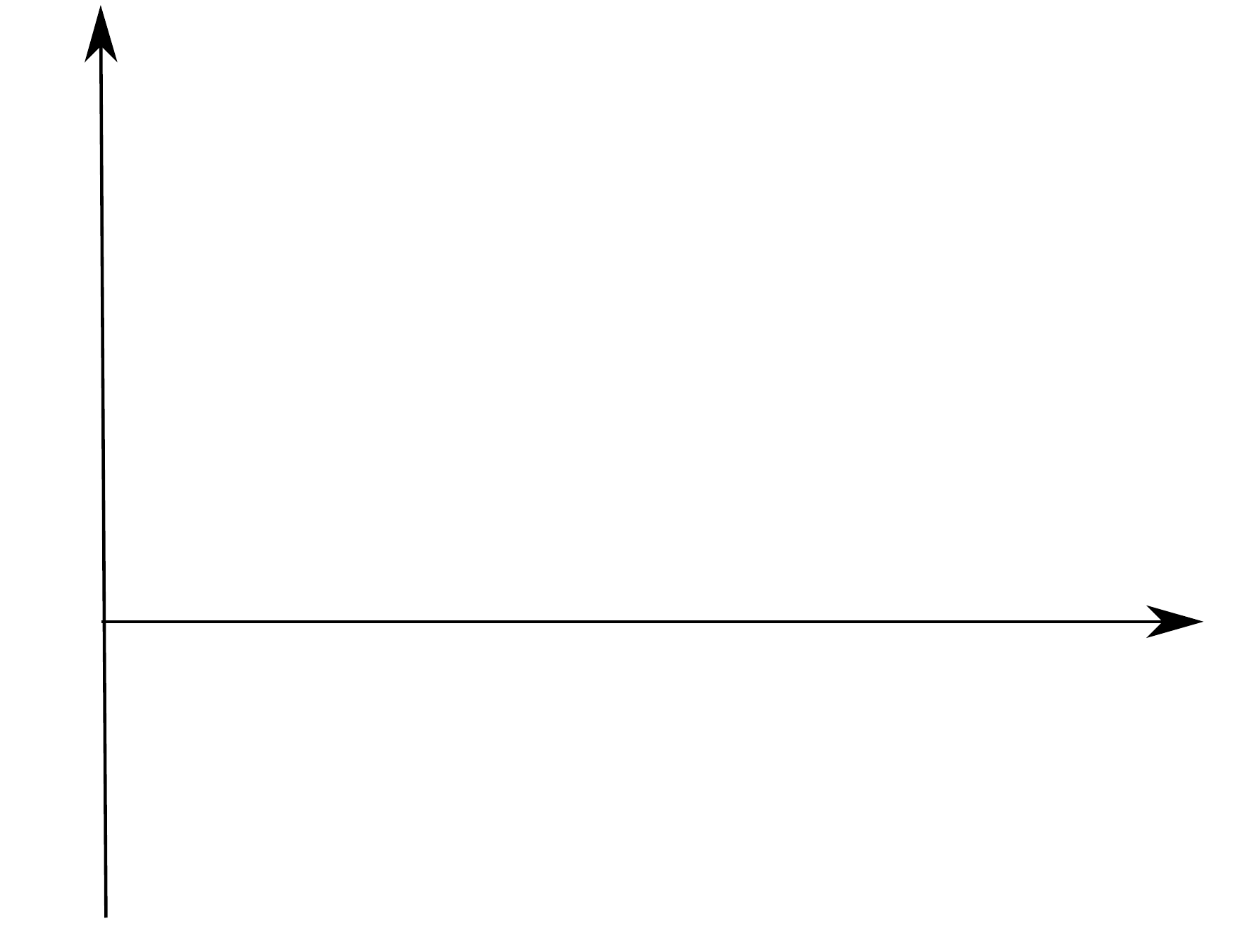}}%
    \put(0.96772346,0.22079401){\color[rgb]{0,0,0}\makebox(0,0)[lt]{\lineheight{1.25}\smash{\begin{tabular}[t]{l}$u$\end{tabular}}}}%
    \put(0.01130342,0.72501102){\color[rgb]{0,0,0}\makebox(0,0)[lt]{\lineheight{1.25}\smash{\begin{tabular}[t]{l}$w$\end{tabular}}}}%
    \put(0,0){\includegraphics[width=\unitlength,page=2]{rw.pdf}}%
    \put(0.53075474,0.21973246){\color[rgb]{0,0,0}\makebox(0,0)[lt]{\lineheight{1.25}\smash{\begin{tabular}[t]{l}$u_{e}$\end{tabular}}}}%
    \put(0,0){\includegraphics[width=\unitlength,page=3]{rw.pdf}}%
  \end{picture}%
\endgroup%
}}
    \caption{Qualitative representation of the vector forcing term $\bm{r}_w$ in
    the phase space $(u,w)$.}\label{fig:r}
\end{figure}

As explained in the previous section we rely on a first order splitting to the
solute equation in the bulk and in the fractures to  separate the advection and
diffusion part from the reaction term. This is motivated by the discontinuous
nature of the reaction term, which benefits from an ad hoc numerical treatment.
Starting from the intermediate solution
$(u_\Omega^{n+\frac{1}{2}},w_\Omega^{n+\frac{1}{2}})$,
$(u_\gamma^{n+\frac{1}{2}},w_\gamma^{n+\frac{1}{2}})$ we integrate, for each
degree of freedom in the porous medium and in the fractures, a system of two
ordinary differential equations. Note that, indeed, after discretization in
space, $u$ and $w$ will be approximated as piecewise constant, thus, with an
abuse of notation, we write for each cell the following system of ODEs
\begin{equation}\label{eq:diff_inclusion}
    d_t [ u, w ]^\top \in\bm{r}_w(u, v; \theta)
    \quad\text{with}\quad \bm{r}_w(u, w; \theta) =[ -r_w(u,w; \theta),\,
    r_w(u,w;\theta)]^\top,
\end{equation}
where $r_w$ is defined as in
\eqref{eq:chimica}. Note that \eqref{eq:diff_inclusion} is written as a
differential inclusion because, strictly speaking, the forcing term is not
defined at $w=0$. It is also important to highlight that the discontinuity
depends on the solution itself and not simply on time. The vector forcing term
$\bm{r}_w$ is represented qualitatively in Figure \ref{fig:r} where we can
observe that it is discontinuous across the line $w=0$ for $u<u_{e}$, in other
words, when the solute concentration is such that we should have precipitate
dissolution.  We can define
\begin{equation*}
\bm{r}_w(u, w; \theta)=\begin{dcases*}
              \bm{r}_w^+(u; \theta) &if $w>0$\\
              \bm{r}_w^-(u; \theta) &if $w<0$
             \end{dcases*},
\end{equation*}
where $\bm{r}_w^+$ and $\bm{r}_w^-$ are defined as
\begin{gather*}
    \bm{r}_w^+(u; \theta)=\lambda^-(\theta)[-(r(u)-1),(r(u)-1)]^\top, \\
    \bm{r}_w^-(u;
    \theta)=\lambda^-(\theta)[-\max(r(u)-1,0),\max(r(u)-1,0)]^\top.
\end{gather*}
Equation \eqref{eq:diff_inclusion} is integrated numerically with an explicit
scheme (Explicit Euler or two-stages Runge-Kutta) combined with an event
location strategy. For the sake of simplicity we will describe the procedure in
the case of the EE scheme. In particular at each step we need to:
\begin{enumerate}
 \item detect if and when the solution reaches the discontinuity threshold, i.e.
 the line $w=0$. Note that this instant, denoted as $\overline{t}$, usually does
 not coincide with $t^n$ or $t^{n+1}$;
 \item at $\overline{t}$ stop and restart the numerical integration with a suitable
 forcing term: $\bm{r}^+_w$, $\bm{r}^-_w$ or a convex combination of the two.
\end{enumerate}

\subsubsection{Detection of the event}

\begin{figure}[bt]
\centerline{\resizebox{0.5\textwidth}{!}{\fontsize{0.8cm}{2cm}\selectfont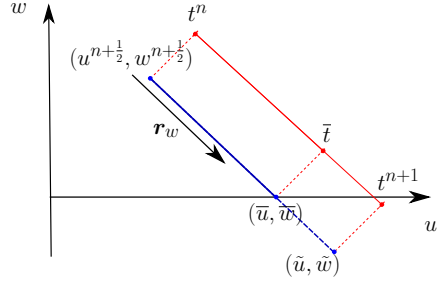}}
\caption{splitting of an integration step across the discontinuity.}\label{fig:disc_int}
\end{figure}

Following \cite{DL} we perform a tentative integration step starting from $t^n$
and the initial conditions $(u^{n+\frac{1}{2}}, w^{n+\frac{1}{2}})$ to obtain
\begin{gather*}
    \begin{aligned}
    \tilde{u}=u^{n+\frac{1}{2}}-\Delta t
    r_w(u^{n+\frac{1}{2}},w^{n+\frac{1}{2}})\\
    \tilde{w}=w^{n+\frac{1}{2}}+ \Delta t r_w(u^{n+\frac{1}{2}},w^{n+\frac{1}{2}})
    \end{aligned}.
\end{gather*}
We then check whether $\tilde{w}<0$: if this is the case it means that in the
$n$-th integration step we are crossing the discontinuity threshold. We proceed
searching for the exact time of the event by defining
$w(\xi)=w^{n+\frac{1}{2}}+\xi \Delta t r_w(u^{n+\frac{1}{2}},w^{n+\frac{1}{2}})$
and search for $\overline{\xi}$ such that $w(\overline{\xi})=0$ by means of a
suitable iterative method.  Once we have detected the time of transition
$\overline{t}=t^n+\overline{\xi}\Delta t$ we perform the following steps:
\begin{enumerate}
 \item integrate the equations from $t^n$ to $\overline{t}$ with the ``old''
 value of the right hand side obtaining $\overline{u}$ and $\overline{w}$ as
 shown in Figure \ref{fig:disc_int};
 \item integrate from
 $t^n+\overline{\xi}\Delta t$ to $t^{n+1}$ with a new value of the right hand
 side.
\end{enumerate}

\subsubsection{Choice of the right hand side}

To determine the behavior of the solution at the discontinuity we let $\bm{n}$
be the normal to the surface of discontinuity in the phase space, in our case
$\bm{n}=[0,\,1]^T$, and observe the sign of $\bm{r}_w^\pm\cdot\bm{n}$.
Suppose for the sake of simplicity that we start from the ``$+$'' region, i.e.
$w^{n+\frac{1}{2}}>0$ and $u<u_{e}$: thus  $\bm{r}^+_w \cdot\bm{n}<0$ meaning
that precipitate is decreasing and we are approaching the discontinuity. On the
opposite side,  $\bm{r}^-_w \cdot\bm{n}=0$. With these conditions the solution,
after $\overline{t}$, should \emph{slide} on the discontinuity, i.e. we finish
the integration of the time step with forcing term $\bm{r}_w^-$, which is
null until the solute reaches its equilibrium concentration. In other
words, stops if the precipitate is not present.

Note that, if $u>u_{e}$, $\bm{r}^\pm_w\cdot\bm{n}>0$, and there is no
discontinuity across $w=0$ because in this case the net reaction rate yields an
increase of precipitate.

\subsection{Spatial discretization}

The discretization in space of the problems in Subsection
\ref{subsec:reduced_complete_model} is not the main focus of this work and
relies on a standard approach.
Since equations are in mixed-dimensions, the numerical schemes considered are
applied in different dimensions.

The main request for the discretization of the flow problem
\eqref{eq:problem_darcy_time_discr} is to obtain a reliable approximation of Darcy velocity that is
locally mass conservative. To solve problem \eqref{eq:problem_darcy_time_discr} we consider two different schemes depending on the geometrical properties of the grids, the Raviart-Thomas of lowest order, see
\cite{Raviart1977,Roberts1991,Martin2005,Boffi2013}, and the mixed virtual
element method of lowest degree, see for example
\cite{Brezzi2014,Benedetto2014,BeiraoVeiga2016,Benedetto2016,Fumagalli2016a,Fumagalli2017a}.
Both schemes handle in an accurate way strong variations of the
permeability tensor which is a typical situation in the underground. For the
numerical solution of problems \eqref{eq:problem_heat_time_discr} and
\eqref{eq:problem_solute_time_discr} we consider an upstream approximation for the
advective part and a two-point flux approximation for the diffusion component,
see \cite{Eymard2000,Jaffre2002,Faille2005,Droniou2013}.

For the coupling between the fracture and the porous media, for simplicity we
assume a conforming strategy meaning that the element of the fracture grids are
composed by faces or edges of the porous media elements neighbouring
the fracture. Other choices are possible that relax some of the geometrical constraints
posed by this approach, see for example \cite{Flemisch2016,Scialo2017}.

\section{Numerical examples}\label{sec:numerical_examples}

In this section we present three numerical examples to show the performances of
the previously introduced mathematical models and splitting scheme. In
particular, in the examples reported in Subsection \ref{subsec:one_d_test} we
validate the good properties of the splitting scheme of Scheme
\ref{scheme:splitting} presented in Section \ref{sec:splitting}. The next examples, presented in Subsection
\ref{subsec:example_single_fracture} and \ref{subsec:example_multiple_fracture}, consider the full problem
with single and multiple fractures, respectively. In these examples, we present
the relevant physical effects the proposed model is able to reproduce by increasing
the geometrical complexity.
These later examples are developed with the library PorePy, a simulation
tool for fractured and deformable porous media written in Python, see \cite{Keilegavlen2019a}.
The scripts associated are freely accessible.

\subsection{1D tests}\label{subsec:one_d_test}

\begin{figure}[tb]
    \centering
    \includegraphics[width=0.33\textwidth]{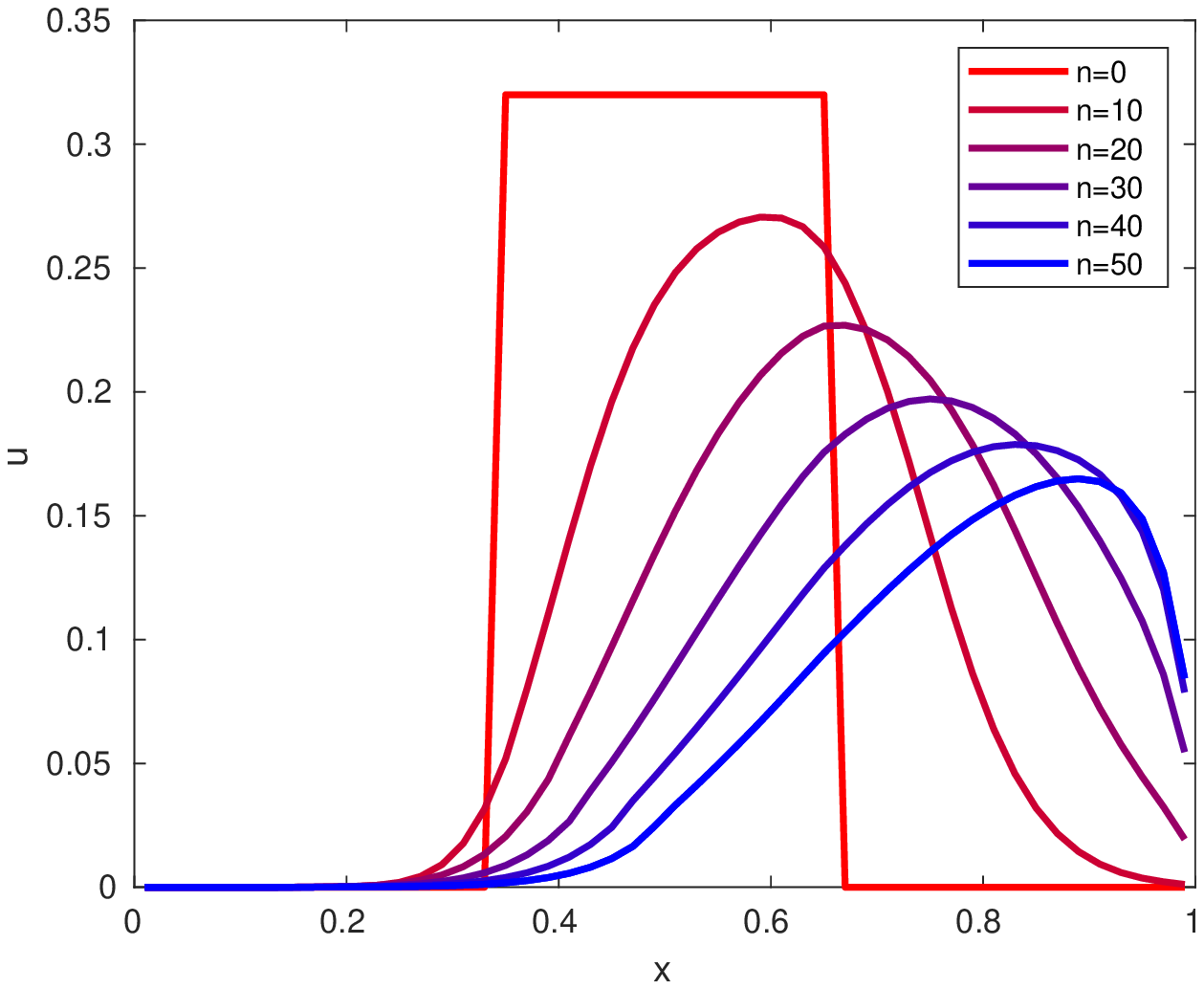}\includegraphics[width=0.33\textwidth]{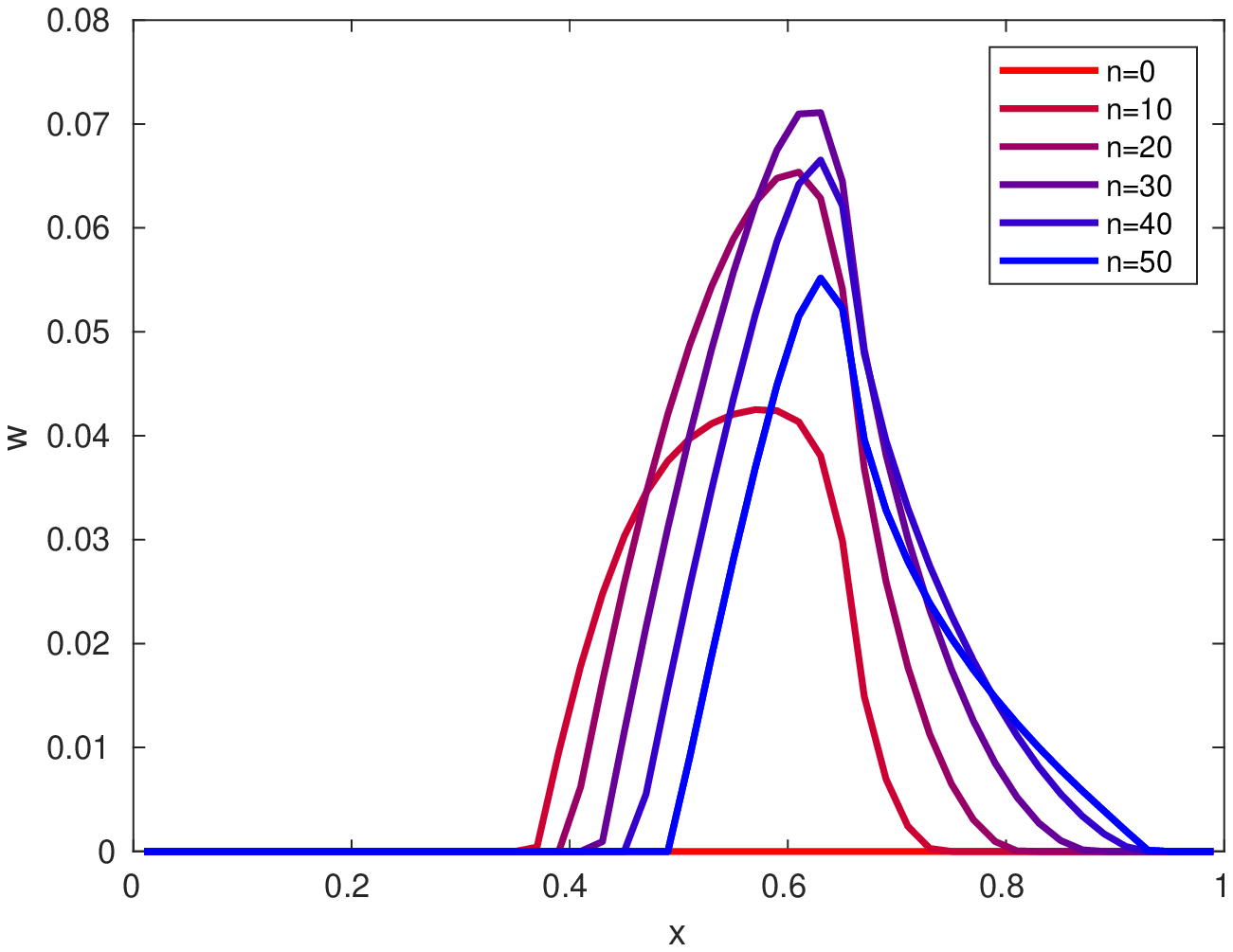}\\
    \includegraphics[width=0.33\textwidth]{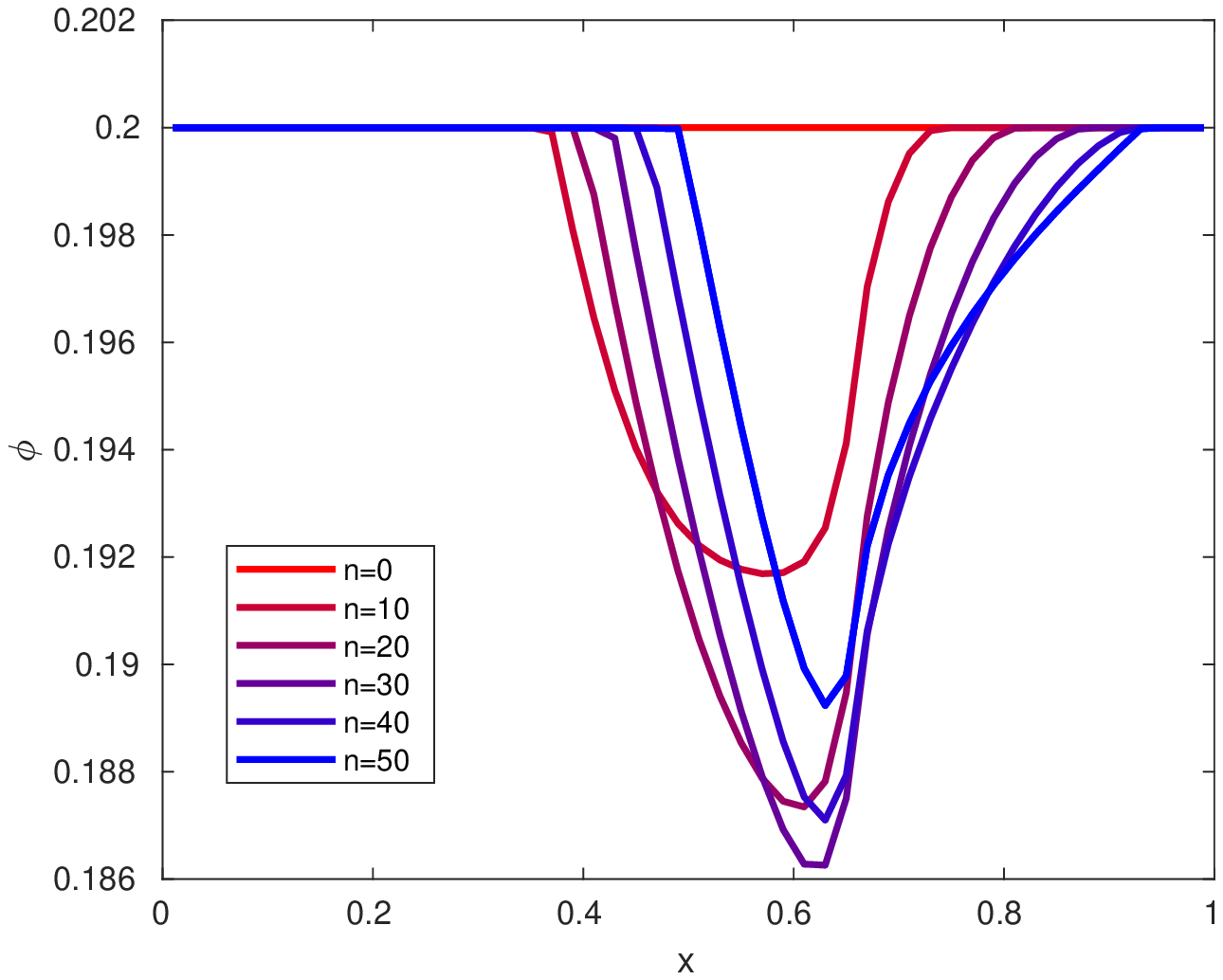}\includegraphics[width=0.33\textwidth]{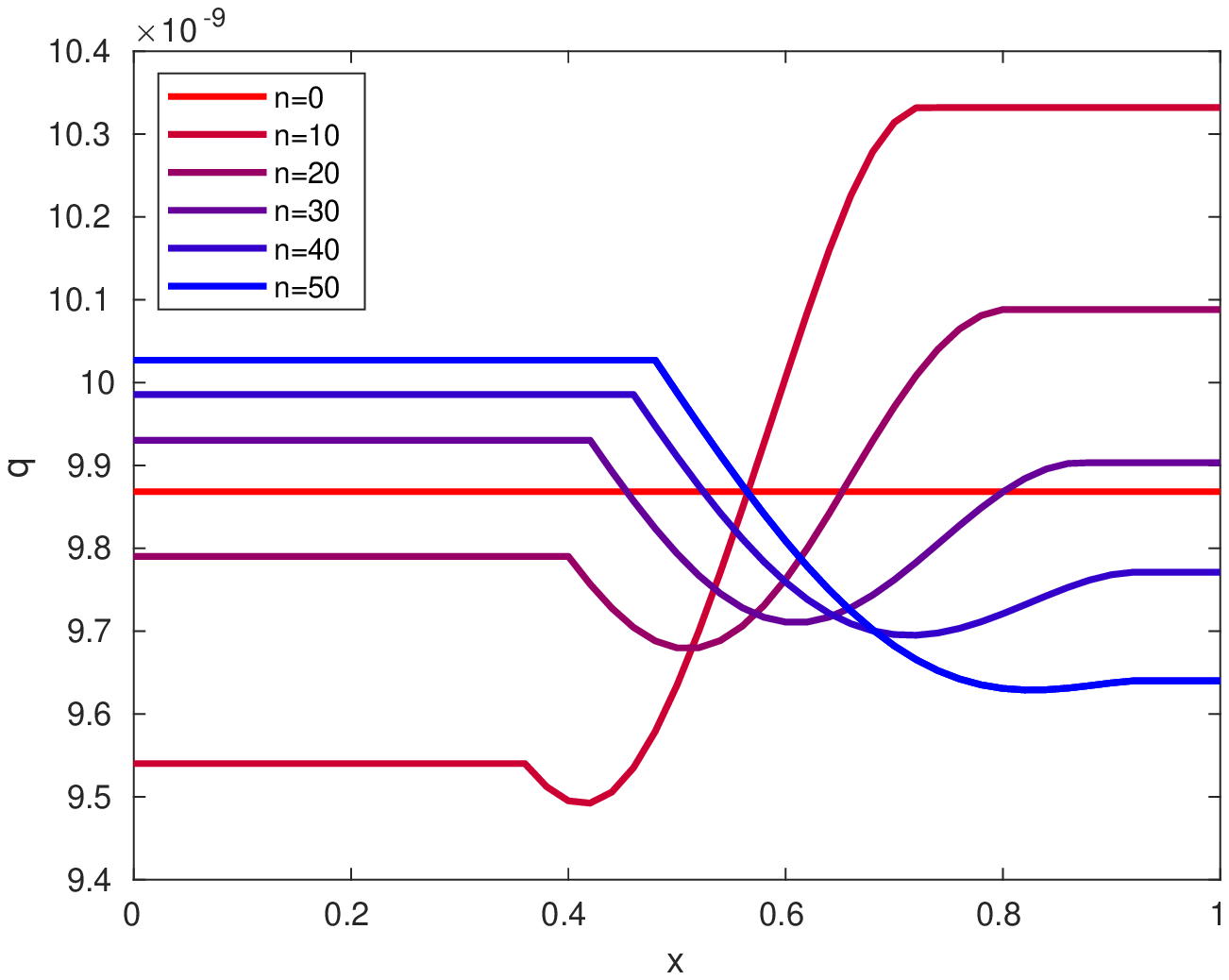}\includegraphics[width=0.33\textwidth]{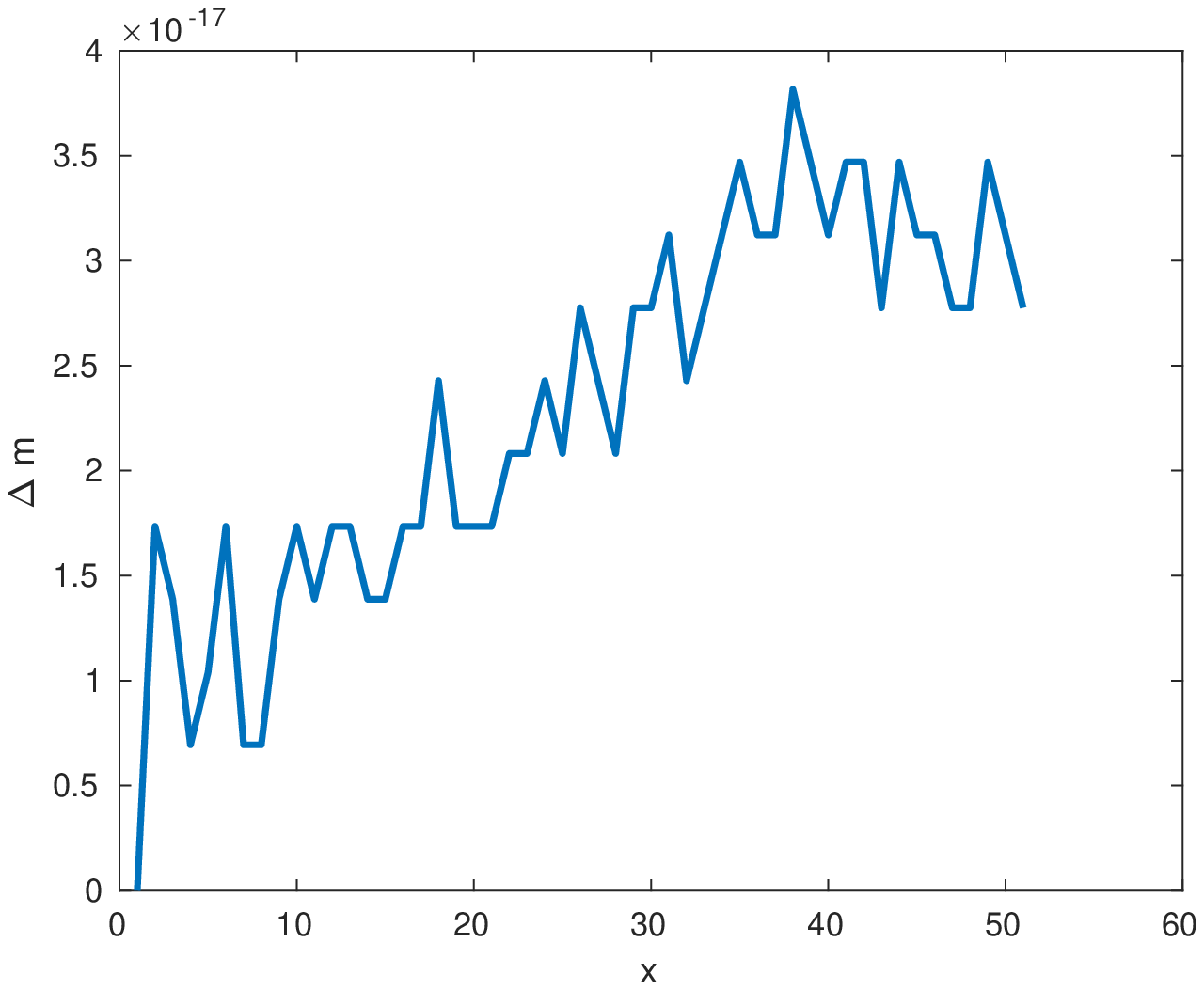}
    \caption{Precipitation and dissolution with a given advection field in 1D at
    different time steps. Top row: solute e precipitate concentrations. Bottom:
    porosity distribution, Darcy velocity and mass balance over time.}%
    \label{fig:test1D}
\end{figure}

Let us consider first a simple one dimensional test case on a domain
$\Omega=(0,1)$ without any fracture. The goal is to test the reliability of the
algorithm presented in Section \ref{sec:splitting}. Solute concentration is set,
at the initial time, to $u=2u_{e}$ in the central part of the domain, while
$w=0$ everywhere at $t=0$.  The advection/diffusion ratio is about 10, and
$CFL=\frac{q\Delta t}{\Delta x}\simeq 8\cdot 10^{-2}$. The reaction rate is such that
$\mathbb{D}a=0.2$, and the coefficient in
\eqref{eq:permeability_aperture} $\eta_\Omega=1$. A pressure drop imposed at the
boundaries results in a Darcy velocity $q\simeq 10^{-8}$, however, as a result of
porosity changes this value will be perturbed.  In Figure \ref{fig:test1D} we
can observe the evolution of $u$, $w$, $\phi$, $q$ and the corresponding mass
balance. We observe that $u$ decreases due to precipitation, and is advected
towards the right boundary. The precipitate grows, initially, but is later
``washed away'' by water with a solute concentration that is lower than the
equilibrium one. The evolution of porosity reflects that of $w$. The quantity
$\Delta m$ is computed as the mass loss/gain that is the difference between mass
of $u$ and $w$ inside the domain, the outflow of $u$ and the initial mass: we
observe that mass is conserved within machine precision.

We use a similar, but simplified setup to evaluate the operator splitting error
by comparing the solution obtained with a monolithic method with the one
obtained with the first order splitting for advection-diffusion and reaction. In
this case we neglect the changes in porosity, consider a constant, given Darcy
velocity and take a linear reaction rate $r(u)=u$.  The value of Darcy velocity
changes at fixed reaction rate to obtain different Damkh\"oler numbers. We can
observe in table \ref{tab:err_splitting} that \textit{i)}
the splitting error decreases linearly with $\Delta t$ for $CFL<1$, \textit{ii)}
however, for $CFL$ close to 1 the splitting error does not decrease with the
correct rate, and \textit{iii)}
the largest errors are obtained for $\mathbb{D}a=1$, i.e. when advection and reaction occur at the same speed.

\begin{table}[tb]
    \centering
    \begin{tabular}{|c|c|c|c|c|} \hline
        $N_t$ $\diagdown$ $Da$ ($CFL$)  &  $0.1 \, \, (8.64 \cdot 10^{-1})$ & $1
        \, \,
        (8.64\cdot 10^{-2})$ & $10 \, \, (8.64\cdot 10^{-3})$ & $100 \, \, (8.64\cdot
        10^{-4})$\\ \hline
        $50$  &   $1.4 \cdot 10^{-3}$ & $3.2 \cdot 10^{-3}$  & $1.4 \cdot
        10^{-3}$&
        $4.9759\cdot 10^{-4}$\\ \hline
        $100$   &  $1.1 \cdot 10^{-3}$& $1.8 \cdot 10^{-3}$ & $7.4362\cdot
        10^{-4}$  &  $2.7832\cdot 10^{-4}$\\ \hline
        $200$ &    $7.6036 \cdot 10^{-4}$ & $1.0 \cdot 10^{-3}$ & $3.7707\cdot
        10^{-4}$  &  $1.4577\cdot 10^{-4}$ \\ \hline
    \end{tabular}
    \caption{$L^\infty$ norm of the difference between the solutions obtained
    with and without the splitting at the final time. We consider different
    advection velocities and increase the number of time steps
    $N_t$. }\label{tab:err_splitting}
\end{table}

Finally, we want to show the impact of Damkh\"oler number on the precipitate
distribution in the presence of a point source to mimic the effect of a fracture
in 2D/3D and predict whether the effect of fracture flow will result in a local
or more diffused change in the porous matrix properties. To this aim  we start
from clean water ($u=0$) and inject a prescribed concentration $u=2u_{e}$ at
$x=0.5$. The velocity field is given, $q=\overline{q} \sign (x-0.5)$.

\begin{figure}
    \centering
    \begin{tabular}{cccc}
        $Da=0.006$ & $Da=0.066$ & $Da=0.662$ & $Da=6.62$\\
        \includegraphics[width=0.2\textwidth]{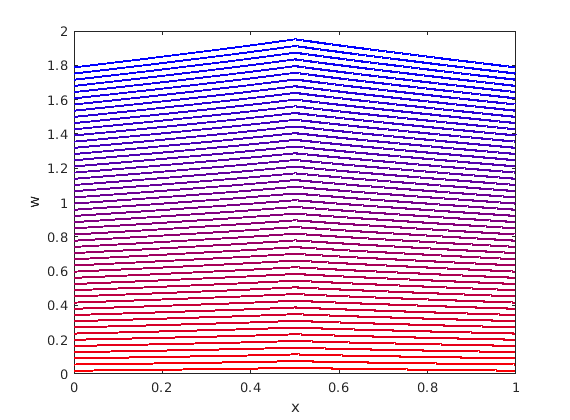} & \includegraphics[width=0.2\textwidth]{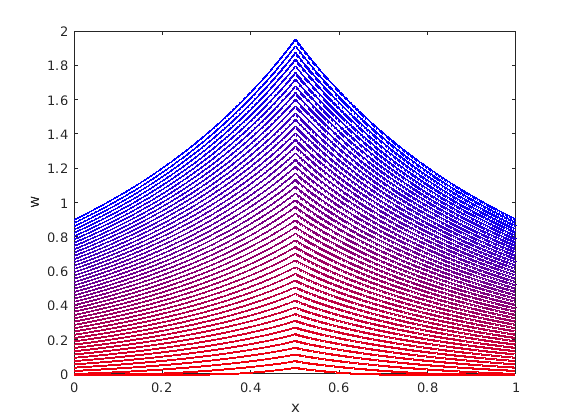} & \includegraphics[width=0.2\textwidth]{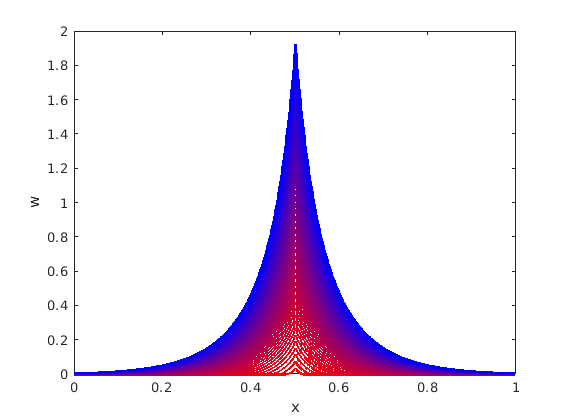} & \includegraphics[width=0.2\textwidth]{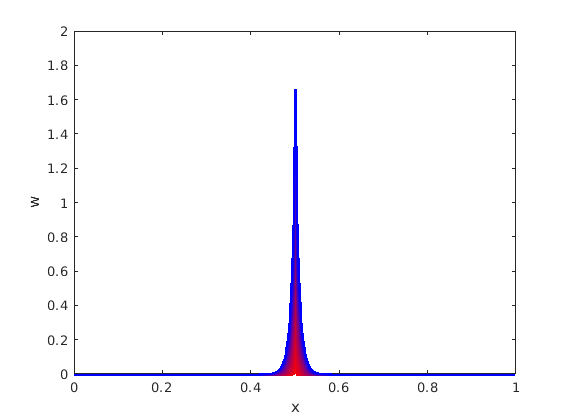}
    \end{tabular}
    \caption{Precipitate concentration for different Darcy velocities,
    corresponding to a point-wise injection of water with $u=2u_{e}$ at $x=0.5$.}\label{fig:Da1}
\end{figure}

If we consider different characteristic Darcy velocities or, in other words,
different Damkh\"ohler numbers, we obtain the results represented in Figure
\ref{fig:Da1}: high advection (small $\mathbb{D}a$) correspond to a flat solute
concentration profile and thus a uniform precipitation, whereas if
advection is very slow compared to reaction, precipitation is focused around the
injection point.

\begin{figure}
    \centering
    \begin{tabular}{cccc}
        $Da=0.006$ & $Da=0.066$ & $Da=0.662$ & $Da=6.62$\\
        \includegraphics[width=0.2\textwidth]{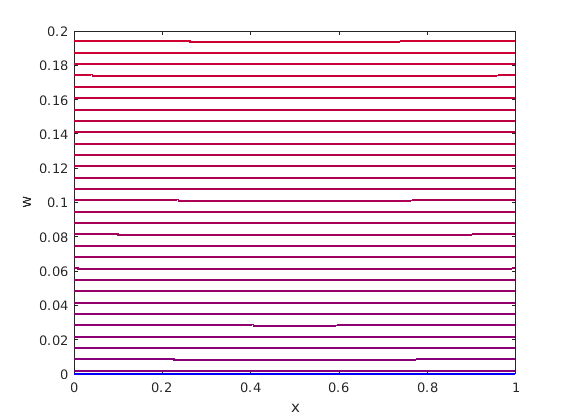} & \includegraphics[width=0.2\textwidth]{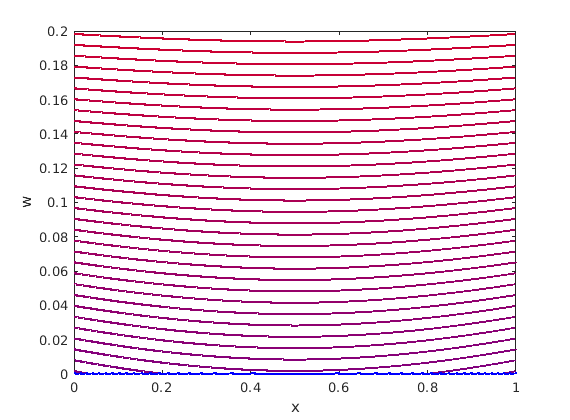} & \includegraphics[width=0.2\textwidth]{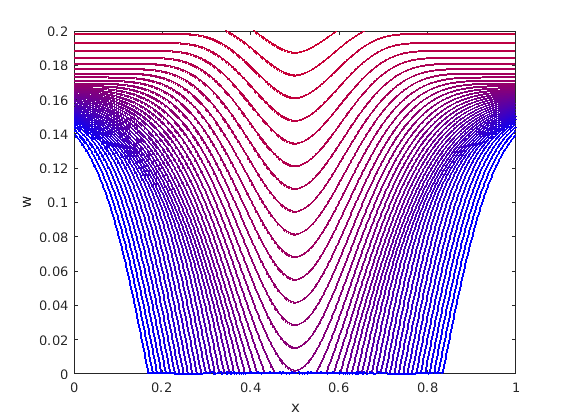} & \includegraphics[width=0.2\textwidth]{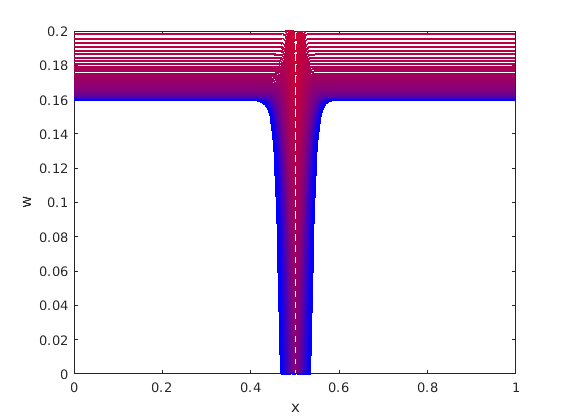}
    \end{tabular}
    \caption{Precipitate concentration for different Darcy velocities,
    corresponding to a point-wise injection of clean water in a porous matrix with $w=2u_{e}$ at $x=0.5$.}\label{fig:Da3}
\end{figure}

Conversely, if we inject clean water into the ``well'' and we start from the
initial constant concentrations $u=0$, $w=2 u_{e}$ we obtain the results in
Figure \ref{fig:Da3}, where, for large advection velocities the precipitate is
eroded uniformly in the whole domain, whereas for large $\mathbb{D}a$ the effect
is concentrated around the injection point.

\subsection{Single fracture}\label{subsec:example_single_fracture}

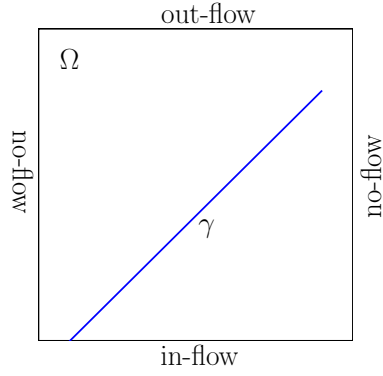
\begin{figure}
    \centering
    \resizebox{0.33\textwidth}{!}{\fontsize{1cm}{2cm}\selectfont
\begingroup%
  \makeatletter%
  \providecommand\color[2][]{%
    \errmessage{(Inkscape) Color is used for the text in Inkscape, but the package 'color.sty' is not loaded}%
    \renewcommand\color[2][]{}%
  }%
  \providecommand\transparent[1]{%
    \errmessage{(Inkscape) Transparency is used (non-zero) for the text in Inkscape, but the package 'transparent.sty' is not loaded}%
    \renewcommand\transparent[1]{}%
  }%
  \providecommand\rotatebox[2]{#2}%
  \newcommand*\fsize{\dimexpr\f@size pt\relax}%
  \newcommand*\lineheight[1]{\fontsize{\fsize}{#1\fsize}\selectfont}%
  \ifx\svgwidth\undefined%
    \setlength{\unitlength}{311.96932983bp}%
    \ifx\svgscale\undefined%
      \relax%
    \else%
      \setlength{\unitlength}{\unitlength * \real{\svgscale}}%
    \fi%
  \else%
    \setlength{\unitlength}{\svgwidth}%
  \fi%
  \global\let\svgwidth\undefined%
  \global\let\svgscale\undefined%
  \makeatother%
  \begin{picture}(1,0.98417801)%
    \lineheight{1}%
    \setlength\tabcolsep{0pt}%
    \put(0,0){\includegraphics[width=\unitlength,page=1]{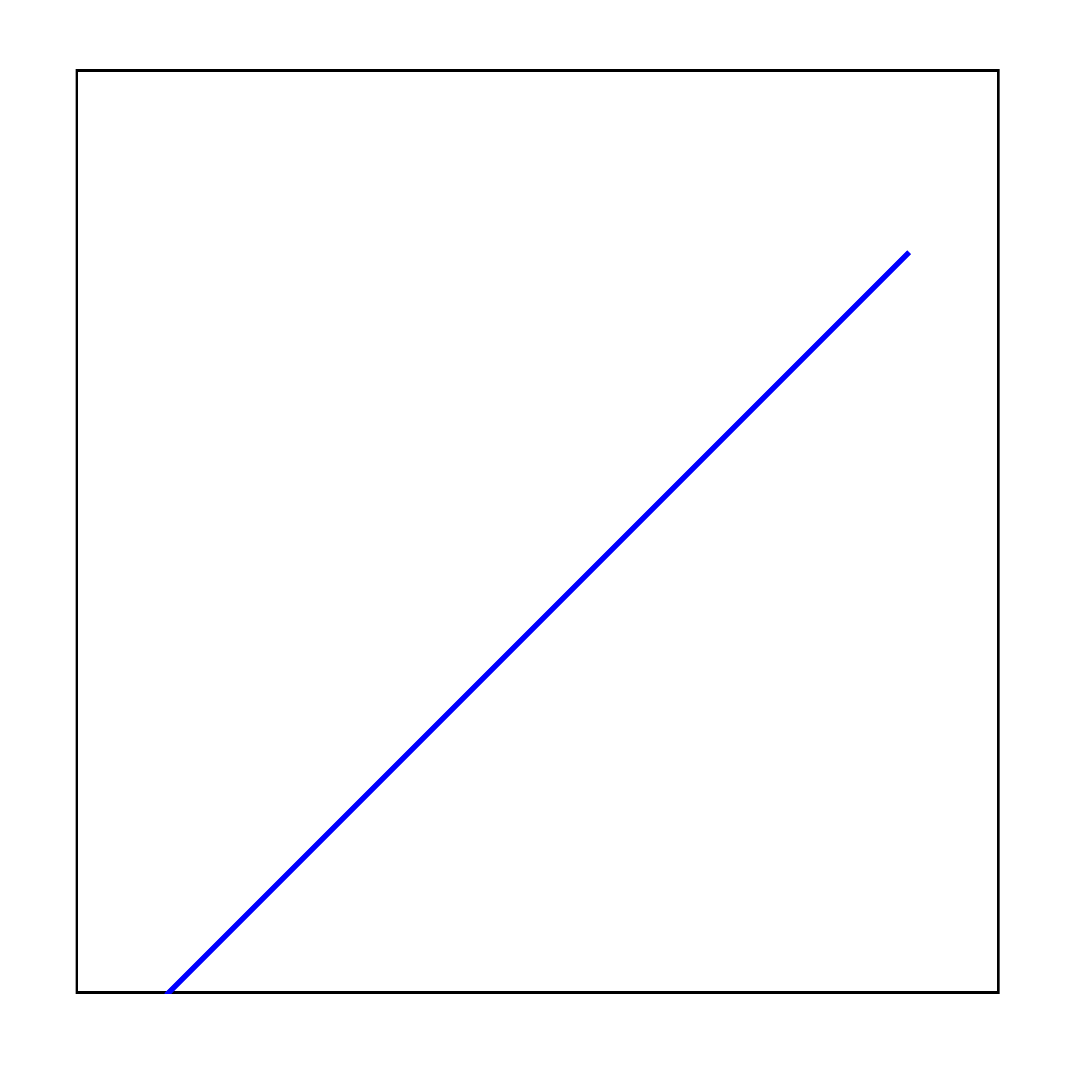}}%
    \put(0.40318924,0.93547031){\color[rgb]{0,0,0}\makebox(0,0)[lt]{\lineheight{1.25}\smash{\begin{tabular}[t]{l}out-flow\end{tabular}}}}%
    \put(0.40068498,0.0009078){\color[rgb]{0,0,0}\makebox(0,0)[lt]{\lineheight{1.25}\smash{\begin{tabular}[t]{l}in-flow\end{tabular}}}}%
    \put(0.9990922,0.40283447){\color[rgb]{0,0,0}\rotatebox{90}{\makebox(0,0)[lt]{\lineheight{1.25}\smash{\begin{tabular}[t]{l}no-flow\end{tabular}}}}}%
    \put(0.00090779,0.64233492){\color[rgb]{0,0,0}\rotatebox{-90}{\makebox(0,0)[lt]{\lineheight{1.25}\smash{\begin{tabular}[t]{l}no-flow\end{tabular}}}}}%
    \put(0.50403155,0.36449245){\color[rgb]{0,0,0}\makebox(0,0)[lt]{\lineheight{1.25}\smash{\begin{tabular}[t]{l}$\gamma$\end{tabular}}}}%
    \put(0.12781315,0.81059003){\color[rgb]{0,0,0}\makebox(0,0)[lt]{\lineheight{1.25}\smash{\begin{tabular}[t]{l}$\Omega$\end{tabular}}}}%
  \end{picture}%
\endgroup%
}
    \caption{Domain $\Omega$ and fracture $\gamma$ for the example of Subsection
    \ref{subsec:example_single_fracture}.}%
    \label{fig:example_single_fracture_domain}
\end{figure}

We consider a single fracture in the unit square domain $\Omega = (0,
1)^2$. The fracture starts from $(0.1, 0)$ and ends at $(0.9, 0.8)$, see Figure
\ref{fig:example_single_fracture_domain} for a schematic representation of the computational domain.
As shown in the figure, we assume that the left and right boundaries are
impervious while on the bottom and top we set in-flow and out-flow conditions,
respectively.
In this example the values of data are artificial and meant to highlight some phenomena, we thus omit units of measure to avoid confusions.
We consider two different settings, the common data are defined in Table
\ref{tab:single_fracture}.

\begin{table}
    \centering
    \begin{tabular}{|c|c|c|c|c|c|}
        \hline
        $\eta_\Omega = 0.5$ & $\phi_{\Omega, 0}=0.2$ & $\eta_\gamma = 2$ &
        $\epsilon_0=10^{-2}$ & $k_0=1$ \\ \hline
        $\mu=1$ &
        $f=0$ & $q_{\partial \Omega}=0$ &
        $p_{\partial \Omega}^{\rm out-flow} = 0$ & $p_{\partial\Omega}^{\rm
        in-flow}=1$\\\hline $k_{\gamma, 0} = 10^2$ & $\kappa_{\gamma, 0} =10^2$ &
        $\kappa_{\gamma, 0} = 10^2$ & $\mu=1$ & $f_\gamma=0$\\\hline
        $q_{\partial \gamma} = 0$ & $p_{\partial \gamma}^{\rm
        out-flow}=0$ & $p_{\partial \gamma}^{\rm in-flow}=1$ &
        $\lambda_w=1$ & $\lambda_s=10^{-1}$\\\hline
        $\rho_w c_w = 1$ & $\rho_s c_s = 1$ & $j=0$ & $\theta_0=1$ &
        $\tau_{\partial \Omega} = 0$\\\hline
        $\theta_{\partial \Omega}^{\rm out-flow}=0$ & $\theta_{\partial \Omega}^{\rm
        in-flow} = 1.5$ &
        $j_\gamma=0$ & $\theta_{\gamma, 0}=1$ & $\tau_{\partial \gamma} = 0$\\\hline
        $\theta_{\partial \gamma}^{\rm out-flow}=0$ & $\theta_{\partial
        \gamma}^{\rm in-flow}=1.5$ &
        $d=1$ & $u_0=0$ & $\chi_{\partial \Omega}=0$\\\hline
        $u_{\partial \Omega}^{\rm out-flow}=0$ &
        $d_\gamma=10^{-1}$ & $\delta_\gamma=10^{-1}$ & $u_{\gamma, 0}=0$ &
        $\chi_{\partial \gamma}=0$\\\hline
        $u_{\partial \gamma}^{\rm out-flow}=0$
        &$\lambda^-(\theta)=10e^{-\frac{4}{\theta}}$ & $r(u) = u^2$ & & \\\hline
    \end{tabular}
    \caption{Common data for the examples in Subsection
    \ref{subsec:example_single_fracture} and Subsection
    \ref{subsec:example_multiple_fracture}.}%
    \label{tab:single_fracture}
\end{table}
The porous medium is discretized with approximately 10000 triangles and the
fracture with approximately 75 segments. Since the computational grid is made of
triangles we can simply employ the lowest order Raviart-Thomas mixed finite element method for
the solution of the flow problem.

In the solution plots presented in the sequel, for each time step, we
represent on the top left the pressure and Darcy velocity; on the top right
the temperature; on the bottom left the solute multiplied by the porosity or
aperture; on the bottom centre the
precipitate multiplied by the porosity or aperture;
and on the bottom right the porosity and aperture. Note that the
values, apart from the velocity, are rescaled at each time step to highlight
some details of the solution. {Moreover, we choose to represent $\phi u$, $\phi
w$ since these quantities, which correspond to the amount of solute/precipitate
per unit rock volume, are easier to interpret when porosity (and aperture)
changes significantly.}

\subsubsection{Solute injection}\label{subsubsec:example_single_fracture_case_i}

In this first case we consider the following additional data: the final
simulation time is $T = 3$ discretized with 60 time steps;
$w_0=0.3$ and $u_{\partial \Omega}^{\rm inflow} = 2$ for
equation \eqref{eq:problem_solute_precipitate}; $w_{\gamma, 0} = 0.3$ and
$u_{\partial \gamma}^{\rm inflow} = 2$ for \eqref{eq:reduced_solute_precipitate}.
\begin{figure}[p]
    \centering
    \subfloat[Solution at time $0.05$ and time step $n=1$.]
    {\includegraphics[width=.75\textwidth]{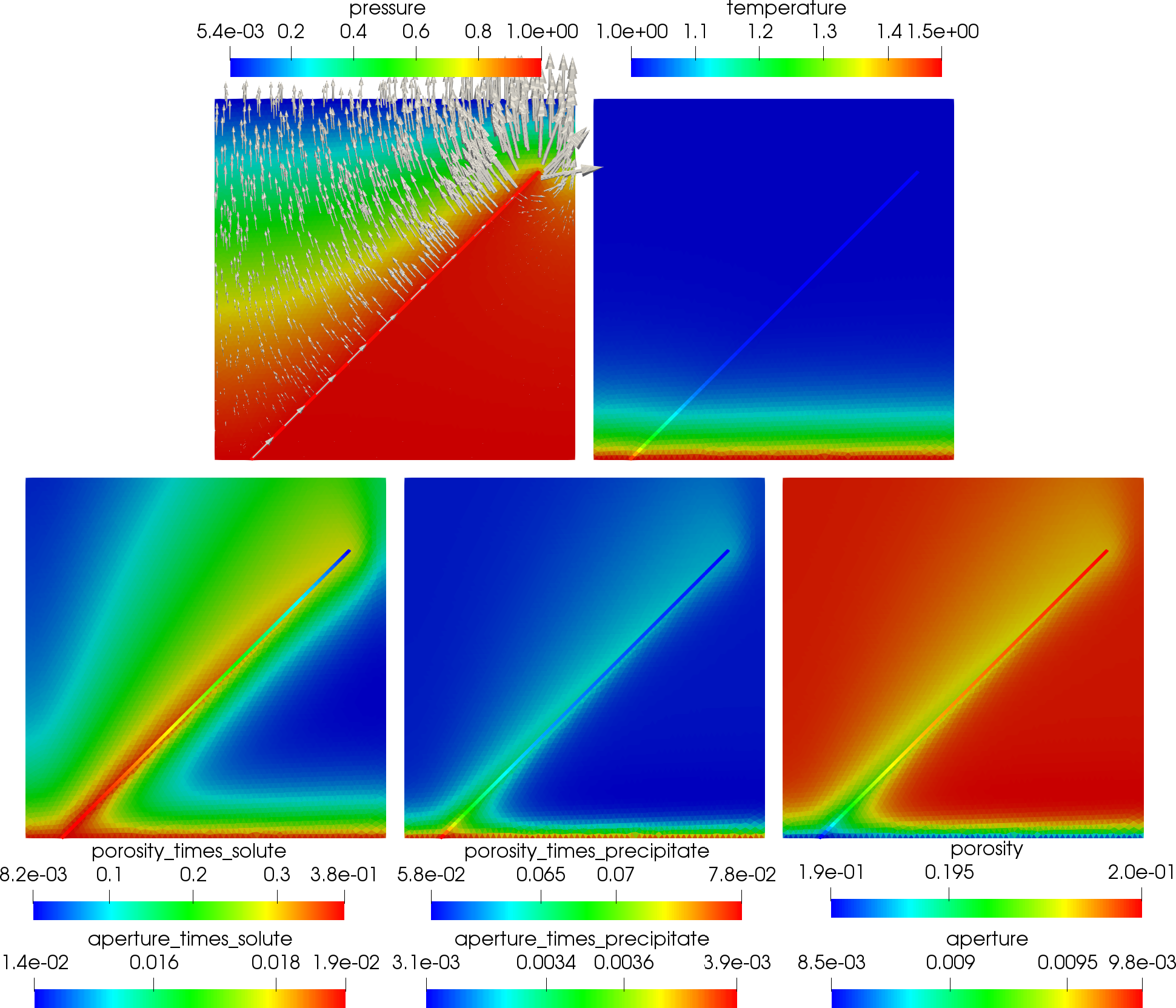}}\\
    \subfloat[Solution at time $2.25$ and time step $n=45$.]
    {\includegraphics[width=.75\textwidth]{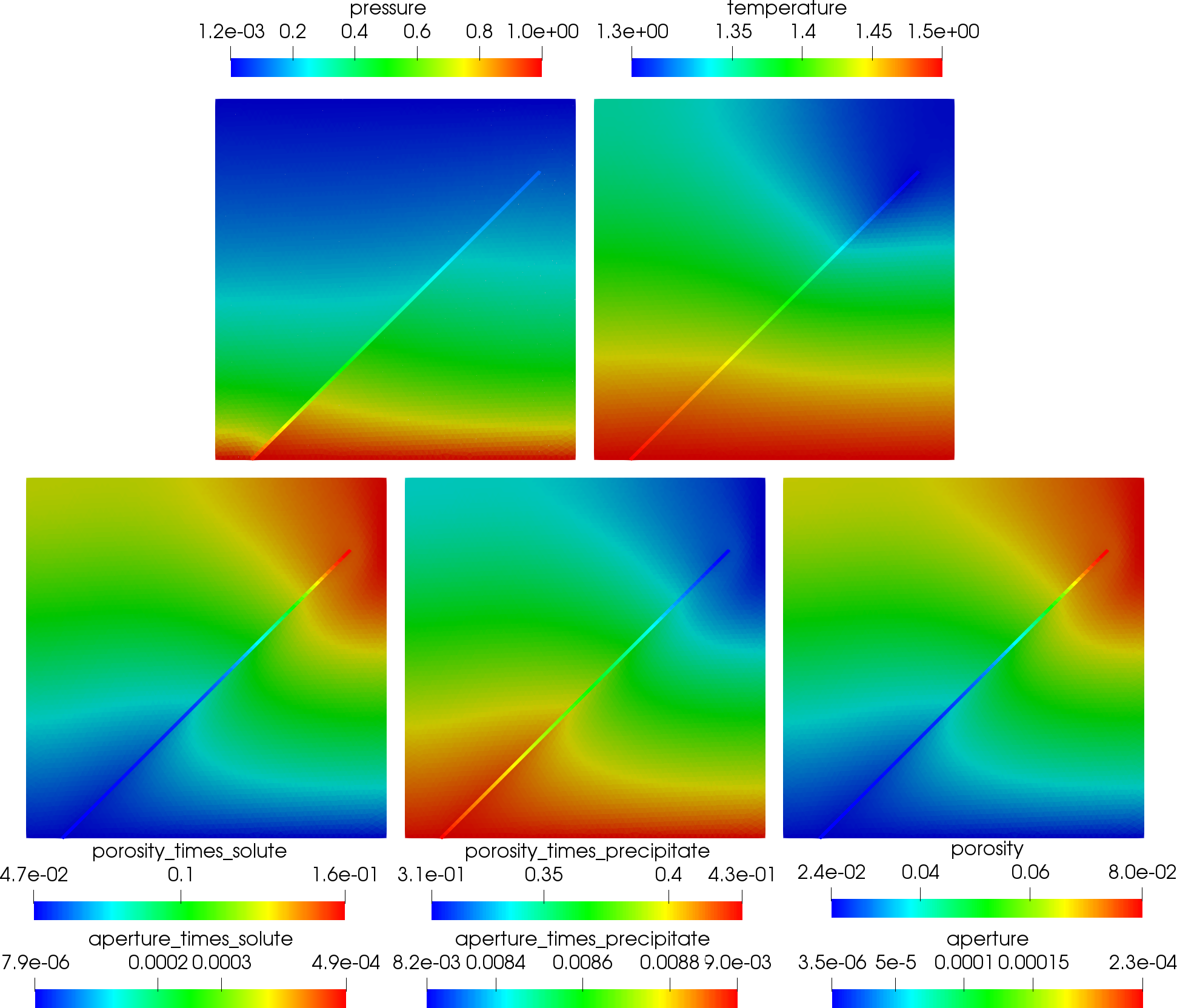}}
    \caption{Numerical solutions of the example in section
    \ref{subsubsec:example_single_fracture_case_i}.}
    \label{fig:solution_example_1_a}
\end{figure}

In this example, we inject in the porous media from the bottom boundary the
solute as well as warm water. As we see in Figure
\ref{fig:solution_example_1_a}, the fracture is highly conductive thus the
solute is transported quickly and starts to react first in the fracture.
The creation of new precipitate tends to block the fracture making it less and
less permeable. The high temperature front accelerates this process making the
fracture almost completely clogged, see the plot of fracture aperture at the bottom of Figure
\ref{fig:solution_example_1_a}. Moreover, also the pressure profile exhibits a
jump across the fracture and the Darcy velocity drops substantially,
both signs of a low fracture permeability. We  also see that the
precipitate accumulates preferably in the fracture. The porosity of the
medium decreases too, first nearby the fracture and after, due to the injection
of solute and temperature front, in the whole porous media.

This example, even with simple data, is able to capture interesting phenomena,
like the clogging of a high permeable fracture.

\subsubsection{Fracture opening}\label{subsubsec:example_single_fracture_case_ii}

In this second case we consider the following additional data: the final
simulation time is $T = 5$ discretized with 100 time steps;
by defining the square $S = \{(x, y) \in \Omega: 0.4
\leq (x, y) \leq 0.6 \}$ we have
$w_0=1$ in $S$ and $w_0=0$ elsewhere and $u_{\partial \Omega}^{\rm in-flow} = 0$ for
equation \eqref{eq:problem_solute_precipitate}; $w_{\gamma, 0} = 1$ in $\gamma
\cap S$ and $w_{\gamma, 0}=0$ otherwise and
$u_{\partial \gamma}^{\rm in-flow} = 0$ for \eqref{eq:reduced_solute_precipitate}.
\begin{figure}[p]
    \centering
    \subfloat[Solution at time $1.15$ and time step $n=23$.]
    {\includegraphics[width=.75\textwidth]{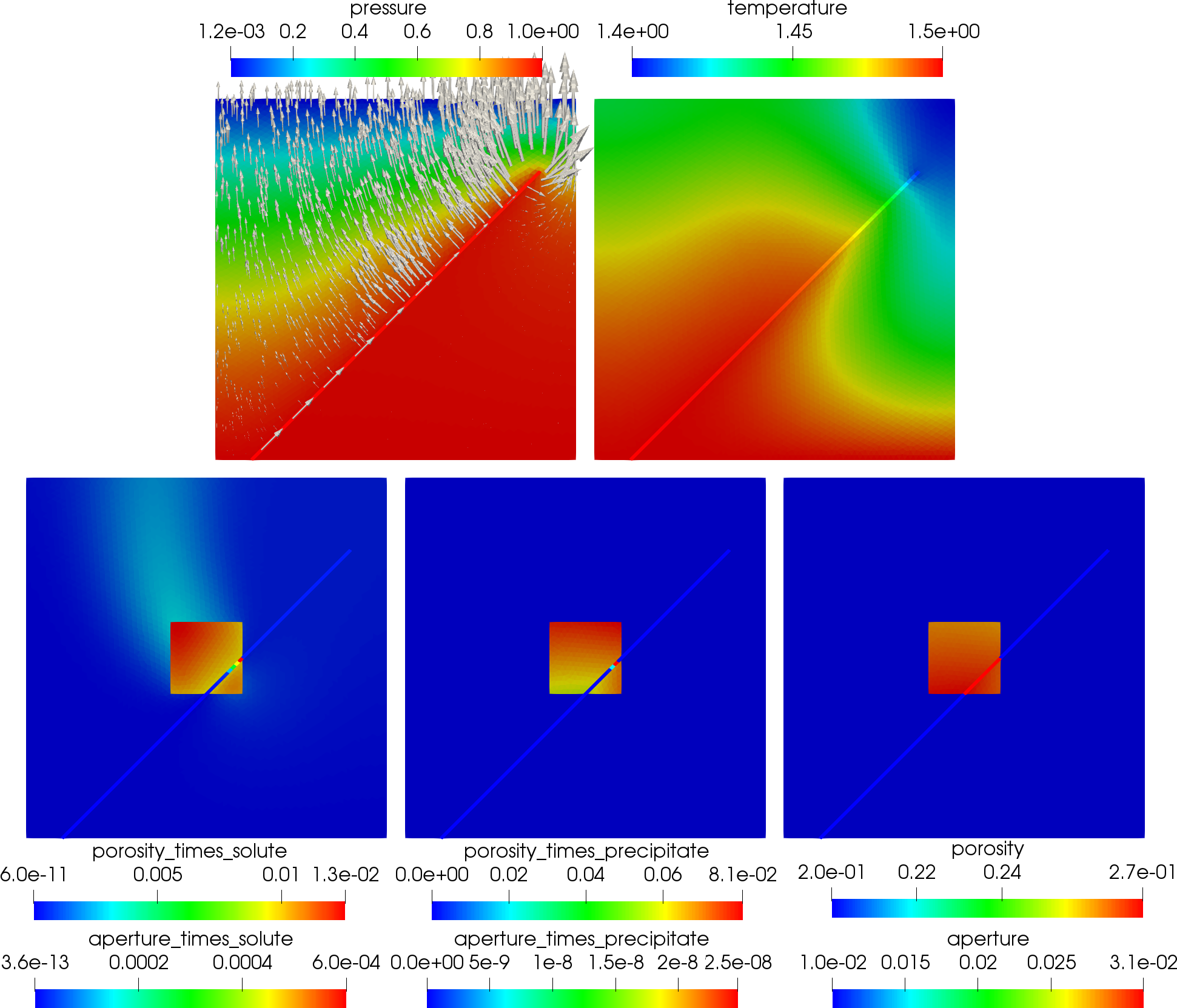}}\\
    \subfloat[Solution at time $1.55$ and time step $n=31$.]
    {\includegraphics[width=.75\textwidth]{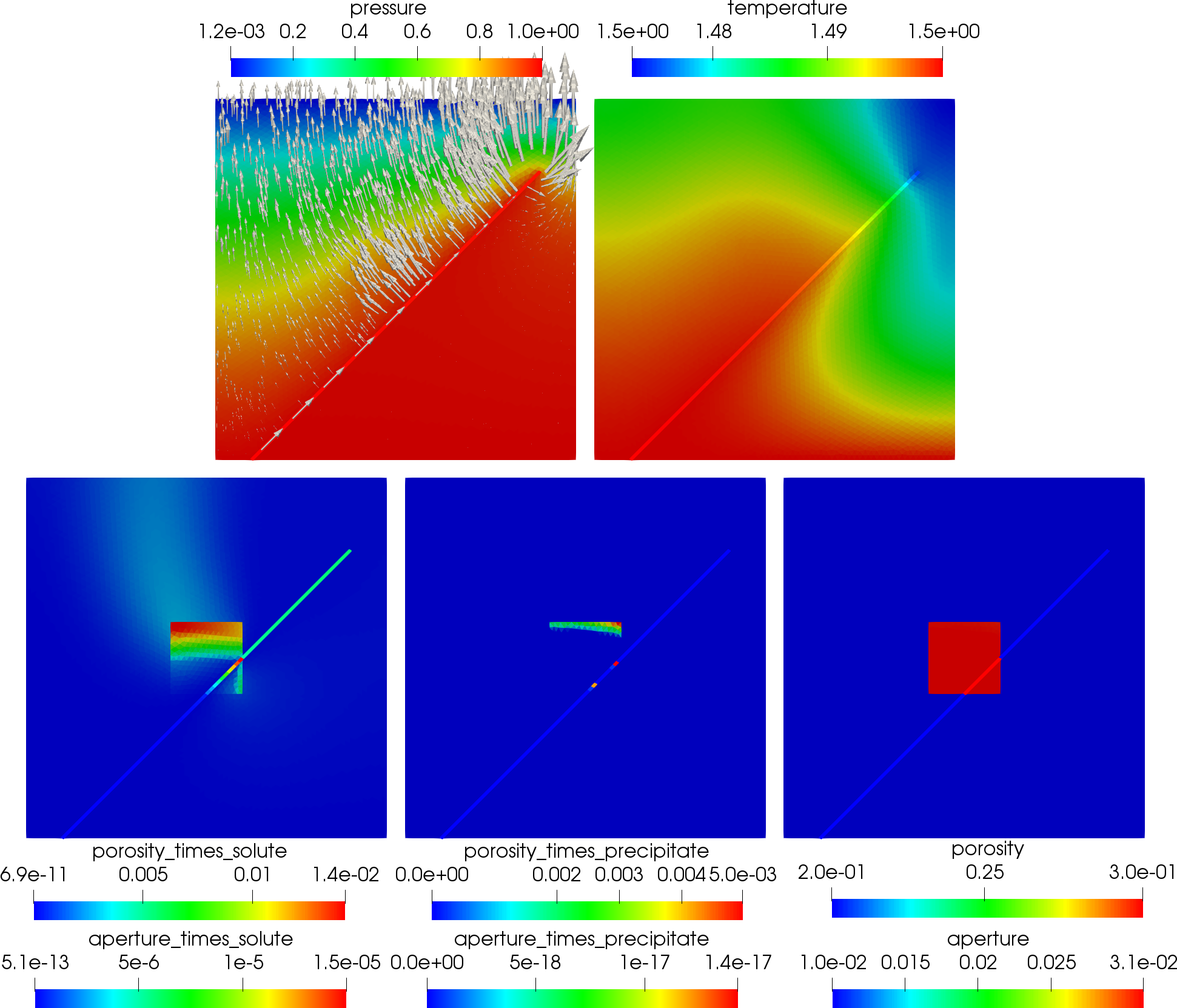}}
    \caption{Numerical solutions of the example in section
    \ref{subsubsec:example_single_fracture_case_i}.}
    \label{fig:solution_example_1_b}
\end{figure}

In this example, we inject clean water (no solute) in the system which has a block of
precipitate in the middle. The precipitate is also contained in the fracture.
Since the fracture has high permeability the clean water dissolves first the
precipitate in the fracture by forming new solute that is released in the
system, as the top of Figure \ref{fig:solution_example_1_b} shows. The fracture
aperture is thus increased as well as the porosity. After some time, as reported
in the bottom of Figure \ref{fig:solution_example_1_b}, the high
temperature front reaches the block of precipitate also through
the porous media. The reaction now becomes faster and dissolves most of the
precipitate forming more solute which is transported upward. The pressure
profile changes only slightly during the simulation but the Darcy velocity
increases visibly in $S$ due to the enhancement of the permeability.

Also in this case, even if the problem setting is simple, we are able to
reproduce important phenomena, like the opening of a fracture due to injection
of clean water.

\subsection{Multiple fractures}\label{subsec:example_multiple_fracture}

\begin{figure}[tb]
    \centering
    \resizebox{0.33\textwidth}{!}{\fontsize{1cm}{2cm}\selectfont
\begingroup%
  \makeatletter%
  \providecommand\color[2][]{%
    \errmessage{(Inkscape) Color is used for the text in Inkscape, but the package 'color.sty' is not loaded}%
    \renewcommand\color[2][]{}%
  }%
  \providecommand\transparent[1]{%
    \errmessage{(Inkscape) Transparency is used (non-zero) for the text in Inkscape, but the package 'transparent.sty' is not loaded}%
    \renewcommand\transparent[1]{}%
  }%
  \providecommand\rotatebox[2]{#2}%
  \newcommand*\fsize{\dimexpr\f@size pt\relax}%
  \newcommand*\lineheight[1]{\fontsize{\fsize}{#1\fsize}\selectfont}%
  \ifx\svgwidth\undefined%
    \setlength{\unitlength}{311.96932983bp}%
    \ifx\svgscale\undefined%
      \relax%
    \else%
      \setlength{\unitlength}{\unitlength * \real{\svgscale}}%
    \fi%
  \else%
    \setlength{\unitlength}{\svgwidth}%
  \fi%
  \global\let\svgwidth\undefined%
  \global\let\svgscale\undefined%
  \makeatother%
  \begin{picture}(1,0.98417801)%
    \lineheight{1}%
    \setlength\tabcolsep{0pt}%
    \put(0,0){\includegraphics[width=\unitlength,page=1]{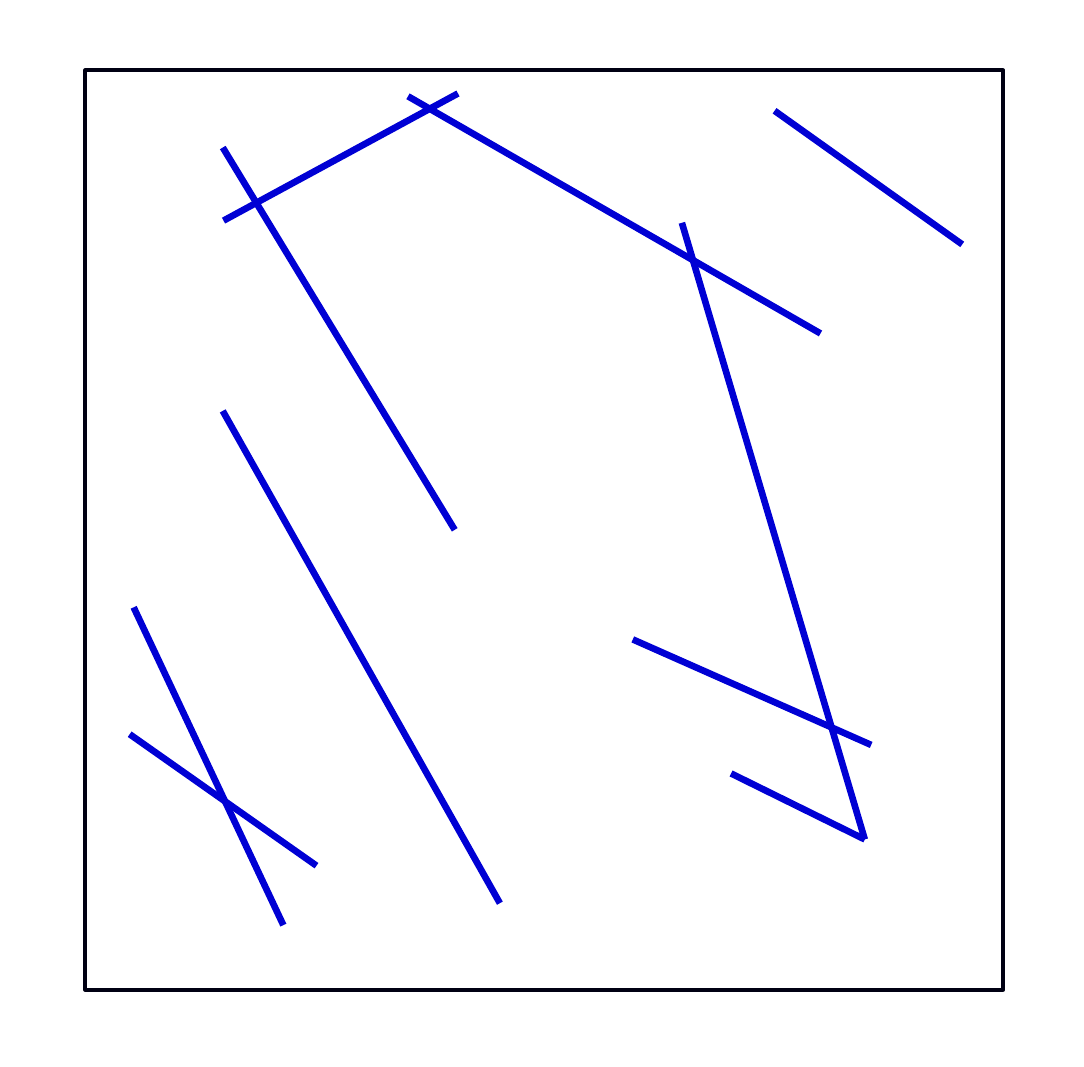}}%
    \put(0.40318921,0.93547031){\color[rgb]{0,0,0}\makebox(0,0)[lt]{\lineheight{1.25}\smash{\begin{tabular}[t]{l}out-flow\end{tabular}}}}%
    \put(0.40068496,0.0009078){\color[rgb]{0,0,0}\makebox(0,0)[lt]{\lineheight{1.25}\smash{\begin{tabular}[t]{l}in-flow\end{tabular}}}}%
    \put(0.9990922,0.40283447){\color[rgb]{0,0,0}\rotatebox{90}{\makebox(0,0)[lt]{\lineheight{1.25}\smash{\begin{tabular}[t]{l}no-flow\end{tabular}}}}}%
    \put(0.00090779,0.64233497){\color[rgb]{0,0,0}\rotatebox{-90}{\makebox(0,0)[lt]{\lineheight{1.25}\smash{\begin{tabular}[t]{l}no-flow\end{tabular}}}}}%
  \end{picture}%
\endgroup%
}%
    \hspace*{0.1\textwidth}%
    \resizebox{0.33\textwidth}{!}{\fontsize{1cm}{2cm}\selectfont
\begingroup%
  \makeatletter%
  \providecommand\color[2][]{%
    \errmessage{(Inkscape) Color is used for the text in Inkscape, but the package 'color.sty' is not loaded}%
    \renewcommand\color[2][]{}%
  }%
  \providecommand\transparent[1]{%
    \errmessage{(Inkscape) Transparency is used (non-zero) for the text in Inkscape, but the package 'transparent.sty' is not loaded}%
    \renewcommand\transparent[1]{}%
  }%
  \providecommand\rotatebox[2]{#2}%
  \newcommand*\fsize{\dimexpr\f@size pt\relax}%
  \newcommand*\lineheight[1]{\fontsize{\fsize}{#1\fsize}\selectfont}%
  \ifx\svgwidth\undefined%
    \setlength{\unitlength}{311.96932983bp}%
    \ifx\svgscale\undefined%
      \relax%
    \else%
      \setlength{\unitlength}{\unitlength * \real{\svgscale}}%
    \fi%
  \else%
    \setlength{\unitlength}{\svgwidth}%
  \fi%
  \global\let\svgwidth\undefined%
  \global\let\svgscale\undefined%
  \makeatother%
  \begin{picture}(1,0.98417801)%
    \lineheight{1}%
    \setlength\tabcolsep{0pt}%
    \put(0,0){\includegraphics[width=\unitlength,page=1]{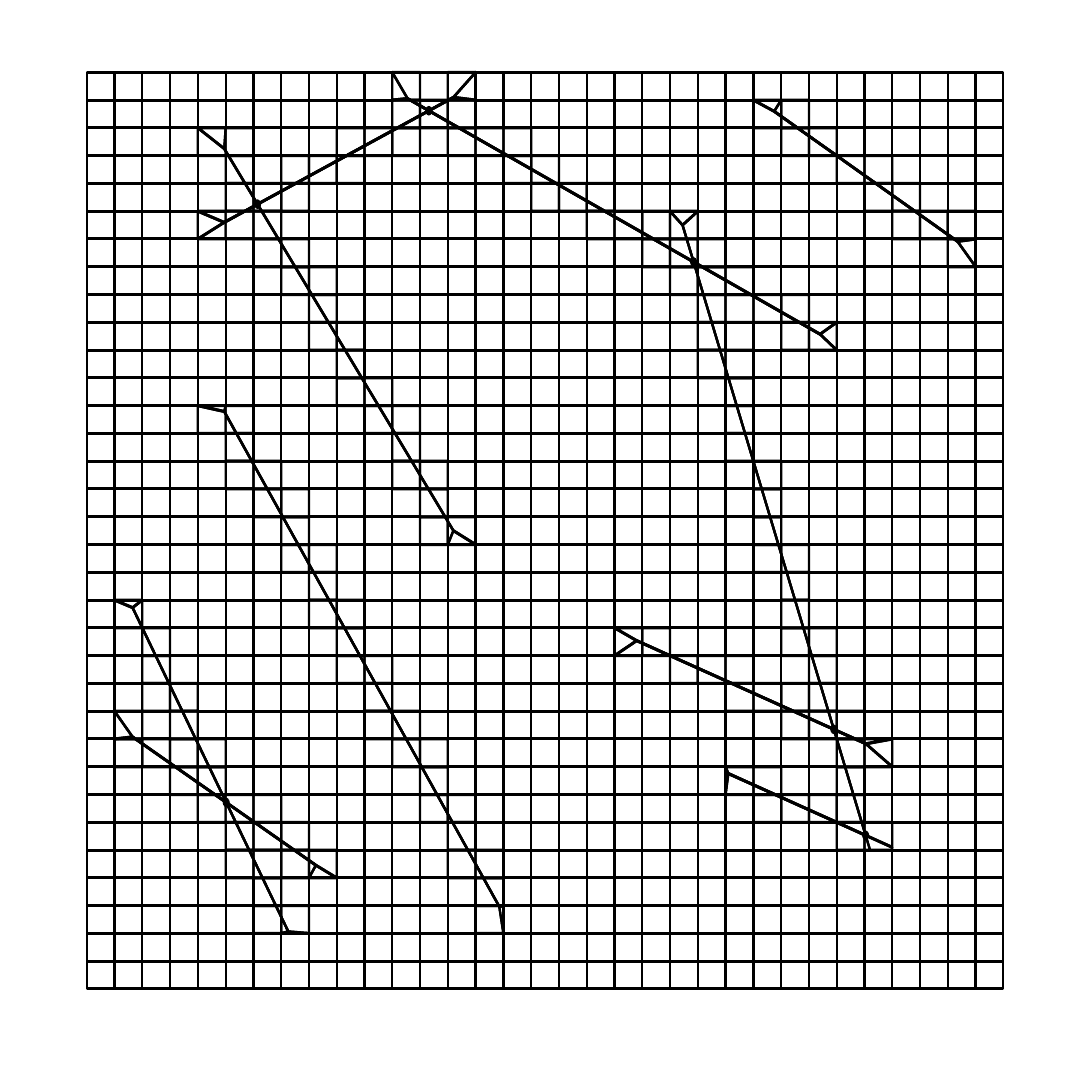}}%
    \put(0.40318921,0.93547031){\color[rgb]{0,0,0}\makebox(0,0)[lt]{\lineheight{1.25}\smash{\begin{tabular}[t]{l}out-flow\end{tabular}}}}%
    \put(0.40068496,0.0009078){\color[rgb]{0,0,0}\makebox(0,0)[lt]{\lineheight{1.25}\smash{\begin{tabular}[t]{l}in-flow\end{tabular}}}}%
    \put(0.9990922,0.40283447){\color[rgb]{0,0,0}\rotatebox{90}{\makebox(0,0)[lt]{\lineheight{1.25}\smash{\begin{tabular}[t]{l}no-flow\end{tabular}}}}}%
    \put(0.00090779,0.64233497){\color[rgb]{0,0,0}\rotatebox{-90}{\makebox(0,0)[lt]{\lineheight{1.25}\smash{\begin{tabular}[t]{l}no-flow\end{tabular}}}}}%
  \end{picture}%
\endgroup%
}
    \caption{On the left, domain $\Omega$ and fracture $\gamma$ for the example of Subsection
    \ref{subsec:example_multiple_fracture}. On the right, the computation grid with
    the cut elements.}%
    \label{fig:example_multiple_fracture_domain}
\end{figure}

In this third case we increase the geometrical complexity by considering a
problem with 10 intersecting fractures. The geometry of the fractures is
taken from the Benchmark 3 of \cite{Flemisch2016a} but with different properties.
The domain $\Omega = (0, 1)^2$ and the fractures are shown in Figure
\ref{fig:example_multiple_fracture_domain}.
For the computational grid, due to the complexity of the fracture network we consider
the procedure discussed in \cite{Fumagalli2020b}: first a Cartesian grid is
constructed and then its cells are cut if a fracture is crossing. Since we may obtain
cells of arbitrary shape, we adopt the lowest order mixed virtual element method for the
discretization of the flow equations. The grid is represented on the right of Figure
\ref{fig:example_multiple_fracture_domain}, the fracture network is discretized
with approximately 200 segments while the porous media grid consists of approximately 1300
cells.

We consider two different settings, in the first the solute is injected from the
bottom and starts to react while transported upward. In the second, the
fractures are filled with precipitate and clean water is injected reacting and
opening the fractures. The data are the same as in the previous case and
reported in Table \ref{tab:single_fracture}, with the
only exception of  $\eta_\gamma = 4$ to emphasize aperture changes.

In this example the solution plots follow the same organization as before.

\subsubsection{Solute injection}\label{subsubsec:example_multiple_fracture_case_i}

In this first case the additional data are the same as in section \ref{subsubsec:example_single_fracture_case_i}. The end time is now set
as $T=2.5$ discretized with 50 time steps.
\begin{figure}[p]
    \centering
    \subfloat[Solution at time $1.25$ and time step $n=25$.]
    {\includegraphics[width=.75\textwidth]{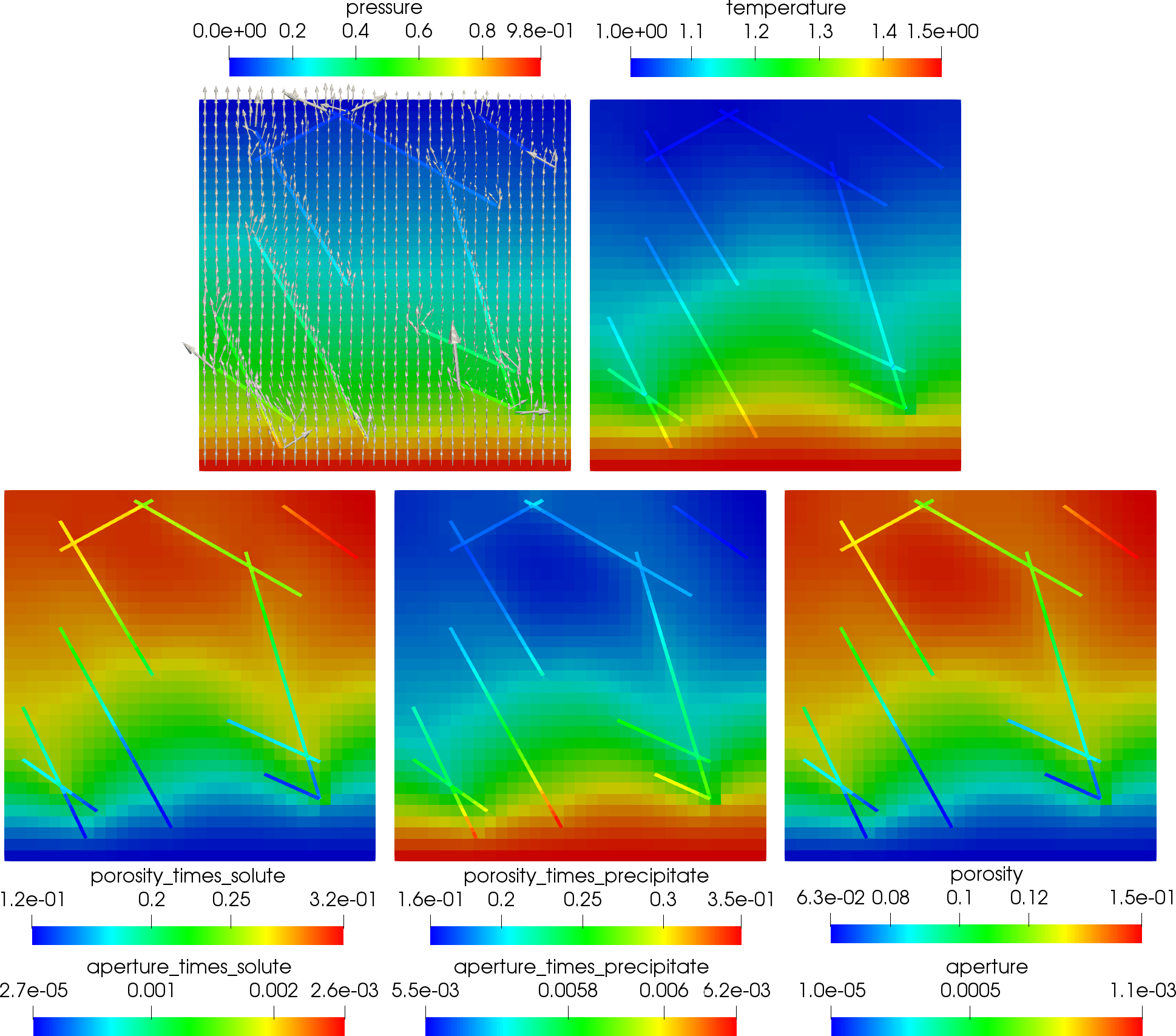}}\\
    \subfloat[Solution at time $2.25$ and time step $n=45$.]
    {\includegraphics[width=.75\textwidth]{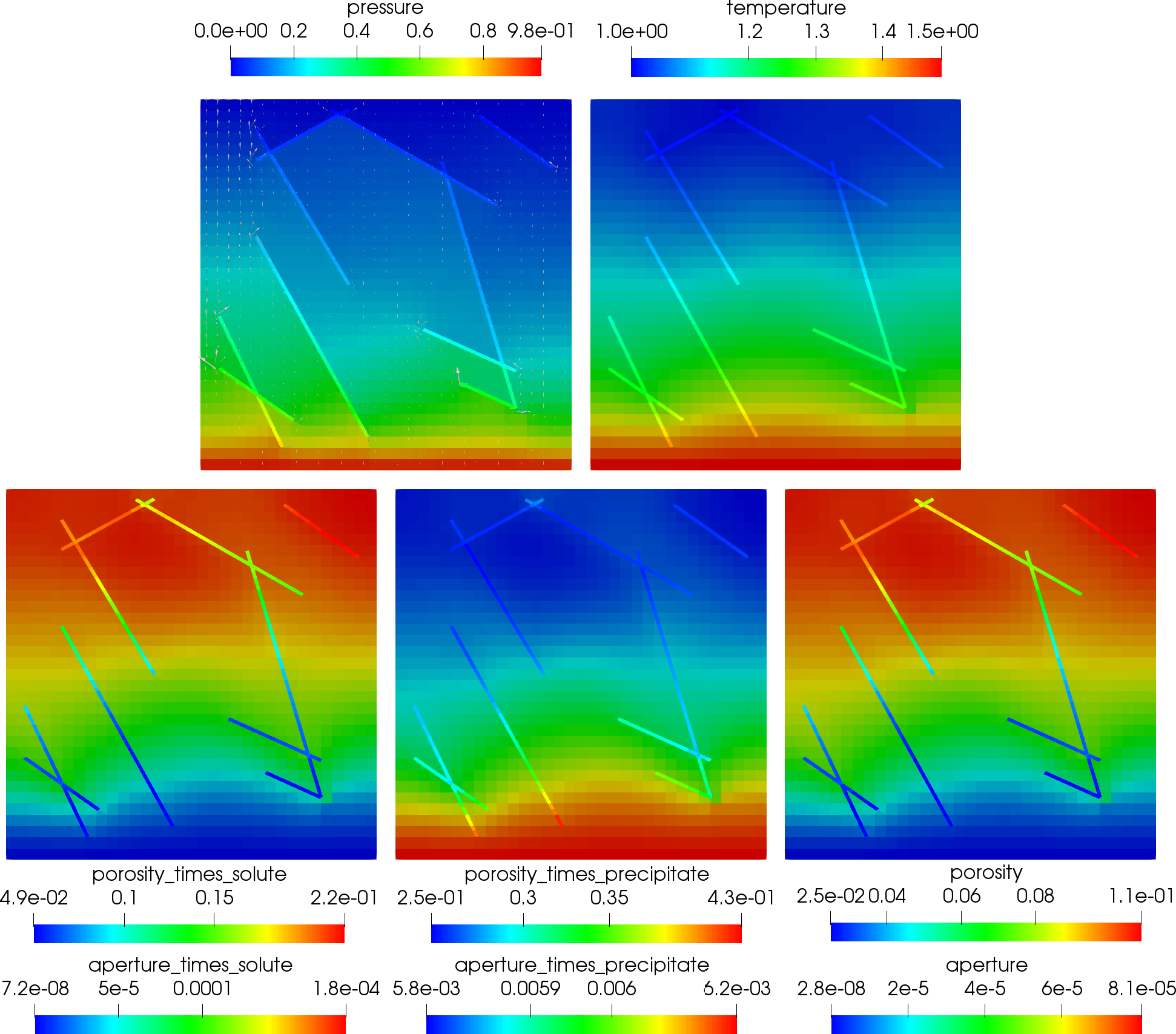}}
    \caption{Numerical solutions of the example in section
    \ref{subsubsec:example_multiple_fracture_case_i}.}
    \label{fig:solution_example_2_a}
\end{figure}

In this example, we inject solute from the bottom of the system along with hot
water. Figure \ref{fig:solution_example_2_a} reports the solutions obtained for
different time steps. As the solute flows in the porous media and fractures it
reacts forming new precipitate. The precipitate occludes the pores and attaches
to the fracture walls, reducing the porosity and fracture aperture. The latter
is reduced more rapidly since the fractures, in the beginning, are more
permeable than the surrounding medium. At the second time step many fractures are
 already almost occluded with an impact on the Darcy velocity which drops and can barely
detected.

This example shows the interaction between a complex network and the occlusion
of the fractures. This in an interesting phenomena that the model is able to
reproduce.

\subsubsection{Opening fractures}\label{subsubsec:example_multiple_fracture_case_ii}

In this case most of the additional data are the same as in
section \ref{subsubsec:example_single_fracture_case_ii}. However, now
the fractures have aperture $\epsilon=10^{-4}$ with an initial precipitate
$w_{\gamma, 0} = 10$. In the porous media the initial precipitate is set to
$w_{0} = 0$.
\begin{figure}[p]
    \centering
    \subfloat[Solution at time $1.3$ and time step $n=26$.]
    {\includegraphics[width=.75\textwidth]{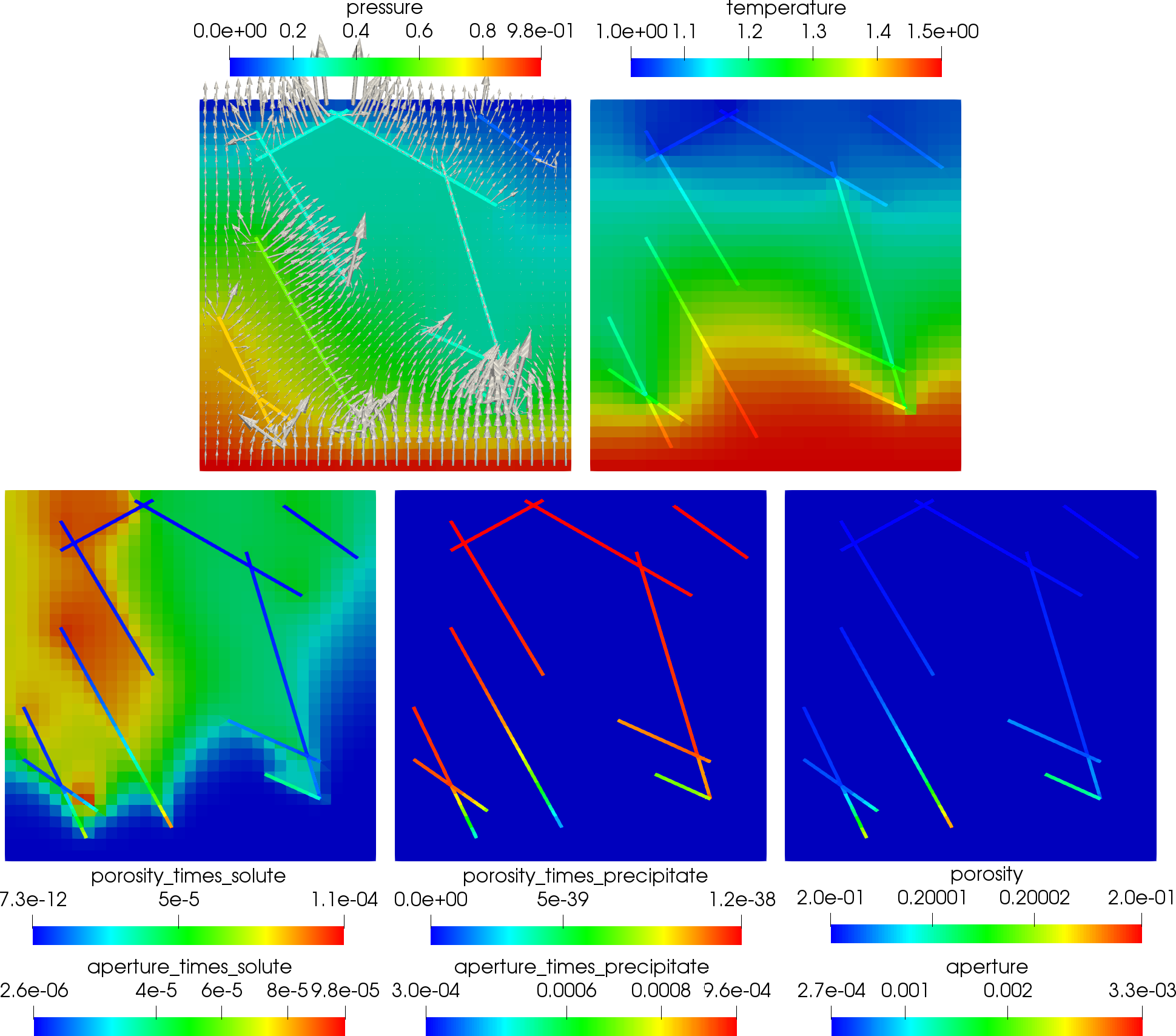}}\\
    \subfloat[Solution at time $3.25$ and time step $n=65$.]
    {\includegraphics[width=.75\textwidth]{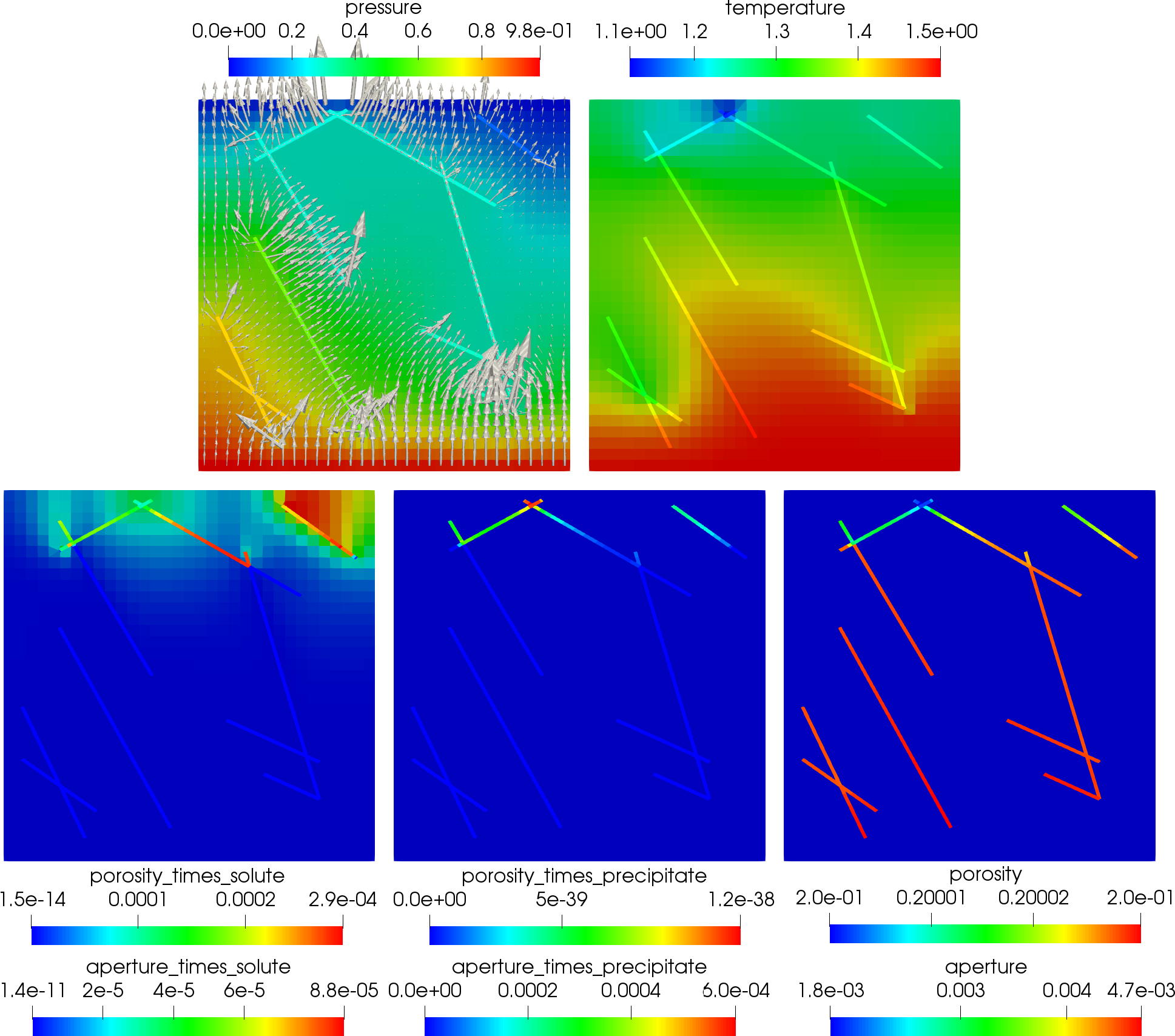}}
    \caption{Numerical solutions of the example in section
    \ref{subsubsec:example_multiple_fracture_case_ii}.}
    \label{fig:solution_example_2_b}
\end{figure}

The numerical solution of this example is represented in Figure
\ref{fig:solution_example_2_b}. We see that in the first time step shown, fractures still have small aperture in most parts, but it starts to increase
as the high temperature front advances. At the same time the precipitate dissolves and
new solute is created and transported upwards. In the second block of
figures, most of the fractures have substantially increased their aperture and
the precipitate is now only present at the top of the domain. Even if not
represented, at the end of the simulation all the precipitate has reacted and
all the solute has been transported away from the domain. The aperture reaches
a stable value of approximately $4.5 \cdot 10^{-3}$ for all fractures.

The possibility to open and thus enhance the flow property of a geothermal
system is crucial in its exploitation. This simple example showed that the model
is also able to capture such scenario.

\section{Conclusion}\label{sec:conclusion}

In this work we have presented a mathematical model to describe in a  fractured
porous media a reactive single-phase flow with thermal effects. The mathematical
model consists of a set of coupled PDEs and ODEs for the evolution of pressure,
Darcy velocity, porosity, temperature, solute and precipitate concentration.
Moreover, this model considers fractures as one co-dimensional manifolds in the
framework of the so-called mixed-dimensional setting, so that fracture aperture
is not any more a geometrical constraint but a model parameter and it can thus
freely vary during the simulation. The resulting system is fully coupled,
non-linear and generally non-smooth due to the modeling of the chemical
reactions involved. To numerically solve the full system, we have introduced a
temporal splitting scheme so that each physical process is solved sequentially
but ensuring, {as proven only experimentally by means of numerical tests}, the
mass conservation of the system. A nested splitting is used to compute the
solute concentration separating the advective-diffusive parts from the reaction.
For the latter a system of discontinuous ordinary differential equations is
solved by using an event driven approach. To conclude, as the numerical examples
showed the presented model is able to capture interesting and physically
relevant phenomena also in presence of complex fracture networks.


\begin{thebibliography}{10}

\bibitem{Agosti2015}
Abramo Agosti, Luca Formaggia, and Anna Scotti.
\newblock Analysis of a model for precipitation and dissolution coupled with a
  darcy flux.
\newblock {\em Journal of Mathematical Analysis and Applications},
  431(2):752--781, 2015.

\bibitem{Agosti2016}
Abramo Agosti, Bianca Giovanardi, Luca Formaggia, and Anna Scotti.
\newblock A numerical procedure for geochemical compaction in the presence of
  discontinuous reactions.
\newblock {\em Advances in Water Resources}, 94:332 -- 344, 2016.

\bibitem{Ahmed2018}
Elyes Ahmed, Alessio Fumagalli, and Ana Budi{\v{s}}a.
\newblock A multiscale flux basis for mortar mixed discretizations of reduced
  darcy-forchheimer fracture models.
\newblock {\em Computer Methods in Applied Mechanics and Engineering},
  354:16--36, 2019.

\bibitem{Ahmed2017}
Elyes Ahmed, J{\'e}r{\^o}me Jaffr{\'e}, and Jean~E. Roberts.
\newblock A reduced fracture model for two-phase flow with different rock
  types.
\newblock {\em Mathematics and Computers in Simulation}, 137:49--70, 2017.
\newblock MAMERN VI-2015: 6th International Conference on Approximation Methods
  and Numerical Modeling in Environment and Natural Resources.

\bibitem{Alboin2002}
Clarisse Alboin, J{\'e}r{\^o}me Jaffr{\'e}, Jean~E. Roberts, and Christophe
  Serres.
\newblock Modeling fractures as interfaces for flow and transport in porous
  media.
\newblock In {\em Fluid flow and transport in porous media: mathematical and
  numerical treatment ({S}outh {H}adley, {MA}, 2001)}, volume 295 of {\em
  Contemp. Math.}, pages 13--24. Amer. Math. Soc., Providence, RI, 2002.

\bibitem{Alboin2000}
Clarisse Alboin, J{\'e}r{\^o}me Jaffr{\'e}, Jean~E. Roberts, Xuewen Wang, and
  Christophe Serres.
\newblock {\em Domain decomposition for some transmission problems in flow in
  porous media}, volume 552 of {\em Lecture Notes in Phys.}, pages 22--34.
\newblock Springer, Berlin, 2000.

\bibitem{Angot2003}
Philippe Angot.
\newblock A model of fracture for elliptic problems with flux and solution
  jumps.
\newblock {\em Comptes Rendus Mathematique}, 337(6):425--430, 2003.

\bibitem{Antonietti2019}
Paola~F. Antonietti, Chiara Facciol\`a, Alessandro Russo, and Marco Verani.
\newblock Discontinuous galerkin approximation of flows in fractured porous
  media on polytopic grids.
\newblock {\em SIAM Journal on Scientific Computing}, 41(1):A109--A138, 2019.

\bibitem{Antonietti}
Paola~Francesca Antonietti, Luca Formaggia, Anna Scotti, Marco Verani, and
  Nicola Verzotti.
\newblock Mimetic finite difference approximation of flows in fractured porous
  media.
\newblock {\em ESAIM: M2AN}, 50(3):809--832, 2016.

\bibitem{Bear1972}
Jacob Bear.
\newblock {\em Dynamics of {F}luids in {P}orous {M}edia}.
\newblock American Elsevier, 1972.

\bibitem{Bear1990}
Jacob Bear and Yehuda Bachmat.
\newblock {\em Introduction to Modeling of Transport Phenomena in Porous
  Media}.
\newblock Theory and Applications of Transport in Porous Media. Springer
  Netherlands, 1990.

\bibitem{BeiraoVeiga2016}
Louren\c{c}o Beir\~{a}o~da Veiga, Franco Brezzi, Luisa~Donatella Marini, and
  Alessandro Russo.
\newblock Mixed virtual element methods for general second order elliptic
  problems on polygonal meshes.
\newblock {\em ESAIM: M2AN}, 50(3):727--747, 2016.

\bibitem{Benedetto2016}
Mat{\'i}as~Fernando Benedetto, Stefano Berrone, Andrea Borio, Sandra
  Pieraccini, and Stefano Scial{\`o}.
\newblock A hybrid mortar virtual element method for discrete fracture network
  simulations.
\newblock {\em Journal of Computational Physics}, 306:148 -- 166, 2016.

\bibitem{Benedetto2014}
Mat{\'i}as~Fernando Benedetto, Stefano Berrone, Sandra Pieraccini, and Stefano
  Scial{\`o}.
\newblock The virtual element method for discrete fracture network simulations.
\newblock {\em Computer Methods in Applied Mechanics and Engineering},
  280(0):135--156, 2014.

\bibitem{Berre2020}
Inga Berre, Wietse~M. Boon, Bernd Flemisch, Alessio Fumagalli, Dennis
  Gl\"{a}ser, Eirik Keilegavlen, Anna Scotti, Ivar Stefansson, Alexandru
  Tatomir, Konstantin Brenner, Samuel Burbulla, Philippe Devloo, Omar Duran,
  Marco Favino, Julian Hennicker, I-Hsien Lee, Konstantin Lipnikov, Roland
  Masson, Klaus Mosthaf, Maria Giuseppina~Chiara Nestola, Chuen-Fa Ni, Kirill
  Nikitin, Philipp Sch\"{a}dle, Daniil Svyatskiy, Ruslan Yanbarisov, and
  Patrick Zulian.
\newblock Verification benchmarks for single-phase flow in three-dimensional
  fractured porous media.
\newblock Technical report, arXiv:2002.07005 [math.NA], 2020.

\bibitem{Berre2019b}
Inga Berre, Florian Doster, and Eirik Keilegavlen.
\newblock Flow in fractured porous media: A review of conceptual models and
  discretization approaches.
\newblock {\em Transport in Porous Media}, 130(1):215--236, 2019.

\bibitem{Boffi2013}
Daniele Boffi, Franco Brezzi, and Michel Fortin.
\newblock {\em Mixed Finite Element Methods and Applications}.
\newblock Springer Series in Computational Mathematics. Springer Berlin
  Heidelberg, 2013.

\bibitem{Boon2018}
Wietse~M. Boon, Jan~M. Nordbotten, and Ivan Yotov.
\newblock Robust discretization of flow in fractured porous media.
\newblock {\em SIAM Journal on Numerical Analysis}, 56(4):2203--2233, 2018.

\bibitem{Brezzi2014}
Franco Brezzi, Richard~S. Falk, and Donatella~Luisa Marini.
\newblock Basic principles of mixed virtual element methods.
\newblock {\em ESAIM: M2AN}, 48(4):1227--1240, 2014.

\bibitem{Chave2018}
Florent Chave, Daniele~A. Di~Pietro, and Luca Formaggia.
\newblock A hybrid high-order method for darcy flows in fractured porous media.
\newblock {\em SIAM Journal on Scientific Computing}, 40(2):A1063--A1094, 2018.

\bibitem{Chave2019}
Florent Chave, Daniele~A. Di~Pietro, and Luca Formaggia.
\newblock A hybrid high-order method for passive transport in fractured porous
  media.
\newblock {\em GEM - International Journal on Geomathematics}, 10(1):12, 2019.

\bibitem{Cote2005}
Jean C\^ot\'e and Jean-Marie Konrad.
\newblock A generalized thermal conductivity model for soils and construction
  materials.
\newblock {\em Canadian Geotechnical Journal}, 42(2):443--458, 2005.

\bibitem{Dreuzy2013}
Jean-Raynald de~Dreuzy, G{\'e}raldine Pichot, Baptiste Poirriez, and Jocelyne
  Erhel.
\newblock Synthetic benchmark for modeling flow in 3d fractured media.
\newblock {\em Computers \& Geosciences}, 50:59 -- 71, 2013.
\newblock Benchmark problems, datasets and methodologies for the computational
  geosciences.

\bibitem{DL}
Lopez~L. Dieci~L.
\newblock Sliding motion in discontinuous differential systems: Theory and a
  computational apporach.
\newblock 2008.

\bibitem{Droniou2013}
J{\'e}r{\^o}me Droniou.
\newblock Finite volume scheme for diffusion equations: Introduction to and
  review of modern methods, April 2013.

\bibitem{Eymard2000}
Robert Eymard, Thierry Gallou\"et, and Rapha\`ele Herbin.
\newblock Finite volume methods.
\newblock In P.~G. Ciarlet and J.~L. Lions, editors, {\em Solution of Equation
  in $\mathcal{R}^n$ (Part 3), Techniques of Scientific Computing (Part 3)},
  volume~7 of {\em Handbook of Numerical Analysis}, pages 713--1018. Elsevier,
  2000.

\bibitem{Faille2002}
Isabelle Faille, Eric Flauraud, Fr{\'e}d{\'e}ric Nataf, Sylvie
  P{\'e}gaz-Fiornet, Fr{\'e}d{\'e}ric Schneider, and Fran\c{c}oise Willien.
\newblock A {N}ew {F}ault {M}odel in {G}eological {B}asin {M}odelling.
  {A}pplication of {F}inite {V}olume {S}cheme and {D}omain {D}ecomposition
  {M}ethods.
\newblock In {\em Finite volumes for complex applications, {III}
  ({P}orquerolles, 2002)}, pages 529--536. Hermes Sci. Publ., Paris, 2002.

\bibitem{Faille2014a}
Isabelle Faille, Alessio Fumagalli, J{\'e}r{\^o}me Jaffr{\'e}, and
  Jean~Elisabeth Roberts.
\newblock Model reduction and discretization using hybrid finite volumes of
  flow in porous media containing faults.
\newblock {\em Computational Geosciences}, 20(2):317--339, 2016.

\bibitem{Faille2005}
Isabelle Faille, Fr{\'e}d{\'e}ric Nataf, Laurent Saas, and Fran\c{c}oise
  Willien.
\newblock {F}inite {V}olume {M}ethods on {N}on-{M}atching {G}rids with
  {A}rbitrary {I}nterface {C}onditions and {H}ighly {H}eterogeneous {M}edia.
\newblock In {\em Domain Decomposition Methods in Science and Engineering},
  volume~40 of {\em Lecture Notes in Computational Science and Engineering},
  pages 243--250. Springer Berlin Heidelberg, 2005.

\bibitem{Flemisch2016a}
Bernd Flemisch, Inga Berre, Wietse Boon, Alessio Fumagalli, Nicolas Schwenck,
  Anna Scotti, Ivar Stefansson, and Alexandru Tatomir.
\newblock Benchmarks for single-phase flow in fractured porous media.
\newblock {\em Advances in Water Resources}, 111:239--258, Januray 2018.

\bibitem{Flemisch2016}
Bernd Flemisch, Alessio Fumagalli, and Anna Scotti.
\newblock {\em A Review of the XFEM-Based Approximation of Flow in Fractured
  Porous Media}, volume~12 of {\em SEMA SIMAI Springer Series}, chapter
  Advances in Discretization Methods, pages 47--76.
\newblock Springer International Publishing, Cham, 2016.

\bibitem{Formaggia2012}
Luca Formaggia, Alessio Fumagalli, Anna Scotti, and Paolo Ruffo.
\newblock A reduced model for {D}arcy's problem in networks of fractures.
\newblock {\em {ESAIM}: {M}athematical {M}odelling and {N}umerical {A}nalysis},
  48:1089--1116, 7 2014.

\bibitem{Frih2011}
Najla Frih, Vincent Martin, Jean~Elisabeth Roberts, and Ai~Sa{\^a}da.
\newblock Modeling fractures as interfaces with nonmatching grids.
\newblock {\em Computational Geosciences}, 16(4):1043--1060, 2012.

\bibitem{Frih2008}
Najla Frih, Jean~E. Roberts, and Ali Saada.
\newblock Modeling fractures as interfaces: a model for {F}orchheimer
  fractures.
\newblock {\em Computers and Geosciences}, 12(1):91--104, 2008.

\bibitem{Faille2014}
Alessio Fumagalli and Isabelle Faille.
\newblock A double-layer reduced model for fault flow on slipping domains with
  hybrid finite volume scheme.
\newblock {\em SIAM Journal on Scientific Computing}, 77:1--26, June 2018.

\bibitem{Fumagalli2016a}
Alessio Fumagalli and Eirik Keilegavlen.
\newblock Dual virtual element method for discrete fractures networks.
\newblock {\em SIAM Journal on Scientific Computing}, 40(1):B228--B258, 2018.

\bibitem{Fumagalli2017a}
Alessio Fumagalli and Eirik Keilegavlen.
\newblock Dual virtual element methods for discrete fracture matrix models.
\newblock {\em Oil \& Gas Science and Technology - Revue d'IFP Energies
  nouvelles}, 74(41):1--17, 2019.

\bibitem{Scialo2017}
Alessio Fumagalli, Eirik Keilegavlen, and Stefano Scial\`o.
\newblock Conforming, non-conforming and non-matching discretization couplings
  in discrete fracture network simulations.
\newblock {\em Journal of Computational Physics}, 376:694--712, 2019.

\bibitem{Fumagalli2012d}
Alessio Fumagalli and Anna Scotti.
\newblock A numerical method for two-phase flow in fractured porous media with
  non-matching grids.
\newblock {\em Advances in Water Resources}, 62, Part C(0):454--464, 2013.
\newblock Computational Methods in Geologic CO2 Sequestration.

\bibitem{Fumagalli2012a}
Alessio Fumagalli and Anna Scotti.
\newblock A {R}educed {M}odel for {F}low and {T}ransport in {F}ractured
  {P}orous {M}edia with {N}on-matching {G}rids.
\newblock In Andrea Cangiani, Ruslan~L. Davidchack, Emmanuil Georgoulis,
  Alexander~N. Gorban, Jeremy Levesley, and Michael~V. Tretyakov, editors, {\em
  Numerical Mathematics and Advanced Applications 2011}, pages 499--507.
  Springer Berlin Heidelberg, 2013.

\bibitem{Fumagalli2012g}
Alessio Fumagalli and Anna Scotti.
\newblock An {E}fficient {XFEM} {A}pproximation of {D}arcy {F}lows in
  {A}rbitrarily {F}ractured {P}orous {M}edia.
\newblock {\em {O}il and {G}as {S}ciences and {T}echnologies - {R}evue d'{IFP}
  {E}nergies {N}ouvelles}, 69(4):555--564, April 2014.

\bibitem{Fumagalli2020a}
Alessio Fumagalli and Anna Scotti.
\newblock Reactive flow in fractured porous media.
\newblock In {\em Finite Volumes for Complex Applications IX proceedings}.
  Springer, 2020.
\newblock Accepted.

\bibitem{Fumagalli2020b}
Alessio Fumagalli, Anna Scotti, and Luca Formaggia.
\newblock Performances of the mixed virtual element method on complex grids for
  underground flow.
\newblock Accepted in SEMA SIMAI Springer Series. Available at arXiv:2002.11974
  [math.NA], 2020.

\bibitem{Giovanardi2015}
Bianca Giovanardi, Anna Scotti, Luca Formaggia, and Paolo Ruffo.
\newblock A general framework for the simulation of geochemical compaction.
\newblock {\em Computational Geosciences}, 19(5):1027--1046, Oct 2015.

\bibitem{Guldberg1864}
C.M. Guldberg and P.~Waage.
\newblock {\em Studies concerning affinity}.
\newblock WileyBlackwell, 1864.

\bibitem{Helmig1997}
Rainer Helmig.
\newblock {\em Multiphase flow and transport processes in the subsurface: a
  contribution to the modeling of hydrosystems.}
\newblock Springer-Verlag, Berlin, Germany, 1997.

\bibitem{Jaffre2002}
J{\'e}r{\^o}me Jaffr{\'e}, Vincent Martin, and Jean~E. Roberts.
\newblock Generalized cell-centered finite volume methods for flow in porous
  media with faults.
\newblock In {\em Finite volumes for complex applications, {III}
  ({P}orquerolles, 2002)}, pages 343--350. Hermes Sci. Publ., Paris, 2002.

\bibitem{Jaffre2011}
J{\'e}r{\^o}me Jaffr{\'e}, Mokhles Mnejja, and Jean~E. Roberts.
\newblock A discrete fracture model for two-phase flow with matrix-fracture
  interaction.
\newblock {\em Procedia Computer Science}, 4:967--973, 2011.

\bibitem{Keilegavlen2019a}
Eirik Keilegavlen, Runar Berge, Alessio Fumagalli, Michele Starnoni, Ivar
  Stefansson, Jhabriel Varela, and Inga Berre.
\newblock Porepy: An open-source software for simulation of multiphysics
  processes in fractured porous media.
\newblock Technical report, arXiv:1908.09869 [math.NA], 2019.

\bibitem{Knabner2014}
Peter Knabner and Jean~Elisabeth Roberts.
\newblock Mathematical analysis of a discrete fracture model coupling darcy
  flow in the matrix with darcy-forchheimer flow in the fracture.
\newblock {\em ESAIM: Mathematical Modelling and Numerical Analysis},
  48:1451--1472, 9 2014.

\bibitem{Knabner1995}
Peter Knabner, C.J. van Duijn, and S.~Hengst.
\newblock An analysis of crystal dissolution fronts in flows through porous
  media. part 1: Compatible boundary conditions.
\newblock {\em Advances in Water Resources}, 18(3):171--185, 1995.

\bibitem{Martin2005}
Vincent Martin, J{\'e}r{\^o}me Jaffr{\'e}, and Jean~Elisabeth Roberts.
\newblock Modeling {F}ractures and {B}arriers as {I}nterfaces for {F}low in
  {P}orous {M}edia.
\newblock {\em SIAM J. Sci. Comput.}, 26(5):1667--1691, 2005.

\bibitem{Morales2010}
Fernando Morales and Ralph~E. Showalter.
\newblock The narrow fracture approximation by channeled flow.
\newblock {\em Journal of Mathematical Analysis and Applications},
  365(1):320--331, 2010.

\bibitem{Morales2012}
Fernando Morales and Ralph~E. Showalter.
\newblock Interface approximation of darcy flow in a narrow channel.
\newblock {\em Mathematical Methods in the Applied Sciences}, 35(2):182--195,
  2012.

\bibitem{Nordbotten2018}
Jan~Martin Nordbotten, Wietse Boon, Alessio Fumagalli, and Eirik Keilegavlen.
\newblock Unified approach to discretization of flow in fractured porous media.
\newblock {\em Computational Geosciences}, 23(2):225--237, 2019.

\bibitem{Nordbotten2011}
Jan~Martin Nordbotten and Micheal~A. Celia.
\newblock {\em Geological Storage of CO2: Modeling Approaches for Large-Scale
  Simulation}.
\newblock Wiley, 2011.

\bibitem{Raviart1977}
Pierre-Arnaud Raviart and Jean-Marie Thomas.
\newblock A mixed finite element method for second order elliptic problems.
\newblock {\em Lecture Notes in Mathematics}, 606:292--315, 1977.

\bibitem{Roberts1991}
Jean~E. Roberts and Jean-Marie Thomas.
\newblock Mixed and hybrid methods.
\newblock In {\em Handbook of numerical analysis, {V}ol.\ {II}}, Handb. Numer.
  Anal., II, pages 523--639. North-Holland, Amsterdam, 1991.

\bibitem{Sandve2012}
Tor~Harald Sandve, Inga Berre, and Jan~Martin Nordbotten.
\newblock An efficient multi-point flux approximation method for {D}iscrete
  {F}racture-{M}atrix simulations.
\newblock {\em Journal of Computational Physics}, 231(9):3784--3800, 2012.

\bibitem{Schwenck2015}
Nicolas Schwenck, Bernd Flemisch, Rainer Helmig, and BarbaraI. Wohlmuth.
\newblock Dimensionally reduced flow models in fractured porous media:
  crossings and boundaries.
\newblock {\em Computational Geosciences}, 19(6):1219--1230, 2015.

\bibitem{Stefansson2018}
Ivar Stefansson, Inga Berre, and Eirik Keilegavlen.
\newblock Finite-volume discretisations for flow in fractured porous media.
\newblock {\em Transport in Porous Media}, 124(2):439--462, Sep 2018.

\bibitem{Tunc2012}
Xavier Tunc, Isabelle Faille, Thierry Gallou{\"e}t, Marie~Christine Cacas, and
  Pascal Hav{\'e}.
\newblock A model for conductive faults with non-matching grids.
\newblock {\em Computational Geosciences}, 16:277--296, 2012.

\bibitem{Noorden2009}
Tycho~L. van Noorden.
\newblock Crystal precipitation and dissolution in a porous medium: Effective
  equations and numerical experiments.
\newblock {\em Multiscale Modeling \& Simulation}, 7(3):1220--1236, 2009.

\end{thebibliography}

\end{document}